\documentclass{amsart}
\usepackage{amsmath,amssymb,amsthm,fullpage,mathptmx, hyperref,dsfont,framed, graphicx, textcomp} 

\usepackage{tikz, braids}
\usetikzlibrary{decorations.pathreplacing, arrows}
\newcommand{\cH}{\mathcal{H}}

\newcommand{\CC}{\mathcal{C}}
\newcommand{\om}{\omega}

\newcommand{\N}{\mathbb{N}}

\DeclareMathOperator{\FPdim}{FPdim}
 
\DeclareMathOperator{\End}{End}  \DeclareMathOperator{\Res}{Res}
\DeclareMathOperator{\Hom}{Hom} \DeclareMathOperator{\im}{i}
 \DeclareMathOperator{\U}{U}

\newcommand{\one}{\mathbf{1}}
\newcommand{\C}{\mathbb C}
\newcommand{\Z}{\mathbb Z}
\newcommand{\ot}{\otimes}
\newcommand{\B}{\mathcal{B}}

\usepackage[all,knot]{xy}
\xyoption{arc}
\newtheorem{example}{Example}
\newtheorem{definition}{Definition}
\usetikzlibrary{shapes,snakes,positioning,decorations.pathmorphing,calc,backgrounds,chains,fit,matrix,trees,scopes,shapes.geometric,shapes.multipart}
\tikzset{TL/.style={scale=.5, ineqn}}
\tikzset{JWP/.style={scale=.5, ineqn}}
\tikzset{TLnode/.style = {inner sep = 0, minimum size = 1cm}}
\tikzset{TL at/.style = {shift={($(#1)-(.5,.5)$)}, execute at begin scope = {\draw (0,0) rectangle (1,1);}}}
\tikzset{overstrand/.style={preaction={draw=white, -, line width=6pt}}}
\tikzset{normalizedvert/.style={circle,draw, minimum size = 1mm, inner sep = 0, fill=white}}
\tikzset{coordish/.style={minimum size=0, inner sep = 0, fill=black}}
\tikzset{nicelabel/.style={text height=1.5ex,text depth=.25ex}} 
\tikzset{horiznicelabel/.style={text width = .5em}} 
\tikzset{ineqn/.style={baseline = {($(current bounding box.center)-(0,1ex)$)}}}
\tikzset{ctrl/.style={controls = { #1 and #1} }}
\tikzset{fibabctree/.style={execute at begin picture={  
 \draw (-.6,0) -- node[auto,swap]{$a$} (0,0) coordinate(botvert) -- node[auto,swap]{$c$} (.6,0);
 \draw (botvert) -- node[auto,swap]{$b'$} (0,.5) coordinate(topvert);
 \foreach \sign in {+,-}
  \draw (topvert) -- +(\sign 1/6,.5) coordinate(leaf\sign);
}}}
\tikzset{fibacbase/.style={baseline=.5cm,nicelabel, execute at begin picture={  
 \draw (-.6,0) node[below right]{$a$} -- (-.25,0) coordinate(leftshadow) -- (0,0) coordinate(botvert) node[below]{$#1$} -- (.25,0) coordinate(rightshadow) -- (.6,0) node[below left]{$c$};
}}}

\tikzset{combspacing/.style = {row sep = .5cm, column sep = 1cm}}  
\tikzset{comb/.style = {combspacing, ampersand replacement = \&, row 1/.style=coordish, matrix of nodes, nodes in empty cells, nodes=draw}}

\usepackage[T1]{fontenc}

\newtheorem{thm}{Theorem}[section]
\newtheorem*{uthm}{Theorem}
\newtheorem{lem}[thm]{Lemma}

\newtheorem{prop}[thm]{Proposition}
\newtheorem{cor}[thm]{Corollary}
\newtheorem{conj}[thm]{Conjecture}
\newtheorem{defn}[thm]{Definition}
\newtheorem{rmk}[thm]{Remark}
\newtheorem{prob}[thm]{Open problem}
\newtheorem{ques}[thm]{Question}

\newtheorem{ex}[thm]{Exercise}

\DeclareMathAlphabet{\mathcal}{OMS}{cmsy}{m}{n}

\begin{document}

\title{Local unitary representations of the braid group and their applications to quantum computing}

\author{Colleen Delaney}
\email{cdelaney@math.ucsb.edu}
\address{Department of Mathematics,
    University of California,
    Santa Barbara, CA
    U.S.A.}

\author{Eric C. Rowell}
\email{rowell@math.tamu.edu}
\address{Department of Mathematics,
    Texas A\&M University,
    College Station, TX
    U.S.A.}

\author{Zhenghan Wang}
\email{zhenghwa@microsoft.com}
\address{Microsoft Research Station Q and Department of Mathematics,
    University of California,
    Santa Barbara, CA
    U.S.A.}

\maketitle
\begin{abstract}
We provide an elementary introduction to topological quantum computation based on the Jones representation of the braid group.  We first cover the Burau representation and Alexander polynomial.  Then we discuss the Jones representation and Jones polynomial and their application to anyonic quantum computation.  Finally we outline the approximation of the Jones polynomial and explicit localizations of braid group representations.
\end{abstract}

\section{Introduction}

Topological quantum computation (TQC) is based on the storage and manipulation of information in the representation spaces of the braid group, which consist of quantum states of certain topological phases of matter \cite{cbms}.  The most important unitary braid group representations for TQC are the Jones representations \cite{JJ}, which are part of the Temperley-Lieb-Jones (TLJ) theories.  TLJ theories are the most ubiquitous examples of unitary modular categories (UMCs).  The Jones-Wenzl projectors, or idempotents, in TLJ theories can be used to model anyons, which are believed to exist in fractional quantum Hall liquids.  Hence the proper mathematical language to discuss TQC is UMC theory and the associated topological quantum field theory (TQFT).  Both UMC and TQFT are highly technical subjects.  However, the representations of the braid group from UMCs or TQFTs are a more accessible point of entry to the subject.  These notes provide an elementary introduction to some representations of the braid group coming from 
UMCs and TQFTs, and their application to topological quantum computation.  We will use {\it the braid group} $\B_\infty$ to mean the direct limit of all $n$-strand braid groups $\B_n$ for all $n\geq 1$.  Therefore, a representation of the braid group $\B_\infty$ is a compatible sequence of representations of $\B_n$.

Our focus is on the representations of the braid group discovered by Jones in the study of von Neumann algebras \cite{JJ}.  Jones representations are \emph{unitary}, which is important for our application to quantum computing.  These representations also have a hidden \emph{locality} and generically dense images.  Unitarity, locality, and density are important ingredients for the two main theorems (the colloquial terms will be made precise later) that we will present:

\begin{thm}\label{FLWthm} The Jones representation of the braid group at $q=e^{\pm 2\pi i /r}$, $r \ne 1,2,3,4,6$, can be used to construct a universal quantum computer.
\end{thm}

\begin{thm}\label{FKWthm} The Jones polynomial of oriented links at $q=e^{\pm 2\pi i/r}$ can be approximated by a quantum computer efficiently for any integer $r\geq 1$.
\end{thm}

While unitarity and density are easy to understand mathematically, locality is not formally defined in our notes as there are several interpretations (one of which is discussed in Section \ref{localization}).  Essentially, a local representation of the braid group is one coming from a local TQFT, whose locality is encoded in the gluing formula.  A first approximation of locality would mean a sequence of representations of $\B_n$ with a compatible Bratteli diagram of branching rules.

We motivate our study of the Jones representation and its quantum applications with the Burau representation, which belongs to the classical world. The Burau representation leads to the link invariant called the Alexander polynomial, which can be computed in polynomial time on a classical computer. On the other hand, the link invariant corresponding to the Jones representation, the Jones polynomial, is $\#P$-hard to compute on a classical computer, but can be approximated by a quantum computer in polynomial time. This approximation of quantum invariants by a quantum computer is realized by the amplitudes of the physical processes of quasi-particles called anyons, whose worldlines include braids.

The contents of these notes are as follows: In section 2, we cover the Burau representation and Alexander polynomial.  In section 3, we discuss the Jones representation and Jones polynomial.  Section 4 is on anyons and anyonic quantum computation.  In section 5, we explain the approximation of the Jones polynomial.  Section 6 is on an explicit localization of braid group representations.  While  full details are not included, our presentation is more or less complete with the exception of Thm. \ref{2qubits}, which is important for addressing the leakage issue.  An elementary inductive argument for Thm. \ref{2qubits} is possible and we will leave it to interested readers.

\section{Burau Representation and Alexander Polynomial}

\subsection{The braid group}

The \emph{n-strand braid group} $\B_n$ is given by the presentation $$\B_n = \Big \langle \sigma_1, \sigma_2, \ldots, \sigma_{n-1} \Bigm\lvert \sigma_i\sigma_j = \sigma_j \sigma_i \text{ for } |i -j| \ge 2, \sigma_i\sigma_{i+1}\sigma_i = \sigma_{i+1}\sigma_i\sigma_{i+1}, i=1, 2, \ldots, n-1 \Big \rangle.$$

The first type of relation is known as far commutativity and the second is the braid relation.  Using the braid relation, one can check that all of the generators of the $n$-strand braid group lie in the same conjugacy class.  Therefore, each $n$-strand braid group $\B_n$ is generated by a single conjugacy class when $n\geq 3$.

The names of the relations are inspired by the geometric presentation of the braid group, in which we picture braids on $n$ \lq\lq strands'', and the braid generators $\sigma_i$ correspond to crossing the $i$th strand over the $i+1$ strand. Multiplication $bb'$ of two braid diagrams $b$ and $b'$ is performed by stacking $b'$ on top of $b$ and interpreting the result as a new braid diagram.

For example, $\B_3 = \langle \sigma_1, \sigma_2 \mid \sigma_1\sigma_2\sigma_1 = \sigma_2\sigma_1\sigma_2 \rangle$, where $\sigma_1$ braids the first two strands and $\sigma_2$ the latter two. $$\begin{tikzpicture}[scale=.5]
\node (sigma1) {};
\draw (0,-1) node {$\sigma_1=$};
 \braid[number of strands={3}] a_1^{-1}  ;
\node (sigma2) [right=of sigma1, xshift=15mm] {};
\begin{scope}[shift={(sigma2)}]
\draw (0,-1) node {$\sigma_2=$};
 \braid a_2^{-1} ;
 \end{scope}
 \end{tikzpicture}$$

In these notes, we use the \lq\lq right-handed convention'' when drawing braid diagrams, so that the overstrand goes from bottom left to top right.
As a result, $$\begin{tikzpicture}[scale=.5]
\draw (-.5,-1) node {$\sigma_1^{-1}=$};
 \braid[number of strands={3}] a_1  ;
 \end{tikzpicture}$$
Swapping the definitions of $\sigma_1$ and $\sigma_1^{-1}$ would give the \lq\lq left-handed convention".

In the picture presentation, far commutativity expresses the fact that when nonoverlapping sets of strands are braided, the result is independent of the order in which the strands were braided. The braid relation is given by $$\begin{tikzpicture}[scale=.5]
\node (braid1) {};
 \braid[number of strands={3}] a_1^{-1} a_2^{-1} a_1^{-1} ;
\node (braid2) [right=of braid1, xshift=10mm] {};
\begin{scope}[shift={(braid2)}]
 \draw(-.25,-2) node {\Large $=$};
 \braid[number of strands={3}] a_2^{-1} a_1^{-1} a_2^{-1} ;
 \end{scope}
 \end{tikzpicture}.$$

The braid relation is called the Yang-Baxter equation by some authors, but we will reserve use of this phrase because, as will be explained shortly, there is a subtle difference between the two.

Another useful perspective is to identify $\B_n$ with the motion group (fundamental group of configuration space) of $n$ points in the disk $D^2$.  Then the braid relation says: given three distinct points on a line in the disk and we exchange the first and third while keeping the middle one stationary, then the braid trajectories are the same if the first and third points cross to the left or right of the middle point.

The \emph{braid group}, denoted by $\B_{\infty}$, is formed by taking the direct limit of the $n$-strand braid groups with respect to the inclusion maps $\B_n \hookrightarrow \B_{n+1}$ sending $\sigma_i \to \sigma_i$. That is, we identify a braid word in $\B_n$ with the same braid word in $\B_{n+1}$.  In pictures, this inclusion map $\B_{n} \to \B_{n+1}$ adds a single strand after the braid $\sigma$.

\subsection{Representations of the Braid Group}
For applications of braid group representations to quantum computing, the braid group representations should be \emph{unitary} and \emph{local}. Moreover, for reasons that are not {\it a priori} clear, since the images of the braid generators $\sigma_i$ will eventually be interpreted as { \it quantum gates} manipulating {\it quantum bits}, they should be of finite order and have algebraic matrix entries.

Recall that a matrix $U$ is unitary if $U^{\dagger}U=UU^{\dagger}=I.$ We denote by $\U(r)$ the group of $r \times r$ unitary matrices. A precise definition of locality requires interpreting the images of elements of the braid group as quantum gates, and is relegated to section 4 where quantum computation is discussed.

One important way to obtain representations of the braid group is to find solutions to the {\it Yang-Baxter equation}.

\subsubsection{The Yang-Baxter Equation and $R$-matrix}
 Let $V$ be a finite dimensional complex vector space with a specified basis, and let $R: V \otimes V \to V \otimes V$ be an invertible solution to the \emph{Yang-Baxter equation} (YBE):

 $$(R \otimes I)(I \otimes R)(R\otimes I) =(I\otimes R)(R \otimes I)(I \otimes R)$$
 where $I$ is the identity transformation of $V$. We call such a solution to the YBE an \emph{R-matrix} (as opposed to R-operator, since we have a basis with which to work).

Any $R$-matrix gives rise to a (local) representation of the braid group via the identification

$$\begin{tikzpicture}
\draw (-1/4,1/2) node {$\sigma_i \rightarrow$};
\draw (1/4,0)--(1/4,1) ;
\draw (3/8,1/2) node {\tiny $\cdots$};
\draw (1/2,1/4) rectangle (1,3/4) ;
\draw (3/4,1/2) node {R};
\draw (5/8,0)--(5/8,1/4) (5/8,3/4)--(5/8,1);
\draw (7/8,0)--(7/8,1/4) (7/8,3/4)--(7/8,1);
\draw (9/8,1/2) node {\tiny $\cdots$};
\draw (10/8,0)--(10/8,1);
\end{tikzpicture}$$

 For example, in the 3-strand braid group, we can take $$\begin{tikzpicture}
 \draw (-1/4,1/2) node {$\sigma_1 \rightarrow$};
\draw (1/2,1/4) rectangle (1,3/4) ;
\draw (3/4,1/2) node {R};
\draw (5/8,0)--(5/8,1/4) (5/8,3/4)--(5/8,1);
\draw (7/8,0)--(7/8,1/4) (7/8,3/4)--(7/8,1);
\draw (10/8,0)--(10/8,1);
\draw (2, 1/2) node {$=R \otimes id_V$};
\end{tikzpicture}$$
 where $R\otimes id_V$ is a map from $V \otimes V \otimes V$ to itself.

Representations of the braid group arising from $R$-matrices are always local, but rarely unitary. There is a natural tension between these two properties that make finding such a representation difficult, which is illustrated by the following example.
\subsubsection{Locality versus unitarity} The following $R$-matrix is a $4 \times 4$ solution to the Yang-Baxter equation for $V=\C^2$ with the standard basis, and as such is local. $$R = \begin{pmatrix} a & 0 & 0 & 0 \\ 0 & 0 & \bar{a} & 0 \\ 0 & \bar{a} & a - \bar{a}^3 & 0 \\ 0 & 0 & 0 & a \\ \end{pmatrix}$$

However, if $R$ were to be unitary, its columns would have to be orthonormal. In particular, $\|a\|=1$ and $$\langle (0,0,\bar{a}, 0)^T, (0, \bar{a},a-\bar{a}^3, 0 )  \rangle = a(a-\bar{a}^3)=0.$$
Since $\bar{a}=a^{-1}$, unitarity would force $a^4=1$. But then the only possibilities for $a$ are $\pm 1$ or  $\pm i$, and one can check that indeed each of these choices results in $R$ being a unitary matrix. \\

 It is in general difficult to find nontrivial solutions to the Yang-Baxter equation. Historically, the theory of quantum groups was developed to address this problem, but solutions that arise from the theory of quantum groups are rarely unitary. The state of the art is that for $\dim V=1$ and $\dim V=2$, all unitary solutions are known. While a classification for larger dimensions is unknown, there do exist nice examples of $4 \times 4$ and $9 \times 9$ unitary solutions \cite{cbms}.





\subsection{The Burau representation of the braid group}
There are two versions of the Burau representation: the unreduced representation, which denoted by $\tilde{\rho}$, and the reduced representation, for which we reserve the notation $\rho$.

\subsubsection{The unreduced Burau representation}

There is a nice probabilistic interpretation of the unreduced Burau representation that is due to Jones, which we will use as an introduction to the subject \cite{jones}. We start by defining the representation for \emph{positive braids}, braids for which all crossings are right-handed. More precisely, $\sigma$ is positive if it can be written $\sigma=\sigma_{i_k}^{s_k} \cdots \sigma_{i_1}^{s_i}$ where $s_i > 0$ for each $i$.

Imagine the braid diagram as a braided bowling alley with $n$ lanes, where lanes cross over and under one another and at every overcrossing there is a trap door that will open with probability $1-t$ when a ball rolls over it. Of course, due to gravity there is zero probability of a ball on a lower lane jumping up onto a lane crossing over it. Then starting from the bottom of the braid and bowling down lane $i$, it ends up in lane $j$ with some probability, which we can identify as the $ij$th entry of a matrix.

More precisely, for positive braids the unreduced Burau representation $\tilde{\rho}: \B_n \to GL_n(\mathds{Z}[t,t^{-1}])$ can be defined by assigning each $\sigma \in B_n$ to the matrix $\tilde{\rho}(\sigma)$ given by $$\tilde{\rho}(\sigma)_{ij} = \sum_{\text{paths p from i to j}} w(p),$$

where $w(p)$ is the probability corresponding to the path $p$, which is always of the form $t^k(1-t)^l$ for some nonnegative integers $k$ and $l$.

For a concrete example, take the following braid, call it $\sigma$, in $\B_4$.$$
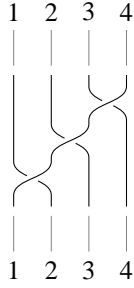
\begin{tikzpicture}[scale=.5]
 \braid[number of strands={4}] (braid)  a_3^{-1} a_2^{-1} a_1^{-1} ;
\node[at=(braid-rev-1-e), pin=south :1] {};
\node[at=(braid-rev-2-e), pin=south :2] {};
\node[at=(braid-rev-3-e), pin=south :3] {};
\node[at=(braid-rev-4-e), pin=south :4] {};
\node[at=(braid-1-s), pin=north :1] {};
\node[at=(braid-2-s), pin=north :2] {};
\node[at=(braid-3-s), pin=north :3] {};
\node[at=(braid-4-s), pin=north :4] {};
 \end{tikzpicture}$$

Note that the labels mark the relative position of the strands, as opposed to the strands themselves. The matrix representing $\sigma$ in $GL_4(\mathds{Z}[t,t^{-1}])$ is then given by

$$\tilde{\rho}(\sigma)=\begin{pmatrix} 1-t & t(1-t) & t^2(1-t) & t^3 \\ 1 & 0 & 0 & 0 \\  0 &  1 & 0 & 0 \\ 0 & 0 & 1 & 0 \\ \end{pmatrix}$$

How to represent left-handed crossings is then already determined, for once we define the representation of a right handed crossing $$\begin{tikzpicture}[scale=.5]
\draw (1/2,-3/4) node {$\tilde{\rho} ($};
\draw (2.25, -3/4) node {$)$};
 \braid[number of strands={2}] a_1^{-1}  ;
 \draw (4.5,-3/4) node {$=\begin{pmatrix} 1-t & t \\ 1 & 0 \end{pmatrix}$};
 \end{tikzpicture},$$

using that $\sigma_1\sigma_1^{-1}=1$ and that $\tilde{\rho}$ is a group homomorphism it follows that $\tilde{\rho}(\sigma_1)\tilde{\rho}(\sigma_1^{-1})=I.$ Therefore the representation of the $\sigma_1^{-1}$ must be given by the inverse of $\tilde{\rho}(\sigma_1),$

 $$\begin{tikzpicture}[scale=.5]
\draw (1/2,-3/4) node {$\tilde{\rho} ($};
\draw (2.25, -3/4) node {$)$};
 \braid[number of strands={2}] a_1  ;
 \draw (5.5,-3/4) node {$=\begin{pmatrix} 0 & 1 \\ \bar{t} & 1-\bar{t}  \\ \end{pmatrix}$};
 \end{tikzpicture},$$where $\bar{t}=1/t$.
Thus left-handed crossings are assigned a factor of $\bar{t}$ for an overcrossing and $1-\bar{t}$ for an undercrossing. The generators $\sigma_i$ of $\B_n$ and their inverses can be represented by extending the construction in the natural way. Then the representation of an arbitrary braid $b=\sigma_{i_k}^{s_k}\cdots \sigma_{i_1}^{s_1}$ is given by multiplying the representations of the constituent $\sigma_{i_k}$ in the braid word. This defines the unreduced Burau representation of the braid group.

As a fun example we introduce the following braid $b=\sigma_3^{-1}\sigma_2^2 \sigma_{3}^{-1}\sigma_1^{-1}$, once drawn by Gauss (see e.g. \cite[Figure 2]{Epple}).
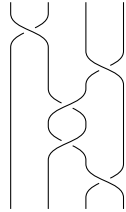
$$\begin{tikzpicture}[scale=.5]
 \braid[number of strands=4] (braid)  a_1 a_3 a_2^{-1} a_2^{-1}  a_3 ;

 \end{tikzpicture} $$

The unreduced Burau representation of the Gauss braid is given by

$$\begin{pmatrix} 0 & 1 & 0 & 0 \\ t\bar{t} + (1-t)^2\bar{t} & t(1-\bar{t}) + (1-t)^2(1-\bar{t}) &0 & (1-t)t  \\ 0 & 0 & \bar{t} & (1-\bar{t}) \\ \bar{t}^2(1-t) & \bar{t}(1-t)(1-\bar{t}) & (1-\bar{t})\bar{t} & (1-\bar{t})^2 +\bar{t}t\end{pmatrix},$$
which has been left unsimplified to make the individual contributions from paths more transparent.

The unreduced Burau representation of a braid $b \in \B_n$ has several properties worth mentioning.

\begin{enumerate}
\item When $t=1$, $\tilde{\rho}(b)$ is a permutation matrix. This allows one to interpret $\tilde{\rho}(b)$ as a deformation of a permutation matrix.
\item The representation $\tilde{\rho}$ is reducible.
\item There exists an invariant row vector (row eigenvector) of $\tilde{\rho}(b)$, independent of $b \in \B_n$.
\end{enumerate}
The first property is clear from the construction of the unreduced Burau representation. The second and third properties are closed related, and we prove them below.

One of the nice aspects of the probabilistic interpretation of the Burau representation is that it is an immediate consequence of the definition that the entries in each row of a matrix $\tilde{\rho}(b)$ should sum to one, since probability must be conserved. Put another way,
$$\tilde{\rho}(\sigma) \begin{pmatrix} 1 \\ 1 \\ \vdots \\ 1\\ \end{pmatrix} = \begin{pmatrix} 1 \\ 1 \\ \vdots \\ 1\\ \end{pmatrix}.$$ That is, there is a one-dimensional subspace which is invariant under $\tilde{\rho}(\sigma)$, for any $\sigma \in \B_n$. Therefore the (unreduced) Burau representation is reducible, and that we can obtain another representation by restricting to the orthogonal subspace span$\{ (1,1,\ldots,1)\}^{\perp}$. This is one way to define the \emph{reduced} Burau representation.

Since the determinant of a matrix is equal to the determinant of its transpose, if $\det (I - \tilde{\rho}(b))=0$ for all $b \in \B_n$, then

$$\det((I-\tilde{\rho}(b))^{T})= \det(I - \tilde{\rho}(b)^{T})=0$$
for all $b \in \B_n$. In other words, since $\tilde{\rho}(b)$ has an eigenvector $\tilde{\rho}(b)^{T}$ has an eigenvector $v$ with eigenvalue $1$, $\tilde{\rho}(b)^{T}v=v$ for some $v\ne 0$. Taking the transpose of this matrix equation, we find

$$v^T \tilde{\rho}(b) = v^{T}.$$

This shows that $\tilde{\rho}(b)$ has an invariant row vector $v^{T}$, proving the third property.

In fact, this row vector takes the form

$$v^{T}=(1,t,t^2, \ldots, t^{n-1}).$$

Observing that
$$(1, \ldots, 1, t^{i}, t^{i+1}, 1 \ldots, 1) \begin{pmatrix}I_{i-1} & 0 & 0 \\ 0 & \begin{pmatrix} 1-t & t \\ 1 & 0 \\ \end{pmatrix} & 0 \\ 0 & 0 & I_{n-i-1}
\end{pmatrix}= (1,\ldots, 1, t^{i} - t^{i+1} + t^{i+1}, t^{i+1}, 1, \ldots, 1) $$
it follows that $v^T$ defines an invariant row vector for the representations of the braid group generators $\tilde{\rho}(\sigma_i)$, and hence for all $\tilde{\rho}(b)$, $b \in \B_n$.

\subsubsection{The reduced Burau representation}

An alternative approach to defining the reduced Burau representation yields an explicit basis.

We explicitly find a basis for an invariant subspace of $\tilde{\rho}(\B_n)$ by looking for eigenvalues and eigenvectors of $\tilde{\rho}(\sigma_i)$. We have seen already that $(1,1,\ldots, 1)^T$ is an eigenvector with eigenvalue 1. One can check that another eigenvector is given by $(0, \ldots, 0, \underbrace{-t}_{i}, \underbrace{1}_{i+1}, 0 \ldots, 0 )^{T}$, corresponding to eigenvalue $-t$.

\begin{prop} Let $v_i = (0, \ldots, 0, \underbrace{-t}_{i}, \underbrace{1}_{i+1}, 0 \ldots, 0 )^{T}$. Then span$\{v_1, \ldots, v_{n-1}\}$ is an invariant subspace of $\tilde{\rho}(b)$ for all $b \in \B_n$.

\end{prop}

Proof. The image of each $v_i$ under $\tilde{\rho}(b)$ can be written as a linear combination of the $v_j$.
$$\tilde{\rho}(\sigma_i)v_i= \begin{pmatrix}I_{i-1} & 0 & 0 \\ 0 & \begin{pmatrix} 1-t & t \\ 1 & 0 \\ \end{pmatrix} & 0 \\ 0 & 0 & I_{n-i-1}
\end{pmatrix} =\begin{pmatrix} 0 \\ \vdots \\ 0  \\t^2 \\ -t \\ 0 \\ \vdots \\ 0\end{pmatrix}=  -tv_i$$

Similar calculations show
$$\tilde{\rho}(\sigma_i)v_{i-1} =v_{i-1}+v_i$$
$$\tilde{\rho}(\sigma_i)v_{i+1}= -tv_i + v_{i+1}$$

and
$$\tilde{\rho}(\sigma_i)v_j = v_j$$ when $|j-i|\ge 2$.

This verifies that the subspace spanned by the $v_i$ is invariant, leading to the following definition.

\begin{defn} The reduced Burau representation $\rho:\B_n \to GL_{n-1}(\mathds{Z}[t,t^{-1}])$ is given by $\rho(b)=\tilde{\rho}(b) \Bigm\lvert_{\text{span}\{v_i\}}$.
\end{defn}

\subsubsection{Unitary Burau representations}

Keeping in mind that we are looking for unitary representations of the braid group, it is natural to ask for which $t \in \mathds{C}^{*}$ the reduced Burau representation $\rho:\B_n \to GL_{n-1}(\mathds{Z}[t,t^{-1}])$ is unitary. \\

One can check that the matrix representation corresponding to a braid group generator fails to be unitary for any choice of $t$.
For example, consider the braid \hspace{1mm}$\begin{tikzpicture}[baseline=-.4cm, scale=.4]
 \braid[number of strands=2] (braid)   a_1^{-1} ;
 \end{tikzpicture} $ \hspace{.05mm} with corresponding (unreduced) Burau matrix representation $\begin{pmatrix} 1 - t & t \\ 1 & 0 \\ \end{pmatrix}$. Then either by direct computation or by noting that (since we can safely ignore an overall phase factor of -1 without affecting unitarity) for no choice of $t$ does this matrix admit the familiar parametrization of elements of SU(2) as $\begin{pmatrix} a & b \\ -\bar{b} & \bar{a} \\ \end{pmatrix}$ where $a,b \in \mathds{C}$ and  $|a|^2 + |b|^2 = 1.$ It follows that the unreduced Burau representation is never unitary.


However, the situation is not completely hopeless. We end with a theorem that tells us how to obtain a unitary representation from the reduced Burau representation.

\begin{thm} Let $t=s^2$ where $s \in \mathds{C}^{*}$, and define $P_{n-1} = \begin{pmatrix} 1 & 0 & \cdots & 0 \\
0 & s &  & \\
\vdots & & \ddots &  \vdots& \\
0 & \cdots & & s^{n-1} \\
\end{pmatrix}$, $J_{n-1} = \begin{pmatrix} s+s^{-1} & -1 &\cdots & 0 \\
-1 & s+s^{-1} & \ddots  & \vdots\\
\vdots & \ddots & \ddots &  -1 & \\
0 & \cdots & -1 & s+s^{-1} \\
\end{pmatrix}$, and $\rho_s(b)=P_{n-1}\rho(b)(P_{n-1})^{-1}$, where $\rho$ is the reduced Burau representation and $b \in \B_n$.

Then $\rho_s$ is unitary with respect to the Hermitian matrix $J_{n-1}$. That is,

$$(\rho_s(b))^{\dagger} J_{n-1}\rho_s(b) = J_{n-1}$$

\end{thm}

Moreover, for those $s \in \mathds{C}^{*}$ for which $J_{n-1}(s)$ can be written as $J_{n-1}(s) = X^{\dagger}X$ for some matrix $X$, $X\rho_s(b)X^{-1}$ gives a unitary representation.

\begin{ex} Find all $s$ such that $J_{n-1}(s)$ can be written $J_{n-1}(s)=X^{\dagger}X$.
\end{ex}

There remain basic questions about the Burau representation to which the answers are not yet known.
\begin{prob} The Burau representation is faithful for $n=1,2,3$, and is not faithful for $n\ge 5$ \cite{bigelow}. What about when $n=4$?
\end{prob}


\subsection{The Alexander polynomial}
The reduced Burau representation of the braid group can be used to study links by using it to define an invariant called the \emph{Alexander polynomial}. The existence of invariants which are both powerful and computable is essential to the classification of any mathematical object. Of course, Nature conspires so that these two characteristics are often hard to satisfy simultaneously. While the Alexander polynomial will be computable in polynomial time, we will see that it is not quite sensitive enough to distinguish between certain types of knots.

\subsubsection{From braids to links}

There is a natural way to turn a braid into a link by identifying the top and bottom strands in an order-preserving manner. This operation is called a \emph{braid closure}.

For example, the closure of the Gauss braid is the connect sum of two Hopf links.  $$
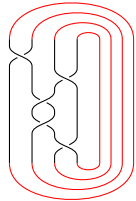
\begin{tikzpicture}[scale=.30]
\node (gauss braid) {};
 \braid[number of strands=4] (braid)  a_1 a_3 a_2^{-1} a_2^{-1}  a_3 ;
\draw[red] (4,0) to [out=90, in=90] (5,0)--(5,-5.5) to [out=-90, in=-90] (4,-5.5);
\draw[red] (3,0)  to [out=90, in=90] (5.5,0)--(5.5,-5.5) to [out=-90, in=-90] (3,-5.5);
\draw[red] (2,0)  to [out=90, in=90] (6,0)--(6,-5.5) to [out=-90, in=-90] (2,-5.5);
\draw[red] (1,0)  to [out=90, in=90] (6.5,0)--(6.5,-5.5) to [out=-90, in=-90] (1,-5.5);
\end{tikzpicture}$$

\subsubsection{The Markov moves}
Two other operations on braids, conjugation and stabilization, also known as the Markov moves of type I and type II, respectively. We will see that performing either type of operation on a braid does not change the link that is obtained from the braid closure.

\textbf{I.} Let $b, g \in \B_n$. Then conjugation of the braid $b$ by the braid $g$ is given by

$$b \mapsto gbg^{-1}.$$

The best way to see that $\hat{b}=\widehat{gbg^{-1}}$ is through a picture:
$$
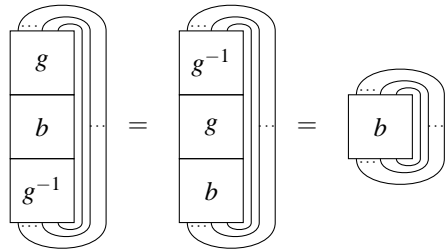
\begin{tikzpicture}[scale=.85]
\node (diagram) {};
\draw (0,0) rectangle (1,1);
\draw (1/2,1/2) node {$b$};
\draw (0,1) rectangle (1,2);
\draw (1/2,3/2) node {$g$};
\draw (1/2,-1/2) node {$g^{-1}$};
\draw (0,0) rectangle (1,-1);
\draw  (3/4,2) to [out=90, in=90] (9/8,2)--(9/8,-1) to [out=-90,in=-90] (3/4,-1);
\draw  (1/2,2) to [out=90, in=90] (10/8,2)--(10/8,-1) to [out=-90,in=-90] (1/2,-1);
\draw (5/16,-17/16) node {\tiny $\cdots$};
\draw (5/16,33/16) node {\tiny $\cdots$};
\draw (1/8,2)  to [out=90, in=90] (12/8,2)--(12/8,-1) to [out=-90,in=-90] (1/8,-1);
\draw (11/8,1/2) node {\tiny $\cdots$};
\draw (2, 1/2) node {$=$};
\node (diagram2) [right=of diagram, xshift=10mm] {};

\begin{scope}[shift={(diagram2)}]
\draw (0,0) rectangle (1,1);
\draw (1/2,1/2) node {$g$};
\draw (0,1) rectangle (1,2);
\draw (1/2,3/2) node {$g^{-1}$};
\draw (1/2,-1/2) node {$b$};
\draw (0,0) rectangle (1,-1);

\draw  (3/4,2) to [out=90, in=90] (9/8,2)--(9/8,-1) to [out=-90,in=-90] (3/4,-1);
\draw  (1/2,2) to [out=90, in=90] (10/8,2)--(10/8,-1) to [out=-90,in=-90] (1/2,-1);
\draw (5/16,-17/16) node {\tiny $\cdots$};
\draw (5/16,33/16) node {\tiny $\cdots$};
\draw (1/8,2)  to [out=90, in=90] (12/8,2)--(12/8,-1) to [out=-90,in=-90] (1/8,-1);
\draw (11/8,1/2) node {\tiny $\cdots$};
\draw (2, 1/2) node {$=$};
\end{scope}
\node (diagram3) [right=of diagram2, xshift=10mm] {};
\begin{scope}[shift={(diagram3)}]
\draw (0,0) rectangle (1,1);
\draw (1/2,1/2) node {$b$};
\draw  (3/4,1) to [out=90, in=90] (9/8,1)--(9/8,0) to [out=-90,in=-90] (3/4,0);
\draw  (1/2,1) to [out=90, in=90] (10/8,1)--(10/8,0) to [out=-90,in=-90] (1/2,0);
\draw (5/16,-1/16) node {\tiny $\cdots$};
\draw (5/16,17/16) node {\tiny $\cdots$};
\draw (1/8,1)  to [out=90, in=90] (12/8,1)--(12/8,0) to [out=-90,in=-90] (1/8,0);
\draw (11/8,1/2) node {\tiny $\cdots$};
\end{scope}
\end{tikzpicture}$$

\textbf{II.} Let $b \in \B_n$, and let $\B_n$ be embedded in $\B_{n+1}$ in the usual way, by adding a rightmost strand. Then stabilization of the braid $b$ is given by the map $\B_n \to \B_{n+1}$

$$b \mapsto b\sigma_n^{\pm 1}.$$

That is, we add a rightmost $n+1$st strand to $b$ to identify it as a braid in $\B_{n+1}$, and then we braid its $n$th and $n+1$st strands with either an over- or under-crossing. Once again, a picture best demonstrates that $\hat{b}=\widehat{b\sigma_n}$.

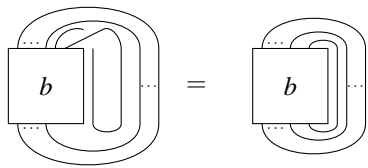
\begin{figure}[h!]
\centering
\begin{tikzpicture}[scale=1]
\node(diagram1) {};
\draw (0,0) rectangle (1,1);
\draw (1/2,1/2) node {$b$};

\draw  (1/2,1) to [out=90, in=90] (14/8,1)--(14/8,0) to [out=-90,in=-90] (1/2,0);
\draw (3/4,1) -- (10/8,1.25) to [out=30, in=90] (12/8,1)--(12/8,0) to [out=-90, in=-90] (9/8,0);
\draw (9/8,1) --(9/8,0);
\draw (5/16,-1/16) node {\tiny $\cdots$};

\draw (5/8,1) to [out=90, in =160] (1,10/8);
\draw (5/16,17/16) node {\tiny $\cdots$};

\draw (1/8,1)  to [out=90, in=90] (16/8,1)--(16/8,0) to [out=-90,in=-90] (1/8,0);
\draw (15/8,1/2) node {\tiny $\cdots$};
\node (diagram2) [right=of diagram1,xshift=20mm] {};
\begin{scope}[shift={(diagram2)}]
\draw (-3/4,1/2) node {$=$};
\draw (0,0) rectangle (1,1);
\draw (1/2,1/2) node {$b$};
\draw  (3/4,1) to [out=90, in=90] (9/8,1)--(9/8,0) to [out=-90,in=-90] (3/4,0);
\draw  (1/2,1) to [out=90, in=90] (10/8,1)--(10/8,0) to [out=-90,in=-90] (1/2,0);
\draw (5/16,-1/16) node {\tiny $\cdots$};
\draw (5/16,17/16) node {\tiny $\cdots$};
\draw (1/8,1)  to [out=90, in=90] (12/8,1)--(12/8,0) to [out=-90,in=-90] (1/8,0);
\draw (11/8,1/2) node {\tiny $\cdots$};
\end{scope}
\end{tikzpicture}
\end{figure} A similar picture shows $\hat{b}=\widehat{b\sigma_n^{-1}}$.

This move introduces a twist in the braid closure, and hence can be undone by a Reidemeister move of type 1, so it doesn't change the link $\hat{b}$.

Not only does manipulating a braid by Markov moves not change the braid closure, but whenever two braid closures agree, their corresponding braids can be related by a finite number of Markov moves.

\begin{thm}[Markov] Consider the map from the set of all braids to the set of all links given by
$$\{\B_n\} \to \{\text{links}\}$$
$$b \mapsto \hat{b}$$

If $\widehat{b_1}=\widehat{b_2}$ as links, then $b_1$ and $b_2$ are related by a finite number of moves of type I or type II and their inverses.
\end{thm}

It is easy to see that the map $b \to \hat{b}$ fails to be injective. For the simplest possible example, take the braid closure of $\sigma_1$, which gives the unknot.


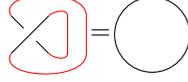
\begin{figure}[h!]
\centering
\begin{tikzpicture}[scale=.5]
\node (closure) {};
\draw (0,0) --( 1,1);
\draw (0,1)--(1/3,2/3);
\draw (2/3,1/3)--(1,0);
\draw[red] (1,0) to [out=-45, in=-90] (3/2,0)--(3/2,1) to [out=90, in=45] (1,1);
\draw[red] (0,0) to [out=-120, in=-90] (2,0)--(2,1) to [out=90, in=120] (0,1);
\node (eq) [right= of closure, xshift = -2mm, yshift =2.5mm] {$=$};
\draw (15/4,1/2) circle (1);
\end{tikzpicture}
\end{figure}

It is also true, although much less trivial to show, that the map is onto. Given any link there exists a finite number of Reidemeister moves that manipulates the link until it is in the form of a closure of a braid.

The Markov theorem gives us a way to study links through braid representations, since any braid invariant that is also invariant under the Markov moves can be improved to a link invariant. 

 \subsubsection{The Alexander polynomial}

In order to be an invariant of links, a quantity must be invariant under the Markov moves of type I and II. From linear algebra, we know that similar matrices have the same determinant. It follows that the determinant of the representation of a braid is invariant under conjugation.

Recall the reduced Burau representation $\rho: \B_n \to GL_{n-1}(\mathds{Z}[t,t^{-1}])$ and define the matrices $M(b)=I - \rho(b)$ and $\tilde{M}(b)=I-\tilde{\rho}(b)$, where $I$ is the identity matrix with appropriate dimensions in each equation.

\begin{defn}
For $b\in \B_n$, the Alexander polynomial is given by

$$\Delta(\hat{b},t) = \frac{\det(M(b))}{1+t+\cdots t^{n-1}}.$$
\end{defn}

This establishes the convention that the Alexander polynomial of the unknot is 1, i.e.
\begin{tikzpicture}[baseline=-2] \node (un) [label=left: $\Delta \Big ($, label=right: $\Big )$] {};
\draw (0,0) circle (1/4);
\node (eq) [right=of un, xshift=-7mm] {$= 1$};
\end{tikzpicture}. We present some results from linear algebra that imply a proof that the Alexander polynomial is a link invariant, and state some of its properties.  The proof that the defined polynomial is indeed a link invariant is simply a check of the invariance of  $\Delta(\hat{b},t)$ under Markov moves using the lemmas below.

\begin{thm}

The Alexander polynomial of a link can be computed in polynomial time by a Turing machine.

\end{thm}

There is a polynomial time algorithm to turn any link into a braid closure.  For a braid $b\in \B_n$, its Burau representation matrix can be computed in poly(n,m) and so can the determinant, where $m$ is the number of elementary braids in $b$.  Note that the size of the Burau representation matrix is only $(n-1)\times (n-1)$ for a braid in $\B_n$.  As a comparison, we will see later the sizes of the Jones representation matrices for braids in $\B_n$ grow as $d^n\times d^n$ for some number $d>1$ as $n\rightarrow \infty$.

\subsubsection{Linear algebraic interlude}

 In order to prove some properties of the Alexander polynomial, it will be useful to apply the following lemma from linear algebra to our present setting.

\begin{lem} Suppose $A$ is an $n \times n$ matrix with the property that there exists a column vector $w=(w_i)^T$ and a row vector $u=(u_j)$ satisfying

\begin{enumerate}
\item $Aw = 0$
\item $uA=0$
\item $w_i \ne 0$, $u_j \ne 0$ for all $i,j$.
\end{enumerate}
That is, $A$ annihilates $w$, $A$ is annihilated by $u$, and the coordinates of $w$ and $u$ are all nonvanishing.
Then

$$(-1)^{i+j} \frac{\det(A(i,j))}{u_iw_j}$$ is independent of $i$ and $j$, where $A(i,j)$ denotes the $i,j$th minor of $A$ - the $(n-1) \times (n-1)$ matrix obtained by the deleting the $i$th row and the $j$th column from $A$.

\end{lem}

The matrix $\tilde{M}$ satisfies the hypotheses of this lemma. Recall that $\tilde{\rho}(b)$ had eigenvector $(1, \ldots, 1)^T$ with eigenvalue $1$ and invariant row vector $(1, t, t^2, \ldots, t^{n-1})$. If we choose $w = (1, \ldots, 1)^T$ and $u=(1, t, t^2, \ldots, t^{n-1})$, it follows that $\tilde{M}w = 0$ and $u\tilde{M}=0$. Evidently the coordinates of both $w$ and $u$ are nonvanishing. While the details are omitted, this leads to the proof of the next lemma.

\begin{lem}\label{detlemma} $$\frac{\det(M(b))}{1+t+\cdots + t^{n-1}}=\det(\tilde{M}(1,1)).$$
\end{lem}

This result gives us the freedom to delete any row and column of the matrix $\tilde{M}$, whose determinant recovers the Alexander polynomial. This makes certain computations easier, like the proof of the skein relation.

\subsection{The Alexander-Conway polynomial, writhe, and skein relation}  Having introduced the Alexander polynomial, one can define a related link invariant - the \emph{Alexander-Conway polynomial} - through a slight renormalization and a quantity called the \emph{writhe} of a braid.

Let $b=\sigma_{i_k}^{s_k}\cdots \sigma_{i_1}^{s_1} \in \B_n$. The writhe or braid exponent is given by $e(b)= \sum_{i=1}^k s_i$.

Taking $z=t^{1/2}-t^{-1/2}$, the Alexander-Conway polynomial is defined as

$$\Delta(\hat{b}, z) = (-t^{1/2})^{n-e(b)-1}\Delta(\hat{b},t)$$

Under this new parametrization the behavior of our knot invariant with respect to left versus right-handed crossing can be expressed in the elegant form of the \emph{skein relation}.

\begin{figure}[h!]
\centering
\begin{tikzpicture}
\node (one) [label=left: $\Delta \Bigg ($, label=right: $\Bigg )$] {};
\draw (-1/4,-1/4) --( 1/4,1/4);
\draw (-1/4,1/4)--(-1/12,1/12);
\draw (1/12,-1/12)--(1/4,-1/4);
\node (minus) [right=of one, xshift=-5mm] {$-$};
\node (two) [right=of minus, label=left: $\Delta \Bigg ($, label=right: $\Bigg )$] {};
\begin{scope}[shift={(two)}]
\draw (1/4,-1/4)--(-1/4,1/4);
\draw (-1/4,-1/4)--(-1/12, -1/12);
\draw ( 1/12, 1/12)--(1/4,1/4);
\end{scope}
\node (eq) [right=of two, xshift=-5mm] {$=$};
\node (three) [right=of eq, label=left: $z$ $ \Delta \Bigg ($, label=right: $\Bigg )$] {};
\begin{scope}[shift={(three)}]
\draw (-1/8,-1/4)--(-1/8,1/4);
\draw (1/8,-1/4)--(1/8,1/4);
\end{scope}
\end{tikzpicture}
\caption{The skein relation satisfied by the Alexander polynomial.}
\end{figure}

\begin{ex} Deduce the skein relation from the definition of the Alexander-Conway polynomial and Lemma \ref{detlemma}.
\end{ex}

\section{Jones Representation and Jones polynomial}
In a manner analogous to how the Alexander polynomial is defined in terms of the Burau representation, another link invariant, the \emph{Jones polynomial}, can be studied alongside the \emph{Jones representation}.
Computing the Alexander polynomial is easy in the sense of complexity theory: there exists a polynomial time algorithm to compute it. This is because obtaining a braid $b\in \B_n$ from a link diagram $L$, computing the Burau representation of $b$, and then taking its determinant can all be done in polynomial time in the number of crossings of $L$.

On the other hand, assuming that $P\ne NP$, that is, assuming the longstanding conjecture that the complexity classes corresponding to polynomial time and nondeterministic polynomial time are distinct, computing the Jones polynomial is hard, in the sense that there exists no polynomial time algorithm.

However, certain values of the Jones polynomial can be approximated in polynomial time by a quantum computer. In this section we introduce the necessary background material for constructing the Jones representation of the braid group: the quantum integers, the Temperley-Lieb/Temperley-Lieb-Jones algebra, and the Temperley-Lieb category. Section 4 covers the application of the Jones representation to quantum computing.



\subsection{Quantum integers}
We should conceptualize the quantum integers as deformations of the integers by  $q$, which we can either think of as generic (a formal variable) or a specific element of $\mathds{C}^{*}$.

\begin{defn}\footnote{There are two conventions in the literature when defining the quantum integers, depending on whether a factor of $1/2$ appears in the exponents; quantum $n$ is sometimes defined as $[n]_q = \frac{q^{n}-q^{-n}}{q-q^{-1}}$.} Let $n \in \mathds{Z}$. Then \emph{quantum } n, denoted $[n]_q$, is given by

$$[n]_q = \frac{q^{n/2} - q^{-n/2}}{q^{1/2}-q^{-1/2}}$$

\end{defn}

For instance, $[1]_q = 1$ and $[2]_q=q^{1/2}+q^{-1/2}$. It is an easy application of L'H\^opital's rule to show that $[n]_q \to n$ in the limit $q \to 1$,  recovering the integers. This shows we can truly think of $[n]_q$ as some deformation of $n$. However, one must take special care when performing arithmetic with quantum integers, since the familiar rules of arithmetic need not apply. However, there is one important relation from integer arithmetic that still holds, the ``quantum doubling" formula.
\begin{prop} $[2][n]=[n+1]+[n-1].$
\end{prop}

This identity will reappear once we have introduced the Temperley-Lieb algebra.

\subsection{The Temperley-Lieb algebra $TL_n(A)$}

Our goal is to find braid group representations with properties that are useful for quantum computation, and towards this end we pass through the \emph{Temperley-Lieb algebra}, or the \emph{Temperley-Lieb-Jones algebra}. To motivate the construction of the Temperley-Lieb algebras, we recall the following theorem that dictates how the the group algebra for a finite group $G$ decomposes into the irreducible representations of $G$ \cite{fultonharris}. 

\begin{thm} Let $G$ be a finite group, and $\mathds{C}[G]=\{\sum a_g g \mid a_g \in \mathds{C}\}$ be the group algebra of $G$ over $\mathds{C}$. Then
$$\mathds{C}[G]\cong \bigoplus_{i} (\dim V_i) V_i$$ where the $V_i$ are a complete set of representatives of the isomorphism classes of finite-dimensional irreducible representations of G.
\end{thm}

To illustrate the theorem we recall the representation theory of $S_3$. There are three irreducible representations: the trivial, sign, and permutation representations, say $U, U'$, and $V$, respectively. Then $\mathds{C}[S_3] = U \oplus U' \oplus 2V$.  Hence, as an algebra, $\mathds{C}[S_3]$ decomposes as $\mathds{C}\oplus\mathds{C} \oplus M_2(\mathds{C})$.

While we can completely describe $\mathds{C}[G]$ when $G$ is finite, when $G$ is infinite, as in the case of $G=\B_n$, we don't have the same luxury. In order to get a handle on $\mathds{C}[\B_n]$ we pass to a finite-dimensional quotient. The first step in this process is to construct the \emph{Hecke algebra}.

\subsubsection{The Hecke algebra $H_n(q)$}

Hereafter we will work in one of two fields, $\mathds{C}$ or $Q(A)$, which we use to denote the field of rational functions in the \emph{Kauffman variable} $A$ over $\mathds{C}$. When we are interested in the \emph{generic} Temperley-Lieb algebra, we work in $Q(A)$, while in general we work in $\mathds{C}$. For now we use $\mathds{F}$ to denote the field $Q(A)$.

The elements of the braid group algebra $\mathds{F}[\B_n]= \{ \sum a_g g \mid g \in \B_n, a_g \in \mathds{F}\}$ are called formal (or quantum) braids. To motivate what relations we should quotient out by, we record a few observations.

Recall the presentation of the braid group
$$\B_n = \Big \langle \sigma_1, \sigma_2, \ldots, \sigma_{n-1} \Bigm\lvert \sigma_i\sigma_j = \sigma_j \sigma_i \text{ for } |i -j| \ge 2, \sigma_i\sigma_{i+1}\sigma_i = \sigma_{i+1}\sigma_i\sigma_{i+1}, i=1, 2, \ldots, n-1 \Big \rangle.$$

Taking the quotient of $\B_n$ by the normal subgroup generated by the $\sigma_i^2$, results in a group isomorphic to $S_n$. Thus there is a surjection of the braid group onto the symmetric group, and we have an exact sequence

$$1 \longrightarrow P\B_n \longrightarrow \B_n \longrightarrow S_n \longrightarrow 1.$$

This implicitly defines $P\B_n$, the \emph{pure braid group} on $n$-strands,  which will be revisited in Section 4. In particular, we can get a representation of the braid group by precomposing with a representation of the symmetric group. However, such a representation will not encode all of the information about the braid group that is need for computation, so it will not be interesting for us. Therefore we need to look for representations which do not factor through $S_n$.

Consider the quotient of $\mathds{F}[\B_n]$ by the relation $ \sigma_i^2 = a\sigma_i +b$, for $i=1, \ldots, n-1$, where $a$ and $b$ are independent of $i$. Note however, that $a$ and $b$ are not independent of one another, since we can rescale by setting $\tilde{\sigma}_i=\sigma_i/a$. Then the relation becomes

$$\tilde{\sigma_i}^2 = \tilde{\sigma_i} + b/a^2$$ In other words, we can just take $a=1$, so that the relation is parametrized by $b$. Taking the quotient of $\mathds{F}[\B_n]$ by this relation defines a Hecke algebra.

\begin{defn} The Hecke algebra $H_n(A)$ is the quotient  ${\mathds{F}}[\B_n]/I$ of the braid group algebra, where $I$ is the ideal generated by $\sigma_i^2 -(A-A^{-3})\sigma_i - A^{-2}$ for $i=1, \ldots, n-1$.
\end{defn}

A presentation on generators and relations of the Hecke algebra elucidates its structure further. Renormalizing via $q=A^{-4}$, and defining a new set of generators by $g_i=A^{-1}\sigma_i$, we can define
$$H_n(q)= \Big \langle g_1, g_2, \ldots, g_{n-1} \Bigm\lvert g_ig_j=g_jg_i \text{ for } |i -j| \ge 2, g_{i+1}g_{i}g_{i+1}=g_ig_{i+1}g_i, \text{ and } g_i^2=q^{-1}g_i+q, i=1, 2, \ldots, n-1 \Big \rangle.$$

Due to the Hecke relation $g_i^2=q^{-1}g_i + q$, $H_n(q)$ (and hence $H_n(A)$) is finite-dimensional.

\subsubsection{A presentation of the Temperley-Lieb algebra on generators and relations}
In order to obtain the Temperley-Lieb algebra, we must pass through one more quotient. Reparametrizing once again, rescaling the generators of $H_n(q)$ by defining $u_i=A\sigma_i - A^2$ and $d=-A^2 - A^{-2}=-[2]_q$, the Hecke algebra relations become

\begin{itemize}
\item $u_iu_j=u_ju_i$ when $|i-j|\ge 2$ \hfill (far commutativity)
\item $u_iu_{i+1}u_i - u_i = u_{i-1}u_iu_{i-1} - u_{i-1}$ \hfill (braid relation)
\item $u_i^2=du_i$ \hfill (Hecke relation)
\end{itemize}
To obtain the Temperley-Lieb algebra, we set the braid relation above to $0$, so impose the following additional relation that

\begin{itemize}
\item $u_iu_{i \pm 1}u_i =u_i$
\end{itemize}

\begin{defn} The generic Temperley-Lieb algebra $TL_n(A)$ is the quotient of the Hecke algebra $H_n(q)/I$, where $I$ is the ideal generated by $u_iu_{i \pm1}u_i-u_i$.
\end{defn}


The generic Temperley-Lieb algebra $TL_n(A)$ is \emph{semisimple} (also called a \emph{multi-matrix algebra}), a direct sum of matrix algebras $M_{n_i}(\mathds{F})$. This is the fact that enables us to work with matrix representations of the braid group, which, if unitary, can be thought of physically as quantum gates. Understanding how $TL_n(A)$ decomposes into matrix algebras is the key to applying the Jones representation to quantum computation.

\begin{thm} If $A$ is generic, then $TL_n(A)$ is semisimple. If $A \in \mathds{C}$, then $TL_n(A)$ is not in general semisimple.
\end{thm}

We will return to the semisimple structure of $TL_n(A)$ after introducing its picture presentation, in which computations can be performed using a graphical calculus.

\subsubsection{A picture presentation of the Temperley-Lieb algebra}
In the graphical calculus the variable $d=-A^2-A^{-2}$ previously defined in the context of the presentation of $TL_n(A)$ with generators and relations takes on an important role. The variable $d$ is called the \emph{loop variable}, for reasons that will soon be clear.

 \begin{defn} A diagram in $TL_n(A)$ is a square with $n$ marked points on the top edge and $n$ marked points on the bottom edge, and these $2n$ boundary points are connected by non-intersecting smooth arcs. In addition, there may be simple closed loops in the diagram.
 \end{defn} 

 An equivalent diagram can be obtained by multiplying by a factor of $d$ for each closed loop removed, and we say two diagrams are the same if they are \emph{d-isotopic}, that is, if they are isotopic and the boundary points are of the respective diagrams paired in the same way.

 An arbitrary element of $TL_n(A)$ is a formal sum of diagrams, where each diagram is a word in the generators $u_i$. The diagram of $u_i$ has a ``cup'' on the top edge connecting the $i$th and $i+1$st marked points, and a ``cap'' on the bottom edge connecting the $i$th and $i+1$st marked points. The $j$th marked point on the top edge is connected to the $j$th marked point on the bottom edge by a ``through strand''.

 \begin{figure}[h!]
\centering
\begin{tikzpicture}
\draw (0,0) rectangle (1,1);
\draw (1/8,1) arc (180:360:1/8);
\draw (1/8,0) arc (180:0:1/8);
\draw (1/2,0) --(1/2,1);
\draw (11/16, 1/2) node{\tiny$\cdots$};
\draw (7/8,0) --(7/8,1);
\draw (1/2,0) node[below] {$u_1$};
\draw (1.25,1/4) node[below] {$,$};
\end{tikzpicture}
\begin{tikzpicture}
\draw (0,0) rectangle (1,1);
\draw (1/8,0) -- (1/8,1);
\draw (1/4,0) arc (180:0:1/8);
\draw (1/4,1) arc (180:360:1/8);
\draw (5/8,0)--(5/8,1);
\draw (12/16,1/2) node{\tiny $\cdots$};
\draw (7/8,0)--(7/8,1);
\draw (1/2,0) node[below] {$u_2$};
\draw (1.5,1/4) node[below] {$, \ldots, $};
\end{tikzpicture}
\begin{tikzpicture}
\draw (0,0) rectangle (1,1);
\draw (1/16,0) -- (1/16,1);
\draw (3/16,0) node{\tiny $\cdots$};
\draw(5/16, 0) -- (5/16, 1);
\draw (3/8,1) arc (180:360:1/8);

\draw (1/2,0) node[below] {$u_i$};
\draw (1.5,1/4) node[below] {$, \ldots, $};

\draw (3/8,0) arc (180:0:1/8);
\draw (11/16,0) -- (11/16,1);
\draw (13/16, 1/2) node{\tiny$\cdots$};
\draw (15/16,0) --(15/16,1);
\end{tikzpicture}
\begin{tikzpicture}
\draw (0,0) rectangle (1,1);
\draw (5/8,1) arc (180:360:1/8);
\draw (5/8,0) arc (180:0:1/8);
\draw (1/2,0) --(1/2,1);
\draw (5/16, 1/2) node{\tiny$\cdots$};
\draw (1/8,0) --(1/8,1);
\draw (1/2,0) node[below] {$u_{n-1}$};

\end{tikzpicture}
\caption{The Temperley-Lieb diagrams of the multiplicative generators $u_i$.}
\end{figure}
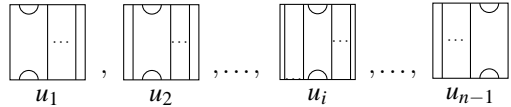

Multiplication of diagrams is performed by vertical stacking followed by rescaling- if $D_1, D_2 \in TL_n(A)$, then $D_1\cdot D_2$ is given by stacking $D_2$ on top of $D_1$ and rescaling to a square.

$$
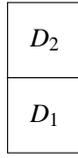
\begin{tikzpicture}
\node (rect) {};
\draw (0,0) rectangle (1,1);
\draw (1/2,1/2) node {$D_1$};
\draw (0,1) rectangle (1,2);
\draw (1/2,3/2) node {$D_2$};
\end{tikzpicture}$$

The Temperley-Lieb relations, far commutativity, the braid relation, and the Hecke relation, can all be verified using the graphical calculus. For example, the Hecke relation becomes

$$\begin{tikzpicture}
\node (rect) {};
\draw (-1/2,1/2) node {$u_i^2=$};
\draw (0,0) rectangle (1,1);
\draw (1/16,0) -- (1/16,1);
\draw (3/16,1/4) node{\tiny $\cdots$};
\draw (3/16,3/4) node{\tiny $\cdots$};
\draw(5/16, 0) -- (5/16, 1);
\draw (3/8,1) arc (180:360:1/8);

\draw[dashed] (0,1/2)--(1,1/2);
\draw (1/2,1/2) circle (1/8);
\draw (3/8,0) arc (180:0:1/8);
\draw (11/16,0) -- (11/16,1);
\draw (13/16, 1/4) node{\tiny$\cdots$};
\draw (13/16, 3/4) node{\tiny$\cdots$};
\draw (15/16,0) --(15/16,1);
\draw (3/2,4/8) node {$= d u_i$};
\end{tikzpicture}$$

\begin{ex} Show that the generators $u_i$ satisfy far commutativity and the braid relation.
\end{ex}

\begin{thm} The diagrammatic algebra for $TL_n(A)$ is isomorphic to the abstract Temperley-Lieb algebra given by generators and relations.
\end{thm}

The proof of this result is made difficult by the diagrammatic algebra being defined up to $d$-isotopy.

As a vector space, $TL_n(A)$ is generated by all the Temperley-Lieb diagrams in $TL_n(A)$, of which there are Catalan number $c_n = \frac{1}{n+1} \binom{2n}{n}$ many up to $d$-isotopy. In order to prove that the set of Temperley-Lieb diagrams forms a basis of $TL_n(A)$ as a vector space, we must show that there are no linear relations among the diagrams. This can be done by introducing an inner product on $TL_n(A)$, defined through the \emph{Kauffman bracket} and a map called the \emph{Markov trace}, which can be thought of as a ``quantum" analogue of the braid closure.

\subsubsection{The Kauffman bracket}
To find finite dimensional representations of the braid group, we begin by looking for an algebra homomorphism from the braid group algebra to finite matrix algebras,
$$\rho: \mathds{F}[\B_n] \to \bigoplus_{m} M_{m \times  m}(\mathds{F}).$$

A finite dimensional representation of $\B_n$ is obtained via the restriction of $\rho$ to the braid group, $\rho \bigm\lvert_{\B_n}$. The Kauffman bracket $\langle \cdot \rangle: \mathds{F}[\B_n] \to TL_n(A)$ is an algebra homomorphism which we can think of as producing \emph{Temperley-Lieb diagrams} by resolving crossings in the braid group algebra.

\begin{figure}[h!]
\center
\begin{tikzpicture}[scale=.5]
\node (one) {};
\draw (1/2,1) --( -1/2,-1) node[below] {\small i};
\draw (-1/2,1)--(-1/6,1/3);
\draw (1/6,-1/3)--(1/2,-1) node[below] {\small i+1};
\node (eq) [right=of one, xshift=-5mm] {$=$ };
\node (two) [right=of eq, xshift=-7.5mm] {$A$};
\begin{scope}[shift = {(two)}]
\draw (1/2,1)--(1/2,-1) node[below] {\small i};
\draw (3/2,1)--(3/2,-1) node [below] {\small i+1};
\end{scope}
\node (plus) [right=of two] {$+$};
\node (three) [right=of plus,xshift=-5mm] {$A^{-1}$};
\begin{scope}[shift = {(three)}]
\draw (1/2,1) arc (-180:0:1/2);
\draw (1/2,-1) node[below] {\small i} arc (180:0:1/2) node [below] {\small  i+1};
\end{scope}
\end{tikzpicture}
\caption{Illustration of the Kauffman bracket.}
\end{figure}

In terms of braid group generators and Temperley-Lieb generators, the Kauffman bracket is expressed by

$$\sigma_i = AI + A^{-1}u_i.$$

\subsubsection{The Markov trace}
  The Markov trace of a diagram is the map Tr: $TL_n(A) \mapsto \mathds{F}$ that sends the diagram $D$ to its tracial closure.

$$\begin{tikzpicture}
\node (diagram) {};
\draw (0,0) rectangle (1,1);
\draw (1/2,1/2) node {$D$};
\draw  (3/4,1) to [out=90, in=90] (9/8,1)--(9/8,0) to [out=-90,in=-90] (3/4,0);
\draw  (1/2,1) to [out=90, in=90] (10/8,1)--(10/8,0) to [out=-90,in=-90] (1/2,0);
\draw (5/16,-1/16) node {\tiny $\cdots$};
\draw (5/16,17/16) node {\tiny $\cdots$};
\draw (1/8,1)  to [out=90, in=90] (12/8,1)--(12/8,0) to [out=-90,in=-90] (1/8,0);
\draw (11/8,1/2) node {\tiny $\cdots$};
\node (eq) [right=of diagram, xshift=4mm,yshift=6mm] {$= d^{\text{ \# loops }}$};
\end{tikzpicture}$$
This defines the Markov trace on a basis of diagrams of $TL_n(A)$, and by extending linearly it is defined on all of $TL_n(A)$. From the trace one can define an inner product $\langle \cdot, \cdot \rangle: TL_n(A) \times TL_n(A)\to \mathds{F}$ given by

$$\langle D_1, D_2 \rangle = \text{Tr} (\overline{D_1}D_2).$$
where the bar over a diagram denotes the diagram obtained by reflecting across the horizontal midline.
$$
\begin{tikzpicture}
\node (diagram) {};
\draw (0,0) rectangle (1,1);
\draw (1/2,1/2) node {$\bar{D}$};
\node (eq) [right=of diagram, xshift=1mm, yshift=5mm] {};
\begin{scope}[shift= {(eq)}]
\draw (-3/16,0) node {$=$};
\draw (0,-1/2) rectangle (1,1/2);
\draw (1/2,0) node {\raisebox{\depth}{\scalebox{1}[-1]{$D$}}};
\end{scope}
\end{tikzpicture}$$


Now the question of whether the diagrams in $TL_n(A)$ are linearly independent can be translated into the question of whether the Gram matrix $(M)_{ij} = \langle \overline{D_i}, D_j \rangle$ has determinant zero.
$M$ is a $c_n \times c_n$ matrix, where $c_n$ is the $n$th Catalan number.
While the details are not provided here, it is possible to express the determinant of $M$ in the closed form

$$\det(M) = \prod_{i=1}^n \Delta_i(d)^{a_{n,i}}$$
where $a_{n,i}= \binom{2n}{n-i-2} + \binom{2n}{n-i} - 2\binom{2n}{n-i-1}$ and $\Delta_i(x)$ is the $i$th Chebyshev polynomial of the second kind, defined recursively by $\Delta_0=1, \Delta_1=x$, and $\Delta_{i+1}=x\Delta_{i} -\Delta_{i-1}$.

Thus whenever the loop variable $d$ is a root of a Chebyshev polynomial appearing in the determinant of the Gram matrix, there is some linear dependence among the Temperley-Lieb diagrams. The following table lists the first few Chebyshev polynomials.
\begin{figure}[h!]
\centering
\begin{tabular}{|c||c|c|c|c|c|}
\hline
$n$ & 0 & 1 & 2 & 3 & 4  \\
\hline
$\Delta_n(x)$ & 1 & $x$ & $x^2-1$ & $x^3-2x$ & $x^4-3x^2+1$  \\
\hline
\end{tabular}
\end{figure}

The roots of $\Delta_4$ are related to the golden ratio $\phi$, which will come up again later. The Chebyshev polynomials are also related to the quantum integers, which we will see later.

\subsection{The Jones polynomial}

To motivate the form that the Jones polynomial takes, we investigate the properties that would be needed for a quantity to give an invariant of a link. Given a braid $b \in \B_n$, we can apply the Kauffman bracket $\langle \cdot \rangle$ to resolve the crossings, resulting in a sum of $2^n$ Temperley-Lieb diagrams. Then $\langle b \rangle \in TL_n(A)$, and we can apply the Markov trace. Thus we can consider their composition $\text{Tr}\langle \cdot \rangle : \B_n \to \{\text{links}\}$. For example,

\begin{figure}[h!]
\centering
\begin{tikzpicture}[scale=.75]
\node (second) {};
\begin{scope}[shift=(second),xshift=-5mm, yshift=13mm]
\draw (1/2,-2/3) node{\large Tr$\langle$};
 \braid a_1^{-1};
 \draw (19/8,-2/3) node{\large $\rangle$};
\end{scope}
\node (third)[right=of second, xshift=17mm, yshift= 5mm] {};
\begin{scope}[shift={(third)}]
\draw (-3/2,0) node {$= A$};
\draw (0,0) circle (1/2);
\draw (0,0) circle (1);
\draw (1.75, 0) node{$+$  $A^{-1}$};
\node (fourth)[right=of third,xshift=8mm] {};
\end{scope}
\begin{scope}[shift={(fourth)}]
\draw (1,-3/4)--(1,3/4) to [out=90, in=90] (-1/4,3/4)--(-1/4,1/4) to [out=-90, in=-90] (1/4,1/4) --(1/4,3/4) to [out=90, in=90] (1/2,3/4)--(1/2,-3/4) to [out=-90, in=-90] (1/4,-3/4)--(1/4,-1/4) to [out=90, in=90] (-1/4,-1/4)--(-1/4,-3/4) to [out=-90, in=-90] (1,-3/4);
\draw (3.5,0) node {$=Ad^2 + A^{-1}d$ $=$ $-A^{3}d$};
\end{scope}
\end{tikzpicture}
\end{figure}
A similar computation shows that the Markov trace of the Kauffman bracket applied to the left handed crossing evaluates to $-A^{-3}$. However, if we calculate the trace of the (topologically equivalent) unknot, the result is $d$. While Tr$\langle \cdot \rangle$ is too sensitive to provide a knot invariant, it can be calibrated by multiplying a factor of $(-A^{-3})^{e(b)}$, where $e(b)$ is the writhe of the braid $b$ introduced in Section 2.

We are finally ready to define the Jones polynomial\footnote{Technically $J(L,q)$ is a Laurent polynomial in $q^{1/2}$, but it is still referred to as a polynomial in the literature.} of a link.

\begin{defn} Let $b \in \B_n$, and let $L=\hat{b}$ be the link obtained from the braid closure of $b$. Then the Jones polynomial $J(L,q)$ of $L$ is given by

$$J(L,q)=\frac{(-A^{-3})^{e(b)} \text{Tr}\langle b \rangle}{d}.$$
\end{defn}

The reason for the factor of $d$ in the denominator is to set the convention that the Jones polynomial of the unknot be equal to 1. It is necessary to point out the the Jones polynomial is not well defined if we try to evaluate at a general link instead of a braid closure - in order to make sense of the Jones polynomial of a link, it must be oriented.


By the Markov theorem, we know that whenever two links arising from braid closures are equal, they are related by a finite number of Markov moves.

Therefore, we must verify that the Jones polynomial is invariant under the Markov moves. This turns out to be an almost identical argument that was used when demonstrating the invariance of the Alexander polynomial under the braid closure.

Let $a,b \in \B_n$. Invariance under conjugation can be seen by sliding diagrams around their tracial strands. The proof-by-picture is identical to that provided for the proof of invariance of the Markov moves under braid closure, except one replaces the braid diagram $b$ by a Temperley-Lieb diagram. Similarly if $a=b\sigma_n^{\pm1}$, then the diagrammatic proof of invariance under stabilization is analogous to that for braids, except now we introduce a factor of $-A^{- 3}$ to correct for the writhe.

\subsubsection{Example: the Jones polynomial of the trefoil knot}
We end this discussion of the Jones polynomial with a famous example - the right-handed trefoil knot, $\widehat{\sigma_1^3}$.

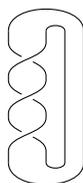
\begin{figure}[h!]
\centering
\begin{tikzpicture}[scale=.5]
\braid a_1^{-1}a_1^{-1}a_1^{-1};
\draw (2,0) to [out=90, in=90] (5/2,0)--(5/2,-7/2) to [out=-90, in=-90] (2,-7/2);
\draw (1,0) to  [out=90, in=90] (3,0)--(3,-7/2) to [out=-90, in=-90] (1,-7/2);
\end{tikzpicture}
\end{figure} There are three crossings, and hence $2^3=8$ terms in the resolution $\langle \sigma_1^3 \rangle$. The terms can be organized by a binary tree of depth 3, where each edge is labeled by a $``+"$ or a $`-"$ according to the term in the Kauffman bracket. We compute the `` -- + --'' term as an example.

\begin{figure}[h!]
\centering
\begin{tikzpicture}[scale=.5]
\node (ppm) {};
\draw (2,0) to [out=90, in=90] (5/2,0)--(5/2,-7/2) to [out=-90, in=-90] (2,-7/2);
\draw (1,0) to [out=-90, in=-90] (2,0);
\draw (1,0) to  [out=90, in=90] (3,0)--(3,-7/2) to [out=-90, in=-90] (1,-7/2);
\draw (1,-7/2) to [out=90, in=90] (2,-7/2);
\draw (1,-5/2)--(1,-1) to [out=90, in=90] (2,-1)--(2,-5/2) to [out=-90, in=-90] (1,-5/2);
\draw (7,-7/4) node {$=-A^{-3}\cdot( A^{-1}\cdot A \cdot A^{-1}) \cdot d^2$};
\end{tikzpicture}
\end{figure}

One can check that $J(\widehat{\sigma_1^3},q)= q + q^3 - q^4$.
It turns out that the Jones polynomial of the left-handed trefoil knot is different. This is an improvement over the Alexander polynomial, which cannot distinguish a knot from its mirror image.

\begin{prob} Does there exist a non-trivial knot with the same Jones polynomial as the unknot? Are there knots which are topologically different from their mirror image but have the same Jones polynomial? How does one interpret the Jones polynomial in terms of classical topology?
\end{prob}

\begin{ex} Let $I_K$ denote an invariant of knots. Then $I_K$ naturally extends to knots with \emph{double points}, via

\begin{center}
\begin{tikzpicture}
\node (zero) [label=left: $I_K \Bigg ($, label=right: $\Bigg )$] {};
\node (eq) [right=of zero, xshift=-5mm] {$=$};
\node (one) [right=of eq, label=left: $I_K \Bigg ($, label=right: $\Bigg )$] {};
\draw (-1/4,-1/4) --( 1/4,1/4);
\draw (-1/4,1/4)--(1/4,-1/4);
\draw[fill] (0,0) circle (.05) {};
\begin{scope}[shift={(one)}]
\draw (-1/4,-1/4) --( 1/4,1/4);
\draw (-1/4,1/4)--(-1/12,1/12);
\draw (1/12,-1/12)--(1/4,-1/4);
\end{scope}
\node (minus) [right=of one, xshift=-5mm] {$-$};
\node (two) [right=of minus, label=left: $I_K \Bigg ($, label=right: $\Bigg )$] {};
\begin{scope}[shift={(two)}]
\draw (1/4,-1/4)--(-1/4,1/4);
\draw (-1/4,-1/4)--(-1/12, -1/12);
\draw ( 1/12, 1/12)--(1/4,1/4);
\end{scope}
\end{tikzpicture}
\end{center} Now let $k$ be a knot, and let $q=e^{h}$, where $h$ is a formal variable. Then the Jones polynomial has the property that
$$J(k,e^h)= \sum v_i h^i$$
where the $v_i$ are knot invariants and the series is potentially infinite. The Alexander polynomial has the same property, and in the series expansion
$$\Delta(k,z) = 1+ c_2 z^2 + c_4z^4 + \cdots$$
the $c_i$ are known as the \emph{Vassiliev invariants}.\\
(a) Show that $v_1=0$. \\
(b) Show that if there are more than three double points, $c_2=v_2=0$. \\
(c) Show that $v_2 = -2 c_2$ 

\end{ex}

\subsection{The generic Jones representation}
While the Alexander polynomial was defined in terms of the Burau representation of the braid group, we were able to formulate a definition of the Jones polynomial that did not depend on the Jones representation. However, in order to approximate the Jones polynomial on a quantum computer, the Jones representation of the braid group must be understood. The definition of the  Jones representation depends on how the Temperley-Lieb algebras decompose into direct sums of matrix algebras.

\begin{defn} The Jones representation $\rho_{k,n}$ at level $k$ of the $n$-strand braid group is given by the image $$\langle \cdot \rangle : \mathds{\mathds{F}}[\B_n] \to TL_n(A)= \bigoplus_{n_i} M_{n_i} (\mathds{F}).$$
\end{defn}

To be precise, the braid group algebra $\mathds{F}[\B_n]$ is resolved into the Temperley-Lieb algebra $TL_n(A)$ via the Kauffman bracket $\langle \cdot \rangle$, and then after identifying the Temperley-Lieb algebra as a direct sum of matrix algebras $\bigoplus_i M_{n_i}(\mathds{F})$, restricting to the braid group gives a representation of $\B_n$.

There is a standard way to decompose an algebra as a direct sum of matrix algebras by finding its \emph{matrix elements}, elements $e_{ij}$ satisfying $e_{ij}e_{kl} = \delta_{jk} e_{il}.$ Before presenting the general method for finding such elements in $TL_n(A)$, we work out the solution for some small values of $n$.

\subsubsection{$TL_2(A)$}
 $TL_2(A)$ is generated as a vector space by the diagrams
\begin{tikzpicture}[scale=.5]
\draw (0,0) rectangle (1,1);
\draw (3/4,0) --(3/4,1);
\draw (1/4,0) --(1/4,1); \end{tikzpicture}
and
\begin{tikzpicture}[scale=.5]
\draw (0,0) rectangle (1,1);
\draw (1/4,1) arc (180:360:1/4);
\draw (1/4,0) arc (180:0:1/4);
\end{tikzpicture}
which we denote by $1$ and $u_1$, respectively. We look for elements $e_1$ and $e_2$ satisfying $e_1^2=e_1$, $e_2^2=e_2$, $e_1e_2=e_2e_1=0$.

For then we could identify
$$e_1 \longrightarrow \begin{pmatrix} 1 & 0 \\ 0 & 0 \end{pmatrix}, e_2 \longrightarrow \begin{pmatrix} 0 & 0 \\ 0 & 1 \end{pmatrix}.$$

We check that taking $e_1 = 1-\frac{1}{d} u_1$ and $e_2 = \frac{1}{d}u_1$ gives two matrix units.

Indeed, the following equations show that idempotency and centrality follow from the Hecke relation. $$e_1^2= (1-\frac{1}{d} u_1)(1-\frac{1}{d} u_1) = 1 -\frac{2}{d} u_1 + \frac{1}{d^2} u_1^2 = 1 - \frac{2}{d} u_1 + \frac{d}{d^2} u_1 = 1-\frac{1}{d} u_1 = e_1,$$

$$e_2^2 = (\frac{1}{d}u_1)^2 = \frac{1}{d^2} u_1^2 = \frac{1}{d}u_1 = e_2, \text{ and }$$
$$e_1e_2 = (1-\frac{1}{d} u_1)\frac{1}{d}u_1 =e_2e_1= \frac{1}{d}u_1 - \frac{1}{d}u_1 = 0.$$

Therefore, the decomposition of a Temperley-Lieb algebra at $n=2$ is given by $TL_2(A) \cong \mathds{F} \oplus \mathds{F}$.

\subsubsection{$TL_3(A)$}
As a vector space $TL_3(A)$ is spanned by  \begin{tikzpicture}[scale=.5]
\draw (0,0) rectangle (1,1);
\draw (1/4,0) --(1/4,1);
\draw (1/2,0) --(1/2,1);
\draw (3/4,0) --(3/4,1);
\end{tikzpicture},
\begin{tikzpicture}[scale=.5]
\draw (0,0) rectangle (1,1);
\draw (1/4,1) arc (180:360:1/8);
\draw (1/4,0) arc (180:0:1/8);
\draw (3/4,0) --(3/4,1);

\end{tikzpicture},
\begin{tikzpicture}[scale=.5]
\draw (0,0) rectangle (1,1);
\draw (1/2,1) arc (180:360:1/8);
\draw (1/2,0) arc (180:0:1/8);
\draw (1/4,0) --(1/4,1);

\end{tikzpicture},
\begin{tikzpicture}[scale=.5]
\draw (0,0) rectangle (1,1);
\draw (1/4,1) to [out=-90,in=90] (3/4,0);
\draw (1/2,1) arc (180:360:1/8);
\draw (1/4,0) arc (180:0:1/8);
\end{tikzpicture}, and
\begin{tikzpicture}[scale=.5]
\draw (0,0) rectangle (1,1);
\draw (3/4,1) to [out=-90,in=90] (1/4,0);
\draw (1/4,1) arc (180:360:1/8);
\draw (1/2,0) arc (180:0:1/8);
\end{tikzpicture}.
By a dimension argument, it is immediate that if $TL_3(A)$ is to be a matrix algebra, then we must have $TL_3(A) \cong \mathds{F} \oplus M_{2}(\mathds{F})$. Thus we will want to look for an idempotent element $p$ to identify with $\begin{pmatrix} 1 &  0 & 0 \\ 0 & 0 & 0 \\ 0 & 0 & 0 \end{pmatrix}$.

We claim that the element
\begin{figure}[h!]
\centering
\begin{tikzpicture}
\node (first) {};
\draw (0,0) rectangle (1,1);
\draw (1/4,0) --(1/4,1);
\draw (1/2,0) --(1/2,1);
\draw (3/4,0) --(3/4,1);
\draw (15/8, 1/2) node {$+$ $\frac{1}{d^2-1}$ $\Bigg ($};
\node (second) [right=of first, xshift=15mm] {};
\begin{scope}[shift={(second)}]
\draw (0,0) rectangle (1,1);
\draw (1/4,1) to [out=-90,in=90] (3/4,0);
\draw (1/2,1) arc (180:360:1/8);
\draw (1/4,0) arc (180:0:1/8);
\end{scope}
\node (third) [right=of second, xshift=4mm]{};
\begin{scope}[shift={(third)}]
\draw (-1/3,1/2) node {$+$};
\draw (0,0) rectangle (1,1);
\draw (3/4,1) to [out=-90,in=90] (1/4,0);
\draw (1/4,1) arc (180:360:1/8);
\draw (1/2,0) arc (180:0:1/8);
\draw (11/8, 1/2) node {$\Bigg )$};
\draw (19/8, 1/2) node {$-$ $\frac{d}{d^2-1}$ $\Bigg ($};
\end{scope}
\node (fourth) [right=of third, xshift=20mm] {};
\begin{scope}[shift={(fourth)}]
\draw (0,0) rectangle (1,1);
\draw (1/4,1) arc (180:360:1/8);
\draw (1/4,0) arc (180:0:1/8);
\draw (3/4,0) --(3/4,1);
\draw (4/3,1/2) node {$+$};
\end{scope}
\node (fifth) [right=of fourth, xshift=4mm]{};
\begin{scope}[shift={(fifth)}]
\draw (0,0) rectangle (1,1);
\draw (1/2,1) arc (180:360:1/8);
\draw (1/2,0) arc (180:0:1/8);
\draw (1/4,0) --(1/4,1);
\draw (11/8, 1/2) node {$\Bigg )$};
\end{scope}
\end{tikzpicture}
\end{figure}\\
has the desired property. We will come to know this element of $TL_3(A)$ as the Jones-Wenzl projector $p_3$.

One can see that the following $\tilde{e}_{ij}$, once properly normalized, give a set of matrix units. 
$$\tilde{e}_{11} = \begin{tikzpicture}[baseline=10] \draw (0,0) rectangle (1,1);
\draw (1/4,1) arc (180:360:1/8);
\draw (1/4,0) arc (180:0:1/8);
\draw (3/4,0) --(3/4,1); \end{tikzpicture}, \hspace{5pt} \tilde{e}_{21} =
\begin{tikzpicture}[baseline=10]
\node (first) {};
\draw (0,0) rectangle (1,1);
\draw (3/4,1) to [out=-90,in=90] (1/4,0);
\draw (1/4,1) arc (180:360:1/8);
\draw (1/2,0) arc (180:0:1/8);
\draw (11/8,1/2) node {$-$ $\frac{1}{d}$};
\node (second) [right=of first, xshift=5mm] {};
\begin{scope}[shift={(second)}]
\draw (0,0) rectangle (1,1);
\draw (1/4,1) arc (180:360:1/8);
\draw (1/4,0) arc (180:0:1/8);
\draw (3/4,0) --(3/4,1);
\end{scope}
 \end{tikzpicture}, \hspace{5pt} \tilde{e}_{12} = \begin{tikzpicture}[baseline=10]
\node (first) {};
\draw (0,0) rectangle (1,1);
\draw (1/4,1) to [out=-90,in=90] (3/4,0);
\draw (1/2,1) arc (180:360:1/8);
\draw (1/4,0) arc (180:0:1/8);
\draw (11/8,1/2) node {$-$ $\frac{1}{d}$};
\node (second) [right=of first, xshift=5mm] {};
\begin{scope}[shift={(second)}]
\draw (0,0) rectangle (1,1);
\draw (1/2,1) arc (180:360:1/8);
\draw (1/2,0) arc (180:0:1/8);
\draw (1/4,0) --(1/4,1);
\draw (11/8, 1/2);
\end{scope}
 \end{tikzpicture}$$

 $$\tilde{e}_{22} =
 \begin{tikzpicture}[baseline=10]
 \node (first) {};
 \draw (0,0) rectangle (1,1);
\draw (1/2,1) arc (180:360:1/8);
\draw (1/2,0) arc (180:0:1/8);
\draw (1/4,0) --(1/4,1);
\draw (11/8,1/2) node {$-$ $\frac{1}{d}$};
\node (second) [right=of first,xshift=5mm] {};
\begin{scope}[shift={(second)}]
\draw (0,0) rectangle (1,1);
\draw (3/4,1) to [out=-90,in=90] (1/4,0);
\draw (1/4,1) arc (180:360:1/8);
\draw (1/2,0) arc (180:0:1/8);
\draw (11/8,1/2) node {$-$ $\frac{1}{d}$};
\end{scope}
\node (third) [right=of second,xshift=5mm]{};
\begin{scope}[shift={(third)}]
\draw (0,0) rectangle (1,1);
\draw (1/4,1) to [out=-90,in=90] (3/4,0);
\draw (1/2,1) arc (180:360:1/8);
\draw (1/4,0) arc (180:0:1/8);
\draw (11/8,1/2) node {$+$ $\frac{1}{d^2}$};
\end{scope}
\node (fourth) [right=of third,xshift=5mm] {};
\begin{scope}[shift={(fourth)}]
\draw (0,0) rectangle (1,1);
\draw (1/4,1) arc (180:360:1/8);
\draw (1/4,0) arc (180:0:1/8);
\draw (3/4,0) --(3/4,1);
\end{scope}
 \end{tikzpicture} $$

Already in the case of $n=3$, our minimal central idempotents are becoming unwieldy. We remark that if $TL_n(A) \cong \bigoplus M_{n_i}(\mathds{F})$, then the dimension of the Temperley-Lieb algebra, namely the Catalan number $c_n$, can be written as a sum of squares. $$\frac{1}{n+1} \binom{2n}{n} = \sum_{i} n_i^2$$ For example, this relation gives $14=1+4+9$, $42= 1+16+25$, and $132 = 1+25 +25 + 81$.

The general theory of how $TL_n(A)$ decomposes into matrix algebras is governed by elements of the Temperley-Lieb algebra called the \emph{Jones-Wenzl projectors}, which are minimal central idempotents. While the details of the general theory are beyond the scope of these notes, studying the Jones-Wenzl projectors is nonetheless essential for using the graphical calculus to compute the Jones representation, for which we compute an example in Section 4.

\subsection{The Jones-Wenzl projectors $p_n$}
The following theorem characterizes the Jones-Wenzl projectors.
\begin{thm} There exists a unique nonzero element $p_n$ in $TL_n(A)$ such that \\
(1) $p_n^2=p_n$\\
(2) $u_ip_n=p_nu_i=0$ for all $i=1, \ldots, n-1$.
\end{thm}
Proof. First we prove uniqueness. Suppose that $p_n$ and $p_n'$ satisfy (1) and (2), and write $p_n=c_0\cdot 1+\sum_{i=1}^{n-1}  c_i u_i$, where we have simply used the fact that $\{1, u_1, \ldots, u_{n-1}\}$ form a basis of $TL_n(A)$ as a vector space. Immediately, the condition $p_n^2=p_n$ forces $c_0=1$, since

$$p_n^2 = (c_0\cdot 1+\sum_{i=1}^{n-1}  c_i u_i)^2 = c_0^2 + \cdots = p_n.$$
Hence $c_0^2=c_0$, and so $c_0=1$. So we can write $p_n=1+\sum_{i=1}^{n-1}  c_i u_i$ and $p_n' = 1+\sum_{i=1}^{n-1}  c_i' u_i$. Then we can check that $p_n'=p_np_n'=p_n$.

As for existence, we provide an explicit construction. We will use the notation
$$\begin{tikzpicture}
\node (first) {};
\draw (-1/2,1/2) node {$p_{n}=$};
\draw (1/4,1/2) rectangle (1/2,3/4);
\draw (3/8, 0) -- (3/8, 2/4);
\draw (3/8, 3/4) -- (3/8, 5/4);
\draw (5/8, 1/4) node {$n$};
\draw (1,1/2) node {$=$};
\node (second) [right=of first] {};
\begin{scope}[shift={(second)}]
\draw (1/4, 0) -- (1/4, 5/4);
\draw (4/8, 1/4) node {$n$};
\end{scope}
\end{tikzpicture}$$
 to represent the $n$th Jones-Wenzl projector.

We define $p_1=$
\begin{tikzpicture}[scale=.5]
\draw (0,0) rectangle (1,1);
\draw (1/2,0)--(1/2,1);
\end{tikzpicture}, $p_2=$ \begin{tikzpicture}[scale=.5]
\draw (0,0) rectangle (1,1);
\draw (3/4,0) --(3/4,1);
\draw (1/4,0) --(1/4,1);
\end{tikzpicture} $-\frac{1}{d}$
\begin{tikzpicture}[scale=.5]
\draw (0,0) rectangle (1,1);
\draw (1/4,1) arc (180:360:1/4);
\draw (1/4,0) arc (180:0:1/4);
\end{tikzpicture} and define $p_n$ inductively by

$$\begin{tikzpicture}
\node (first) {};
\draw (-1/2,1/2) node {$p_{n+1}=$};
\draw (1/4,1/2) rectangle (1/2,3/4);
\draw (3/8, 0) -- (3/8, 2/4);
\draw (3/8, 3/4) -- (3/8, 5/4);
\draw (5/8, 1/4) node {$n$};
\node (second) [right= of first,xshift=-.25cm] {};
\begin{scope}[shift={(second)}]
\draw (1/2,1/2) node {$- $ $\frac{\Delta_{n-1}(d)}{\Delta_n(d)}$};
\draw (3/2, 3/4) rectangle (3,1);
\draw (3/2, 1/4) rectangle (3,1/2);
\draw (2, 0)--(2, 5/4);
\draw (7/4,0)--(7/4,5/4);
\draw (9/4,5/8) node {\tiny $\cdots$};
\draw (10/4,0)--(10/4,5/4);
\draw (11/4,5/4)--(11/4,13/16) to [out=-90,in=-90] (13/4,13/16)--(13/4,5/4);
\draw (11/4,0)--(11/4,7/16) to [out=90,in=90] (13/4,7/16)--(13/4, 0);
\end{scope}
\end{tikzpicture}$$
where $\Delta_n = \text{Tr}( p_n)$ is the Markov trace of the $n$th Jones-Wenzl projector. It is straightforward to verify that this construction results in idempotent objects that annihilate the generators of the Temperley-Lieb algebra.

\begin{prob} Does every diagram $u_i$ appear in $p_{n+1}$?
\end{prob}

\subsubsection{The recursive definition of $p_n$  and the quantum doubling formula}
We can calculate $\Delta_n$ explicitly - later we'll show that $\Delta_n = (-1)^n [n+1]$. For now we will take this as fact. 

We can recover the formula $[2][n]=[n+1]+[n-1]$ from quantum arithmetic that we proved previously.  Recall that $d=-A^2 - A^{-2}=-q^{1/2}-q^{-1/2}=-[2]$. Taking the Markov trace of both sides of the recursive formula for $p_{n+1}$ gives
$$\Delta_{n+1} = d \Delta_n - \frac{\Delta_{n-1}(d)}{\Delta_n(d)} \Delta_n.$$
and hence
$$(-1)^{n+1}[n+2]=(-1)^n[n+1](-[2]) - (-1)^{n-1}[n]$$
after simplifying this becomes
$$[n+2]=[2][n+1]-[n],$$
or
$$[n+2] + [n] = [2][n+1].$$ This is the ``quantum doubling'' formula shifted by one.


\subsection{The non-generic Jones representation of the braid group}
The generic Jones representation of the braid group is given by the image of a braid in the matrix algebra decomposition of the generic Temperley-Lieb algebra. What happens for specific $A \in \mathds{C}^{*}$?

\begin{ques} For which $A$ is the Jones representation unitary? How does one compute the Jones representation and study its properties? How can the Jones representation be used to perform quantum computation?
\end{ques}

Let $A \in \mathds{C}$ and suppose $\rho(\sigma_i)^{\dagger} = \rho(\sigma_i^{-1})$ for $i=1,2,\ldots, n$. Then assuming $u_i^{\dagger} = u_i$, the equation $$\rho(\sigma_i)^{\dagger}\rho(\sigma)= (A^{\dagger} + (A^{-1})^{\dagger}u_i^{\dagger})(A+A^{-1}u_i) = |A|^2 + |A^{-1}|^2 d u_i + u_i (A^{\dagger}A^{-1} + (A^{-1})^{\dagger}A) = 1$$

is satisfied when $|A|=1$. One can also show that $u_i^{\dagger}=u_i$ and $|A|=1$ are actually necessary conditions for $\rho$ to be a unitary representation. That is, $A \in S^1$.

More can be said about what happens for specific value of the Kauffman variable $A$. Recall the recursive definition of Jones-Wenzl projector $p_{n+1}$, which involves the coefficient,  $\frac{\Delta_{n-1}}{\Delta_n}$. Thus if $A$ is a root of $\Delta_n$, the Jones-Wenzl projector $p_{n+1}$ is undefined. The following proposition characterizes when $A$ is a root of $\Delta_n$.  \begin{prop} $\Delta_n = (-1)^n[n+1] = (-1)^n \frac{A^{2n+2}-A^{-2n-2}}{A^2-A^{-2}}$.
\end{prop}

\begin{cor} Whether or not $A \in S^1$, the Jones representation is well-defined when $A$ is not a root of unity.
\end{cor}

To see what can go wrong when $A$ is a root of unity, consider $TL_2(A)=\mathds{C}[1, u_1]$ for $A$ a primitive eighth root of unity. Then $-A^{-2}=A^{2}$ and hence $d=-A^{2}-A^{-2}=0$. But then $TL_2(A) = \mathds{C}[1,x]$ where $x^2=0$, which is not a matrix algebra. For suppose $TL_2(A)$ were a matrix algebra. Then the dimension would force the isomorphism $TL_2(A) \cong \mathds{C} \oplus \mathds{C}$, and there would exist two central idempotents $e_1$ and $e_2$. Let $e_1=a+bx$. Then
$$(a+bx)^2 = a^2 +2abx +b^2x^2 = a^2+2abx=a+bx,$$
which has no consistent solution.

This is illustrative of a general problem that we may not necessarily get a matrix algebra when $A$ is a root of unity. We bypass this difficulty by passing to a quotient of the Temperley-Lieb algebra which is semi-simple, called the \emph{Temperley-Lieb Jones algebra}. Then the non-generic Jones representation is defined in analogy with the generic definition, as the image of the braid group in a matrix algebra decomposition of $TLJ_n(A)$.

\subsection{The Temperley-Lieb-Jones algebra $TLJ_n(A)$}
Let $r \ge 3$, and assume $A$ is either a primitive $4r$th root of unity or a primitive $2r$th root of unity if $r$ is odd.
\begin{defn}  There exists a semisimple quotient of $TL_n(A)$, called the Temperley-Lieb-Jones algebra, denoted by $TLJ_n(A)$, formed by taking the quotient by the $(r-1)$st Jones-Wenzl projector $p_{r-1}$.
\end{defn}
\begin{prob} Given a finite-dimensional algebra, taking the quotient by a Jacobson radical gives a semisimple algebra. Is the Jones-Wenzl quotient the same as the Jacobson quotient?
\end{prob}

\subsection{The Temperley-Lieb category $TLJ(A)$}
The Jones representations of an $n$-strand braid group are determined by how the Temperley-Lieb algebra $TL_n(A)$ or $TLJ_n(A)$ decomposes into a matrix algebras. In order to make sense of the relationship between $TLJ_n(A)$ and anyons, we promote the Temperley-Lieb-Jones algebras $\{TLJ_n\}$ to a \emph{Temperley-Lieb category.}
  The objects of this category will be finite sets of points $a_1, \ldots, a_n$ in the unit interval  $[0,1]$, allowing for the empty set, each point colored by an element of the \emph{label set } $\mathcal{L}=\{0,1, \ldots, k\}$. At each marked point, there is a Jones-Wenzl projector $p_{a_n}$.


Given two objects, which we label $X_{a_1, \ldots, a_n}$ and $X_{b_1, \ldots, b_m}$, where the subscript indicates the integer labeling of the specified point, the morphisms are as follows. If $m-n \not\equiv 0 \mod 2$, then the only morphism between the two objects is the zero morphism. If however, $m+n$ is an even number, then the set of morphisms
is given by the span of all Temperley-Lieb-Jones diagrams connecting the points $X_{a_i}$ and $X_{b_j}$, together with disjoint unions of loops colored by natural numbers. That is, $\text{Hom}(X_{a_1, \ldots, a_n}, X_{b_1, \ldots, b_m}) = \mathds{F}[\text{all \textcolor{red}{colored} Temperley-Lieb diagrams connecting} \sum a_i + \sum b_j]$ modulo
\begin{itemize}
\item  $\bigcirc_i=\Delta_i$
\item relative $d$-isotopy
\item $p_{k+1}=0$.
\end{itemize}
Note that if the $p_{k+1}$ vanishes in $TLJ(A)$, then the recurrence relation for Jones-Wenzl projectors implies that $p_{m}$ vanishes for all $m \geq k+1$.
The following three properties of the category are immediate from the definition.
\begin{prop}
\begin{enumerate}
\item $TLJ(A)$ is a $\mathds{C}$-linear category
\item Hom$(X,X)$ is an algebra for all $X$
\item Hom$(X,Y)$ is a Hom$(X,X) - $Hom$(Y,Y)$ bimodule
\end{enumerate}
\end{prop}

\subsubsection{Trivalent vertices}
The trivalent vertex is the most fundamental part of the Temperley-Lieb category.

The following figure from \cite{cbms} gives the resolution of a labeled trivalent vertex into Temperley-Lieb diagrams.

$$
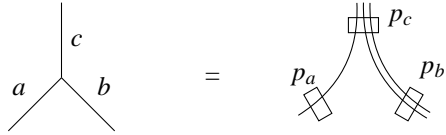
\begin{tikzpicture}
\node[shape=coordinate] (trivalent) {};
\draw (trivalent) -- node[auto,swap] {$c$} (90:1cm);
\draw (trivalent) -- node[auto] {$b$} (-45:1cm);
\draw (trivalent) -- node[auto,swap] {$a$} (225:1cm);
\node [right = of trivalent, minimum size = 2cm] (eq) {=};
\node[shape=coordinate, right = of eq] (v) {};
\begin{scope}[shift = {(v)}]
\draw (210:1cm) to[bend right] (95:1cm);
\draw[shift = {(210:.7cm)}, rotate = 120] (-.2cm,-.1cm) rectangle (.2cm, .1cm) node[auto, anchor= south] {$p_a$} ;
\draw (90:1cm) to[bend right] (-32.5:1cm);
\draw (85:1cm) to[bend right] (-27.5:1cm);
\draw[shift = {(90:.7cm)}] (-.2cm,-.1cm) rectangle (.2cm, .1cm) node[auto, anchor= west] {$p_c$};
\draw[shift = {(-30:.7cm)}, rotate = 60] (-.2cm,-.1cm) rectangle (.2cm, .1cm) node[auto, anchor= south west] {$p_b$} ;
\end{scope}
\end{tikzpicture}$$

 The labeling of the trivalent vertex is subject to the following conditions:
 \begin{enumerate}
 \item $a+b+c$ is even \hfill (``parity'')
 \item $a+b \ge c, b+c \ge a,$ and $c+a\ge b$ \hfill (``triangle inequality'')
 \item $a+b+c \le k$ \hfill (``positive energy condition'')
 \end{enumerate}
For example,
$$
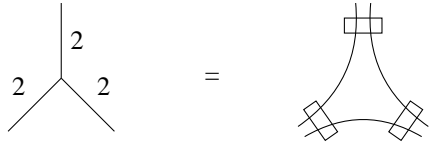
\begin{tikzpicture}
\node[shape=coordinate] (trivalent) {};
\draw (trivalent) -- node[auto,swap] {$2$} (90:1cm);
\draw (trivalent) -- node[auto] {$2$} (-45:1cm);
\draw (trivalent) -- node[auto,swap] {$2$} (225:1cm);
\node [right = of trivalent, minimum size = 2cm] (eq) {=};
\node[shape=coordinate, right = of eq] (v) {};
\begin{scope}[shift = {(v)}]
\draw[shift={(0:-.1cm)}] (220:1cm) to[bend right] (90:1cm); 
\draw[shift={(0:.1cm)}] (90:1cm) to[bend right] (-40:1cm); 
\draw (225:1.1cm) to[bend left] (-45:1.1cm); 
\draw[shift = {(225:.8cm)}, rotate = 125] (-.25cm,-.1cm) rectangle (.25cm, .1cm) node[auto, anchor= south] {} ; 

\draw[shift = {(90:.7cm)}] (-.25cm,-.1cm) rectangle (.25cm, .1cm) node[auto, anchor= west] {}; 
\draw[shift = {(-45:.8cm)}, rotate = 55] (-.25cm,-.1cm) rectangle (.25cm, .1cm) node[auto, anchor= south west] {} ; 

\end{scope}
\end{tikzpicture}$$


We summarize the important definitions into the following:
\begin{defn}
\begin{enumerate}
\item As an algebra, $TLJ_n(A)$ is the Hom space $n$ points on the unit interval, each marked by 1, with itself. We will denote this by Hom$({1}^{\otimes n},{1}^{\otimes n})$. More generally, the shorthand $a$ stands for the object with one point in the unit interval, marked by $a$.
\item The colored Temperley-Lieb-Jones algebra is given by Hom$(a^{\otimes n},a^{\otimes n})$, where $a \in \{0,1, \ldots, k\}.$
\item The Jones representation for $TLJ_n(A)$ is given by its image on $\bigoplus_{n_i} Hom(i, \mathds{1}^{\otimes n})$.
\end{enumerate}
\end{defn}

\subsubsection{Physical interpretation of $TLJ_n(A)$}

We want to have a physical interpretation to go along with our definition of $TLJ_n(A)$. Morphisms in the category, which are Temperley-Lieb-Jones diagrams, depict quantum processes of \emph{anyons}, the quasi-particle excitations of a 2D topological quantum system, such as those theorized to exist in fractional quantum Hall states.

In terms of the mathematical formalism:
\begin{defn} An object $X$ in the Temperley-Lieb-Jones category is simple if $Hom(X,X) \cong \mathds{C}$, in which case we say $X$ is an anyon.
\end{defn}

The number of distinct types of anyons is dictated by $k$, the level of the theory, and each type of anyon has an associated number, called its quantum dimension.

\begin{defn}
The quantum dimension of $a \in \mathcal{L}$, thought of as a representative of an isomorphism class of simple objects, is given by the loop value $d_a = \bigcirc_a$.
\end{defn}

The structure of the category captures the notion of fusion of particles.

\begin{defn}The fusion rules are the collection $\{N_{ab}^c=\dim Hom (a \otimes b, c) \mid a,b,c \in L\}$. More compactly, the fusion rules are implicitly defined through the equation
$$a \otimes b = \bigoplus_c N_{ab}^c c$$
where $c$ runs over the label set $\{0,1,2, \ldots, k\}$.
\end{defn}

Anyons generalize bosons (like photons) and fermions (like electrons) in two dimensions.

Given $n$ indistinguishable particles, with locations $x_1, \ldots, x_n$, then the type of particle is determined by what happens to their wavefunction $\psi(x_1, \ldots, x_n)$ upon interchanging their locations. For bosons, interchanging produces no change, while for fermions, a negative sign is generated by the interchange of any particles. That is

$$\psi(x_1, \ldots, x_i, \ldots, x_j, \ldots, x_n) = \pm \psi (x_1, \ldots, x_j, \ldots, x_i, \ldots, x_n)$$
depending on whether the particles are bosons or fermions. When we allow the wavefunction to be altered by an arbitrary phase $e^{i\theta}$, then we have an \emph{anyon}. For natural reasons one only considers rational phases of the form $e^{i\pi p/q}$ where $p,q \in \mathds{Z}$. The reason allowing an arbitrary phase produces this more general picture in two dimensions is due to the fact that there are no nontrivial knots in $\mathds{R}^4$. More precisely, if $S^1 \hookrightarrow \mathds{R}^4$ is an embedding, then the image of $S^1$ can always be isotoped to the trivial knot.

For topological quantum computation, we are interested in values of $k$, called the \emph{level} of the theory,  for which the corresponding topological phase of matter features anyons which are \emph{nonabelian}. Nonabelian anyons are those for which the representations of the braid group have non-abelian image for $n$ large enough.  Hence anyons can either be abelian or nonabelian, but in order to be useful for anyonic quantum computation they must be nonabelian.

\section{Anyon Systems and Anyonic Quantum Computation}

In this section we describe the algebraic theory of anyon systems, which is given by the Temperley-Lieb-Jones category $TLJ(A)$ for a fixed $A=\pm\im e^{\pm 2\pi i/4r}$, whose associated TQFT is  known as the \emph{Jones-Kauffman theory at level $k$}. The focus will be on two theories, the \emph{Ising theory} and the \emph{Fibonacci theory}.  In this section, we will use anyon and Jones-Wenzl projector interchangeably.  Anyons can be modeled by simple objects in unitary modular categories and Jones-Wenzl projectors represent simple objects in Jones-Temperley-Lieb categories.

Anyons can be harnessed to store and manipulate quantum bits, or \emph{qubits}, leading to a model of quantum computation whose topological nature lends it a special robustness. Braiding the anyons gives a \emph{quantum gate} that acts on qubits via the Jones representation. Given a specific anyon model of level $k$ described by a Temperley-Lieb Jones category $TLJ(A)$, understanding the image of the Jones representation $\rho_{k,n}:\B_n \to TLJ_n(A) \cong \bigoplus_{n_i}M_{n_i}(\mathds{C})$ is tantamount to assessing the power of the anyons to perform quantum computation. We will at least need the images of the braid group representations to be infinite and dense in order for the model to be \emph{universal} for quantum computation by braiding alone, that is, powerful enough to accurately and efficiently perform quantum computation.



This section is organized as follows. First we introduce the Ising and Fibonacci theories. Then for each of the two theories we investigate the dimensions of certain Jones representations by counting admissible labelings of fusion trees, and demonstrate how to encode a qubit with two dimensional representations. Then we show how to compute the Jones representation of the four-strand braid group at level 2, and trivial total charge. After introducing the $R$-symbols and $F$-symbols, we sketch how to compute the Jones representation of the three-strand braid group at level 3, with nontrivial total charge. Finally, with representations for the Ising theory and Fibonacci theory in hand, we present some results about their images and interpret the consequences for their corresponding anyonic models of quantum computation.

\subsection{Introduction}
We begin by setting the parameters of the theory $TLJ(A)$ that will describe our anyon model. Pick an integer $r\geq 3$,  and choose $A \in \{\pm i e^{\pm 2\pi i/4r}\}$ for unitarity. Then the \emph{level} of the theory for this choice of $A$ is $k=r-2$.  For each level, there are $4$ essentially equivalent theories, depending on which of the four choices of $A$ are made. Then the loop variable $d$ can be expressed in terms of the level by the equation
$$d=-A^{2}-A^{-2} = e^{\pm 4\pi i/4r} - e^{\mp 4 \pi i/4r} = 2(\cos \pi /r) = 2 \cos \frac{\pi}{k+2}.$$

The first few levels $k=1,2,3$ then correspond to $d=1, \sqrt{2}, \phi$, where $\phi=\frac{1+\sqrt{5}}{2}$ is the golden ratio.

\subsubsection{Level 1}

As a warmup to the $TLJ(A)$ theories that will be useful for quantum computation, we begin with level $k=1$.
The loop variable becomes $d=2\cos \frac{\pi}{3} = 1$, giving us the freedom to create and destroy loops as we please without having to account for them.

The $(k+1)$st Jones-Wenzl projector that is zero in $TLJ_n(A)$ is given by $p_2= \begin{tikzpicture}[baseline=.125cm,scale=.5]
\node (first) {};
\draw (0,0) rectangle (1,1);
\draw (3/4,0) --(3/4,1);
\draw (1/4,0) --(1/4,1);
\node (second) [right=of first, xshift=-3.5mm] {};
\begin{scope}[shift={(second)}]
\draw (-3/8,1/2) node {$-$};
\draw (0,0) rectangle (1,1);
\draw (1/4,1) arc (180:360:1/4);
\draw (1/4,0) arc (180:0:1/4);
\end{scope}
\end{tikzpicture}$
and hence $\begin{tikzpicture}[baseline=.125cm,scale=.5]
\node (first) {};
\draw (0,0) rectangle (1,1);
\draw (3/4,0) --(3/4,1);
\draw (1/4,0) --(1/4,1);
\node (second) [right=of first,xshift=-3.5mm] {};
\begin{scope}[shift={(second)}]
\draw (-3/8,1/2) node {$=$};
\draw (0,0) rectangle (1,1);
\draw (1/4,1) arc (180:360:1/4);
\draw (1/4,0) arc (180:0:1/4);
\end{scope}
\end{tikzpicture}$ in the level 1 theory. This category is equivalent to the category of super-vector spaces; it describes the trivial free fermion topological theory.

\subsubsection{The Ising and Fibonacci theories}
We first introduce the two anyon models in parallel, choosing
$$A = \begin{cases}  ie^{-2\pi i/16} & k=2  \text{ (Ising)} \\ ie^{2\pi i/20} & k = 3  \text{ (Fibonacci)}\end{cases}.$$
$TLJ(A)$ is a \emph{unitary modular category} (UMC) when $k$ is even, as for the Ising theory, and unitary pre-modular tensor category when $k$ is odd. As for $k=3$, it contains the Fibonacci sub-theory, which is a UMC of rank=$2$.

For level $k=2$ ($r=4$), the simple Jones-Wenzl projectors are $\{p_0, p_1, p_2\}.$ Thought of as anyons, the projectors have an alternative physical labeling $\{1, \sigma, \psi\}$, corresponding to the vacuum or ground state, \emph{Ising anyon}, and \emph{Majorana fermion}, respectively. The fusion rules for the Ising theory, in their most succinct form, are given by $1 \otimes x=x\otimes 1 = x$ for $x \in \{1,\sigma, \psi\}$, $\sigma \otimes \sigma = 1\oplus \psi$, $\sigma \otimes \psi = \psi \otimes \sigma = \sigma$, and $\psi \otimes \psi =1$. The last relation is particularly interesting--it tells us that the Majorana fermion $\psi$ is its own antiparticle. This is the famous Ising theory.

To prove that $1, \sigma,$ and $\psi$ are the only simple objects when $k=2$, we need to compute all $\textrm{Hom}(x,x)$. A nice way to do this is to use the inner product $\langle \cdot, \cdot \rangle$ that was previously defined on the Temperley-Lieb algebra in terms of the Markov trace.  Specializing $A$ to the particular root of unity, we have the following positivity.
\begin{prop}
For $A= \pm \im e^{\pm 2\pi i/4r}$, this inner product is positive definite on all $Hom(X,Y)$.
\end{prop}

For level $k=3$ ($r=5$), the simple Jones Wenzl projectors are $\{p_0, p_1, p_2, p_3\}$. The subset $\{p_0, p_2\}$ or $\{1,\tau\}$ corresponding to the vacuum and the \emph{Fibonacci anyon} generates the Fibonacci subtheory. The Fibonacci fusion rules are given by $1 \otimes \tau=\tau \otimes 1 = \tau$, and $\tau \otimes \tau = 1 \oplus \tau$.

\subsubsection{Notation}\label{TLJmatrix}

Unfortunately, two different things are both denoted by $1$'s: the $TLJ(A)$ label $1 \in \mathcal{L}$, and the ground state $1$ in an anyon system such as for the Ising and Fibonacci theories corresponding to $0 \in \mathcal{L}$.  Typically they will be clear from the context and we will use $\mathcal{L}$ when labeling diagrams to avoid confusion.

Having chosen $A=\pm i e^{\pm 2\pi i/r}$, we consider the Jones representation $\rho_{k,n,i}: \B_n \to TLJ_n(A) \to M_{n_i}(\mathds{C})$. Such a representation is parametrized by the level $k$ of the theory, the number of strands $n$ in the braid group, and the \emph{total charge} $i$.

Define the vector space $V_{k,n,i}$ to be the $\mathds{C}$-span of the \emph{fusion trees}
$$
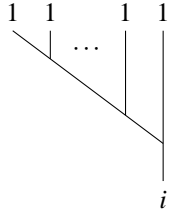
\begin{tikzpicture}
\draw (0,2)--(0,0); 
\draw (-2,2)--(0,.5);
\draw (-3/2, 2)--(-3/2,13/8);
\draw (-1/2,2)--(-1/2,7/8);
\draw (-1,7/4) node {$\cdots$};
\draw (0,-1/4) node {$i$};
\draw (-2,2.25) node {$1$};
\draw (-3/2,2.25) node {$1$};
\draw (-1/2,2.25) node {$1$};
\draw (0,2.25) node {$1$};
\end{tikzpicture}$$
where the internal edges are admissibly labeled by elements of the label set $\mathcal{L}=\{0,1,2, \ldots, k\}$. Physically, an admissible labeling of a fusion tree corresponds to a possible fusion process of the corresponding anyons.

Our ultimate goal is to understand the image of the Jones representation $\rho_{k,n,i}(\B_n)$ in $U(V_{k,n,i})$, the unitary transformations on the vector space $V_{k,n,i}$ and interpret them as quantum gates. As a first step we count admissible labelings of Ising and Fibonacci fusion trees to get the dimension of the representations for small $n$, looking for a two-dimensional representation in which to encode a qubit in order to get single-qubit gates. We will eventually also want a representation of at least dimension four, so that we can produce two-qubit gates. This turns out to be enough to show that Fibonacci anyons can be used for universal quantum computation.

As a warm up to the Ising and Fibonacci theories, we first consider $k=1$.
\subsubsection{$k=1$}

When $k=1$, the label set has two elements, $\mathcal{L}=\{0,1\}$. Depending on whether $n$ is even or odd, by a parity argument there is only one way to admissibly label the fusion tree by elements of $\mathcal{L}$.

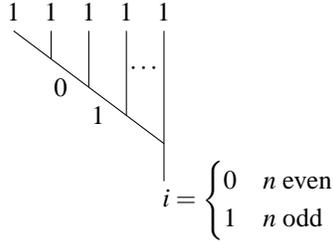
$$\begin{tikzpicture}
\draw (0,2)--(0,0); 
\draw (-2,2)--(0,.5);
\draw (-3/2, 2)--(-3/2,13/8);
\draw (-1/2,2)--(-1/2,7/8);
\draw (-1,2)--(-1, 10/8);
\draw (-1/4,6/4) node {$\cdots$};
\draw (9/8,-1/4) node {$i=\begin{cases} 0 & n \text{ even} \\ 1 & n \text{ odd}\end{cases}$};
\draw (-2,2.25) node {$1$};
\draw (-3/2,2.25) node {$1$};
\draw (-1, 2.25) node {$1$};
\draw (-1/2,2.25) node {$1$};
\draw (0,2.25) node {$1$};

\draw (-11/8, 1.25) node {$0$};
\draw (-7/8,7/8) node {$1$};
\end{tikzpicture}$$ Thus when $k=1$ we have a one-dimensional representation of the braid group $\B_n$.

\subsection{Dimensions of Level 2 representations}
Predictably, the dimension of the representation of $\B_n$ gets more complicated as we increase the level of the theory.  To motivate the general pattern, we work through the first few values of $n$ explicitly, labeling fusion trees with elements of $\mathcal{L}=\{0,1,2\}$.
\subsubsection{$n=2$}
For $n=2$, there are two admissible values of $i$ for the fusion tree, resulting in two one-dimensional representations.

$$
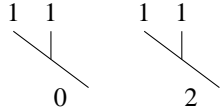
\begin{tikzpicture}
\draw (-3/2, 2)--(-3/2,13/8);
\draw (-2,2.25) node {$1$};
\draw (-3/2,2.25) node {$1$};
\draw (-2,2)--(-1,5/4);
\draw (-11/8, 1.125) node {$0$};
\end{tikzpicture}\hspace{15pt} \begin{tikzpicture}
\draw (-3/2, 2)--(-3/2,13/8);
\draw (-2,2.25) node {$1$};
\draw (-3/2,2.25) node {$1$};
\draw (-2,2)--(-1,5/4);
\draw (-11/8, 1.125) node {$2$};
\end{tikzpicture}$$

\subsubsection{$n=3$}
When $n=3$, the value of $i$ is determined, but there are two different ways to label the edges of the fusion tree consistently, giving a two-dimensional representation.
$$
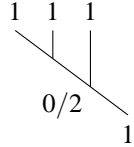
\begin{tikzpicture}
\node (first) {};
\draw (-3/2, 2)--(-3/2,13/8);
\draw (-1,2)--(-1, 10/8);
\draw (-2,2.25) node {$1$};
\draw (-3/2,2.25) node {$1$};
\draw (-1,2.25) node {$1$};
\draw (-2,2)--(-1/2,7/8);
\draw (-11/8, 1) node {$0/2$};
\draw (-1/2, 5/8) node {$1$};
\end{tikzpicture}$$

  \subsubsection{$n=4$}
When $n=4$, there are two distinct values of $i$, and for each value of $i$, two different ways to label the edges of the fusion tree. Therefore we get two separate two-dimensional representations.

$$
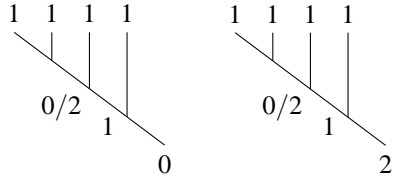
\begin{tikzpicture}
\draw (-1/2,2)--(-1/2,7/8);
\draw (-1/2,2.25) node {$1$};
\draw (-1,2)--(-1, 10/8);
\draw (-1,2.25) node {$1$};
\draw (-3/2, 2)--(-3/2,13/8);
\draw (-2,2.25) node {$1$};
\draw (-3/2,2.25) node {$1$};
\draw (-2,2)--(0,1/2);
\draw (-11/8, 1) node {$0/2$};
\draw (-3/4,3/4) node {$1$};
\draw (0,1/4) node {$0$};
\end{tikzpicture}
\hspace{15pt}
\begin{tikzpicture}
\draw (-1/2,2)--(-1/2,7/8);
\draw (-1/2,2.25) node {$1$};
\draw (-1,2)--(-1, 10/8);
\draw (-1,2.25) node {$1$};
\draw (-3/2, 2)--(-3/2,13/8);
\draw (-2,2.25) node {$1$};
\draw (-3/2,2.25) node {$1$};
\draw (-2,2)--(0,1/2);
\draw (-11/8, 1) node {$0/2$};
\draw (-3/4,3/4) node {$1$};
\draw (0,1/4) node {$2$};
\end{tikzpicture}$$

Both of these representations are isomorphic to $\mathds{C}^2$. The representation $\rho_{2,4,0}$ corresponding to the lefthand fusion tree is presented in the following section. The other representation $\rho_{2,4,2}$ is different, but similar, and is left to the reader as an exercise.

\subsubsection{The Majorana qubit}
We introduce a convenient piece of notation for fusion trees. Often we want to make the identification of a certain fusion tree corresponding to a two-dimensional representation with the standard orthonormal basis vectors $|0\rangle$ and $|1\rangle$ of $\mathds{C}^2$. To make this identification, we typically need to normalize a fusion tree. Instead of carrying around a potentially cumbersome normalization factor along with the fusion trees, we use open circles at the vertices of the fusion tree to indicate that it is normalized.
The usual notation for a qubit is as a superposition of the states $0$ and $1$, $\alpha| 0 \rangle + \beta |1\rangle$, $|\alpha|^2+|\beta|^2=1$. By identifying $|0\rangle$ and $|1\rangle$ with the fusion trees
$$
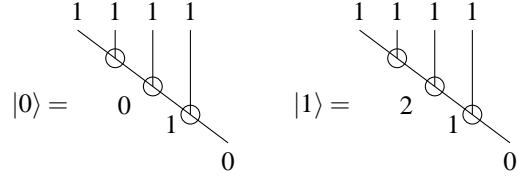
\begin{tikzpicture}
\node (first) {};
\draw (-2.5, 1) node {$|0\rangle=$};
\draw (-1/2,2)--(-1/2,7/8);
\draw (-1,10/8) circle (1/8);
\draw (-3/2,13/8) circle (1/8);
\draw(-1/2, 7/8) circle (1/8);
\draw (-1/2,2.25) node {$1$};
\draw (-1,2)--(-1, 10/8);
\draw (-1,2.25) node {$1$};
\draw (-3/2, 2)--(-3/2,13/8);
\draw (-2,2.25) node {$1$};
\draw (-3/2,2.25) node {$1$};
\draw (-2,2)--(0,1/2);
\draw (-11/8, 1) node {$0$};
\draw (-3/4,3/4) node {$1$};
\draw (0,1/4) node {$0$};
\node (second) [right=of first] {};
\begin{scope}[shift={(second)}, xshift=25mm]
\draw (-2.5, 1) node {$|1\rangle=$};
\draw (-1/2,2)--(-1/2,7/8);
\draw (-1,10/8) circle (1/8);
\draw (-3/2,13/8) circle (1/8);
\draw(-1/2, 7/8) circle (1/8);
\draw (-1/2,2.25) node {$1$};
\draw (-1,2)--(-1, 10/8);
\draw (-1,2.25) node {$1$};
\draw (-3/2, 2)--(-3/2,13/8);
\draw (-2,2.25) node {$1$};
\draw (-3/2,2.25) node {$1$};
\draw (-2,2)--(0,1/2);
\draw (-11/8, 1) node {$2$};
\draw (-3/4,3/4) node {$1$};
\draw (0,1/4) node {$0$};
\end{scope}
\end{tikzpicture}$$
we arrive at the famous \emph{Majorana qubit.}


This type of analysis produces the dimension of the Jones representation of $\B_n$ for $k=2$. In a similar way, we can analyze the dimensions of the Jones representation for the Fibonacci subtheory, by considering what happens when we label the top of the fusion tree by 2's.
\subsubsection{How the Fibonacci theory got its name}\label{fibthy}

Technically, the theory of Fibonacci anyons uses the \emph{colored Jones representation} where instead of considering $\text{Hom}(i,{1}^{\otimes n})$, we replace the label $1$ with another label $2=\tau$ in $\mathcal{L}=\{0,1,2,\ldots, k\}$ and consider $\text{Hom}(i,{a}^{\otimes n})$ for $a \in \mathcal{L}$. In particular, we are looking for a basis of $\text{Hom}(1, \tau^n)$, where $\tau$ is the Fibonacci anyon. This still provides a representation of the $n$-strand 
 group, we have just \lq\lq colored'' the braids by $\tau$.

\begin{rmk} It is possible to obtain the same representation through the uncolored Jones representation with the right choice of Kauffman variable $A$ up to a character because $1\otimes 3=2$.
\end{rmk}

The anyon model $\{1, \tau\}$ is called the Fibonacci theory because the Fibonacci numbers appear as the dimensions of the spaces $\text{Hom}(1,\tau \otimes \cdots \otimes \tau)$. Hereafter we use the $TLJ(A)$ labels and anyon labels interchangeably.

When $n=1$ there is one admissible fusion tree, but $i\ne 0$, and hence $\dim(V_{3,1,0})=0$.
$$\begin{tikzpicture}
\node (first) {};
\node (first) {};
\draw (0,0)--(0,1);
\draw (1/4, 1/8) node {$2$};
\draw (1/4,7/8) node {$2$};
\end{tikzpicture}$$

When $n=2$, there are two ways to label a fusion tree, one of which has trivial total charge, and hence $\dim(V_{3, 2, 0})=1$.
$$
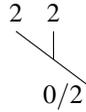
\begin{tikzpicture}
\draw (-3/2, 2)--(-3/2,13/8);
\draw (-2,2.25) node {$2$};
\draw (-3/2,2.25) node {$2$};
\draw (-2,2)--(-1,5/4);
\draw (-11/8, 1.125) node {$0/2$};
\end{tikzpicture}$$

When $n=3$, the image splits into a one-dimensional space isomorphic to $\mathds{C}$ and a two-dimensional space, isomorphic to $\mathds{C}^2$.
$$
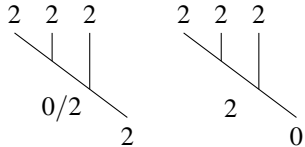
\begin{tikzpicture}
\node (first) {};
\draw (-3/2, 2)--(-3/2,13/8);
\draw (-1,2)--(-1, 10/8);
\draw (-2,2.25) node {$2$};
\draw (-3/2,2.25) node {$2$};
\draw (-1,2.25) node {$2$};
\draw (-2,2)--(-1/2,7/8);
\draw (-11/8, 1) node {$0/2$};
\draw (-1/2, 5/8) node {$2$};
\node (second) [right= of first, xshift=10mm] {};
\begin{scope}[shift={(second)}]
\draw (-3/2, 2)--(-3/2,13/8);
\draw (-1,2)--(-1, 10/8);
\draw (-2,2.25) node {$2$};
\draw (-3/2,2.25) node {$2$};
\draw (-1,2.25) node {$2$};
\draw (-2,2)--(-1/2,7/8);
\draw (-11/8, 1) node {$2$};
\draw (-1/2, 5/8) node {$0$};
\end{scope}
\end{tikzpicture}$$
Evidently $\dim(V_{3, 3, 0})=1$.

Now if $n=4$ and $i=0$, we get  $\dim(V_{3, 4, 0})=2$.
$$
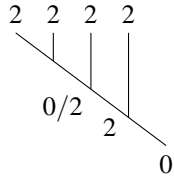
\begin{tikzpicture}
\node (first) {};
\draw (-1/2,2)--(-1/2,7/8);
\draw (-1/2,2.25) node {$2$};
\draw (-1,2)--(-1, 10/8);
\draw (-1,2.25) node {$2$};
\draw (-3/2, 2)--(-3/2,13/8);
\draw (-2,2.25) node {$2$};
\draw (-3/2,2.25) node {$2$};
\draw (-2,2)--(0,1/2);
\draw (-11/8, 1) node {$0/2$};
\draw (-3/4,3/4) node {$2$};
\draw (0,1/4) node {$0$};
\end{tikzpicture}$$

So far the dimensions form the sequence $0,1,1,2 \ldots$, the first few Fibonacci numbers. The Fibonacci fusion rule
$$ \tau \otimes \tau = 1 \oplus \tau$$
governs the dimension of $\text{Hom}(i,\tau^{\otimes n})$.

Using this fusion rule, we can make an observation about how our Fibonacci fusion trees are nested in one another.

$$
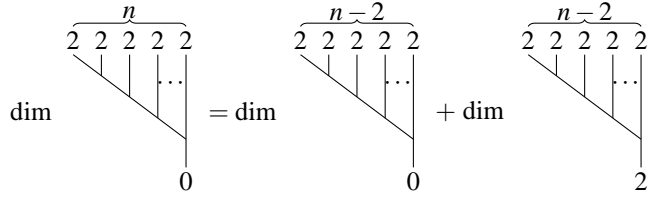
\begin{tikzpicture}[scale=.75]
\node (first) {};
\draw[decorate,decoration=brace] (-2,2.5)--(0,2.5);
\draw (-1,2.75) node {$n$};
\draw (-2.75,1) node {$\dim$};
\draw (0,2)--(0,0); 
\draw (-2,2)--(0,.5);
\draw (-3/2, 2)--(-3/2,13/8);
\draw (-1/2,2)--(-1/2,7/8);
\draw (-1,2)--(-1, 10/8);
\draw (-1/4,6/4) node {$\cdots$};
\draw (0,-1/4) node {$0$};
\draw (-2,2.25) node {$2$};
\draw (-3/2,2.25) node {$2$};
\draw (-1, 2.25) node {$2$};
\draw (-1/2,2.25) node {$2$};
\draw (0,2.25) node {$2$};
\draw (1,1) node {$= \dim$};
\end{tikzpicture}
\begin{tikzpicture}[scale=.75]
\draw[decorate,decoration=brace] (-2,2.5)--(0,2.5);
\draw (-1,2.75) node {$n-2$};
\draw (0,2)--(0,0); 
\draw (-2,2)--(0,.5);
\draw (-3/2, 2)--(-3/2,13/8);
\draw (-1/2,2)--(-1/2,7/8);
\draw (-1,2)--(-1, 10/8);
\draw (-1/4,6/4) node {$\cdots$};
\draw (0,-1/4) node {$0$};
\draw (-2,2.25) node {$2$};
\draw (-3/2,2.25) node {$2$};
\draw (-1, 2.25) node {$2$};
\draw (-1/2,2.25) node {$2$};
\draw (0,2.25) node {$2$};
\draw (1,1) node {$+$  $\dim $};
\end{tikzpicture}
\begin{tikzpicture}[scale=.75]
\draw[decorate,decoration=brace] (-2,2.5)--(0,2.5);
\draw (-1,2.75) node {$n-2$};
\draw (0,2)--(0,0); 
\draw (-2,2)--(0,.5);
\draw (-3/2, 2)--(-3/2,13/8);
\draw (-1/2,2)--(-1/2,7/8);
\draw (-1,2)--(-1, 10/8);
\draw (-1/4,6/4) node {$\cdots$};
\draw (0,-1/4) node {$2$};
\draw (-2,2.25) node {$2$};
\draw (-3/2,2.25) node {$2$};
\draw (-1, 2.25) node {$2$};
\draw (-1/2,2.25) node {$2$};
\draw (0,2.25) node {$2$};
\end{tikzpicture}$$

So the dimensions satisfy the recurrence relation
$$f_{n,0} = f_{n-2,0} + f_{n-2,2} = f_{n-2,0} + f_{n-1,0}$$
which is exactly the relation that defines the Fibonacci numbers $F_n$ and hence
$\dim V_{3,\tau^{\otimes n},0} = F_{n-1}$.

Algebraizing this fusion rule, we get the equation $x^2 = 1 + x$, whose solutions are the golden ratio $\phi$ and its Galois conjugate.
The golden ratio also satisfies the identity $\phi=\phi^{-1}+1$, which will be useful for calculations in the Fibonacci theory. The Fibonacci representation always splits into two subrepresentations as $\rho_{3, \tau^{\otimes n}} = \rho_{3,\tau^{\otimes n}, 1} \oplus \rho_{3,\tau^{\otimes n}, \tau}: \B_n \to U(F_{n-1}) \oplus U(F_n)$, corresponding to the total charge $1$ and total charge $\tau$.

\subsubsection{The Fibonacci qubit}
For the Fibonacci qubit, we identify

$$
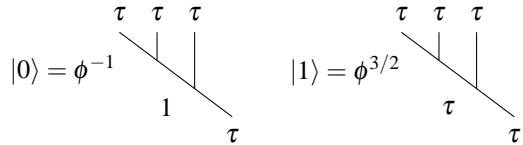
\begin{tikzpicture}
\node (first) {};
\draw (-2.75,1.5) node {$|0\rangle =\phi^{-1}$};
\draw (-3/2, 2)--(-3/2,13/8);
\draw (-1,2)--(-1, 10/8);
\draw (-2,2.25) node {$\tau$};
\draw (-3/2,2.25) node {$\tau$};
\draw (-1,2.25) node {$\tau$};
\draw (-2,2)--(-1/2,7/8);
\draw (-11/8, 1) node {$1$};
\draw (-1/2, 5/8) node {$\tau$};
\node (second) [right= of first, xshift=25mm] {};
\begin{scope}[shift={(second)}]
\draw (-2.75,1.5) node {$|1\rangle =\phi^{3/2}$};
\draw (-3/2, 2)--(-3/2,13/8);
\draw (-1,2)--(-1, 10/8);
\draw (-2,2.25) node {$\tau$};
\draw (-3/2,2.25) node {$\tau$};
\draw (-1,2.25) node {$\tau$};
\draw (-2,2)--(-1/2,7/8);
\draw (-11/8, 1) node {$\tau$};
\draw (-1/2, 5/8) node {$\tau$};
\end{scope}
\end{tikzpicture}$$
or in the new notation

$$
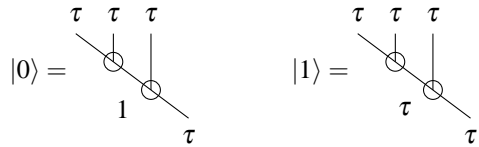
\begin{tikzpicture}
\node (first) {};
\draw (-2.5,1.5) node {$|0\rangle =$};
\draw (-3/2, 2)--(-3/2,13/8);
\draw (-1,2)--(-1, 10/8);
\draw (-1,10/8) circle (1/8);
\draw (-2,2.25) node {$\tau$};
\draw (-3/2,13/8) circle (1/8);
\draw (-3/2,2.25) node {$\tau$};
\draw (-1,2.25) node {$\tau$};
\draw (-2,2)--(-1/2,7/8);
\draw (-11/8, 1) node {$1$};
\draw (-1/2, 5/8) node {$\tau$};
\node (second) [right= of first, xshift=25mm] {};
\begin{scope}[shift={(second)}]
\draw (-2.5,1.5) node {$|1\rangle =$};
\draw (-3/2, 2)--(-3/2,13/8);
\draw (-3/2,13/8) circle (1/8);
\draw (-1,2)--(-1, 10/8);
\draw (-1,10/8) circle (1/8);
\draw (-2,2.25) node {$\tau$};
\draw (-3/2,2.25) node {$\tau$};
\draw (-1,2.25) node {$\tau$};
\draw (-2,2)--(-1/2,7/8);
\draw (-11/8, 1) node {$\tau$};
\draw (-1/2, 5/8) node {$\tau$};
\end{scope}
\end{tikzpicture}$$

\subsubsection{Dense versus sparse qubit encodings}
The qubit encoding above using three Fibonacci anyons is called a \emph{dense encoding.} By raising the number of anyons, like in the two-dimensional representation
 $$
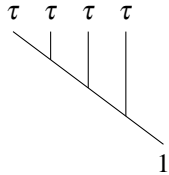
\begin{tikzpicture}
\node (first) {};
\draw (-1/2,2)--(-1/2,7/8);
\draw (-1/2,2.25) node {$\tau$};
\draw (-1,2)--(-1, 10/8);
\draw (-1,2.25) node {$\tau$};
\draw (-3/2, 2)--(-3/2,13/8);
\draw (-2,2.25) node {$\tau$};
\draw (-3/2,2.25) node {$\tau$};
\draw (-2,2)--(0,1/2);
\draw (0,1/4) node {$1$};
\end{tikzpicture}$$
we obtain a \emph{sparse} encoding. While the dense encoding is mathematically easier to work with, the sparse encoding is physically preferable. This is because the total charge $i$ of an anyon system is a boundary condition, and in an experimental set up, letting the boundary condition just correspond to the ground state is energetically more favorable.

Now that we have found two-dimensional representations $\rho_{2,4,0/2}$ and $\rho_{3,\tau^{\otimes 3}, \tau}$, we would like to be able to compute these them explicitly and write down their matrices with respect to an orthonormal basis on the vector spaces $V_{k,n,i}$.

\subsection{Computing Jones Representations and $\rho_{2,4,0}(\sigma_1)$}
To illustrate a general method for computing the Jones representation of the generators $\sigma_i$ of the braid group $\B_n$, we calculate the Jones representations $\rho_{2,4,0}(\sigma_1)$.

Recall the two fusion trees that span $V_{2,4,0}$, shown bellow, which we now call $\tilde{e}_0$ and $\tilde{e}_1$. The first step is to turn them into a basis using Gram-Schmidt orthonormalization.

$$
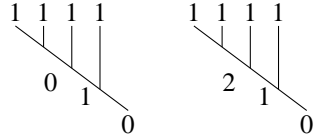
\begin{tikzpicture}[scale=.75]
\node (first) {};
\draw (-1/2,2)--(-1/2,7/8);
\draw (-1/2,2.25) node {$1$};
\draw (-1,2)--(-1, 10/8);
\draw (-1,2.25) node {$1$};
\draw (-3/2, 2)--(-3/2,13/8);
\draw (-2,2.25) node {$1$};
\draw (-3/2,2.25) node {$1$};
\draw (-2,2)--(0,1/2);
\draw (-11/8, 1) node {$0$};
\draw (-3/4,3/4) node {$1$};
\draw (0,1/4) node {$0$};
\node (second) [right=of first] {};
\begin{scope}[shift={(second)}, xshift=15mm]
\draw (-1/2,2)--(-1/2,7/8);
\draw (-1/2,2.25) node {$1$};
\draw (-1,2)--(-1, 10/8);
\draw (-1,2.25) node {$1$};
\draw (-3/2, 2)--(-3/2,13/8);
\draw (-2,2.25) node {$1$};
\draw (-3/2,2.25) node {$1$};
\draw (-2,2)--(0,1/2);
\draw (-11/8, 1) node {$2$};
\draw (-3/4,3/4) node {$1$};
\draw (0,1/4) node {$0$};
\end{scope}
\end{tikzpicture}$$
\subsubsection{Notation}
It will be convenient to introduce another notation for elements of $\text{Hom}(i, 1^{\otimes n})$, in which the need to label every edge of a fusion tree is eliminated. Edges labeled by the ground state $0$ become dashed edges, edges labeled by a $1$ are usual lines, and edges labeled by a $2$ become wavy lines.

$$ 
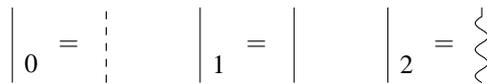
\begin{tikzpicture}
\node (first) {};
\draw (0,0)--(0,1);
\draw (1/4, 1/4) node {$0$};
\draw (3/4,1/2) node {$=$};
\node (second) [right=of first] {};
\begin{scope}[shift={(second)}]
\draw[dashed] (0,0)--(0,1);
\end{scope}
\node (third) [right=of second] {};
\begin{scope}[shift={(third)}]
\draw (0,0)--(0,1);
\draw (1/4, 1/4) node {$1$};
\draw (3/4,1/2) node {$=$};
\end{scope}
\node (fourth) [right=of third] {};
\begin{scope}[shift={(fourth)}]
\draw (0,0)--(0,1);
\end{scope}
\node (fifth) [right=of fourth] {};
\begin{scope}[shift={(fifth)}]
\draw (0,0)--(0,1);
\draw (1/4, 1/4) node {$2$};
\draw (3/4,1/2) node {$=$};
\end{scope}
\node (sixth) [right=of fifth] {};
\begin{scope}[shift={(sixth)}]
\draw[decorate,decoration={coil, aspect=0}] (0,0)--(0,1);
\end{scope}
\end{tikzpicture}$$

Under this new notation $\tilde{e}_0$ and $\tilde{e}_1$ become

$$
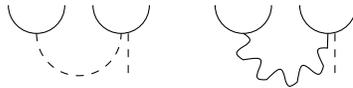
\begin{tikzpicture}[scale=.75]
\node (first) {};
\draw (1,1) arc (-180:0:1/2);
\draw (2.5,1) arc (-180:0:1/2);
\draw[dashed] (3/2,1/2) arc (-180:0:3/4);
\draw[dashed] (3.125, 1/2)--(3.125,-1/4);
\node (second)[right= of first] {};
\begin{scope}[shift={(second)}, xshift=2cm]
\draw (1,1) arc (-180:0:1/2);
\draw (2.5,1) arc (-180:0:1/2);
\draw[decorate,decoration={coil, aspect=0}] (3/2,1/2) arc (-180:0:3/4);
\draw[dashed] (3.125, 1/2)--(3.125,-1/4);
\end{scope}
\end{tikzpicture}$$

Hereafter we will drop the dashed lines labeling the ground state. First we calculate the inner products $\langle \tilde{e}_i, \tilde{e}_j\rangle$ using the graphical calculus.

$\begin{tikzpicture}[scale=.75]
\node (first) {$\langle \tilde{e}_0, \tilde{e}_0 \rangle=$};
\node (second) [xshift=0mm] {};
\begin{scope}[shift={(second)}]
\draw (1,1) arc (-180:0:1/2);
\draw (2.5,1) arc (-180:0:1/2);
\draw (3.5,1) to [out=90, in=90] (4,1)--(4,-1) to [out=-90, in=-90] (3.5,-1);
\draw (2.5,1) to [out=90, in=90] (4.5,1)--(4.5,-1) to [out=-90, in=-90] (2.5,-1);
\draw (2,1) to [out=90, in=90] (5,1)--(5,-1) to [out=-90, in=-90] (2,-1);
\draw (1,1) to [out=90, in=90] (5.5,1)--(5.5,-1) to [out=-90, in=-90] (1,-1);
\draw (1,-1) arc (180:0:1/2);
\draw (2.5,-1) arc (180:0:1/2);
\end{scope}
\node (third) [right=of second,xshift=30mm] {$=d^2=2$};
\end{tikzpicture}$

$\begin{tikzpicture}[scale=.75]
\node (first) {$\langle \tilde{e}_1, \tilde{e}_1 \rangle=$};
\node (second) [xshift=0cm] {};
\begin{scope}[shift={(second)}]
\draw (1,1) arc (-180:0:1/2);
\draw (2.5,1) arc (-180:0:1/2);
\draw[decorate,decoration={coil, aspect=0}] (3/2,1/2) arc (-180:0:3/4);
\draw[decorate,decoration={coil, aspect=0}] (3,-3/2) arc (0:180:3/4);
\draw (3.5,1) to [out=90, in=90] (4,1)--(4,-2) to [out=-90, in=-90] (3.5,-2);
\draw (2.5,1) to [out=90, in=90] (4.5,1)--(4.5,-2) to [out=-90, in=-90] (2.5,-2);
\draw (2,1) to [out=90, in=90] (5,1)--(5,-2) to [out=-90, in=-90] (2,-2);
\draw (1,1) to [out=90, in=90] (5.5,1)--(5.5,-2) to [out=-90, in=-90] (1,-2);
\draw (1,-2) arc (180:0:1/2);
\draw (2.5,-2) arc (180:0:1/2);
\end{scope}
\node (third) [right=of second,xshift=4.5cm] {};
\begin{scope}[shift={(third)}]
\draw (-3/2,0) node {$=$};
\draw (0,0) circle (1) {};
\draw[decorate,decoration={coil, aspect=0}] (4/5,3/5)--(4,3/5);
\draw[decorate,decoration={coil, aspect=0}] (4/5,-3/5)--(4,-3/5);
\draw (4.8,0) circle (1) {};
\end{scope}
\node (fourth) [right=of third, xshift=4.5cm] {};
\begin{scope}[shift={(fourth)}]
\draw (-3/2,0) node {$=$};
\draw (0,0) circle (1) {};
\draw (4/5,3/5)--(4,3/5);
\draw (4/5,-3/5)--(4,-3/5);
\draw (4.8,0) circle (1) {};
\draw (5/4, 1) node {2};
\draw (3/4, 0) node {1};
\draw (1/2, 5/4) node {1};
\draw (6.5,0) node {$=1$};
\end{scope}
\end{tikzpicture}$

$\begin{tikzpicture}[scale=.75]
\node (first) {$\langle \tilde{e}_0, \tilde{e}_1 \rangle=$};
\node (second) [xshift=1cm] {};
\begin{scope}[shift={(second)}]
\draw (1,1) arc (-180:0:1/2);
\draw (2.5,1) arc (-180:0:1/2);
\draw[decorate,decoration={coil, aspect=0}] (3/2,1/2) arc (-180:0:3/4);
\draw (3.5,1) to [out=90, in=90] (4,1)--(4,-2) to [out=-90, in=-90] (3.5,-2);
\draw (2.5,1) to [out=90, in=90] (4.5,1)--(4.5,-2) to [out=-90, in=-90] (2.5,-2);
\draw (2,1) to [out=90, in=90] (0,1)--(0,-2) to [out=-90, in=-90] (2,-2);
\draw (1,1) to [out=90, in=90] (.5,1)--(.5,-2) to [out=-90, in=-90] (1,-2);
\draw (1,-2) arc (180:0:1/2);
\draw (2.5,-2) arc (180:0:1/2);
\end{scope}
\node (third) [right=of second,xshift=3.5cm] {};
\begin{scope}[shift={(third)}]
\draw (-3/2,0) node {$=$};
\draw (0,0) circle (1) {};
\draw[decorate,decoration={coil, aspect=0}] (1,0)--(3.8,0);
\draw (4.8,0) circle (1) {};
\end{scope}
\node (fourth)[right=of third,xshift=3.75cm] {};
\begin{scope}[shift={(fourth)}]
\draw (-1/4,0) node {$= 0$};
\end{scope}
\end{tikzpicture}$

The full details of the calculations for $\langle \tilde{e}_1,\tilde{e}_1\rangle$ and $\langle \tilde{e}_0,\tilde{e}_1 \rangle$ are left to the reader as an exercise, and entails inserting the appropriate Jones-Wenzl projectors at the trivalent vertices and using the graphical calculus.

Since $\langle \tilde{e}_0, \tilde{e}_1\rangle=0$, the choices $e_0 = \frac{1}{\sqrt{2}} \tilde{e}_0$ and $e_1 = \tilde{e}_1$ define an orthonormal basis $\{e_0, e_1\}$ of $V_{2,4,0}$. Then with respect to this basis, $\rho(\sigma_1) = \begin{pmatrix} \langle  e_0,\sigma_1 e_0 \rangle & \langle  e_0, \sigma_1 e_1 \rangle \\ \langle e_1,  \sigma_1 e_0 \rangle & \langle e_1,\sigma_1  e_1 \rangle \end{pmatrix}$.

For example,
$$\begin{tikzpicture}[scale=.5]
\node (first) {};
\draw (-3,-1.5) node {$\langle \sigma_1 e_0, e_0 \rangle = (\frac{1}{\sqrt{2}})^2$};
\braid[number of strands={2}] a_1^-1;
\draw (1,-1.5) arc (-180:0:1/2);
\draw (2.5,-1.5) arc (-180:0:1/2);
\draw (2.5,0)--(2.5,-1.5);
\draw (3.5,0)--(3.5,-1.5);

\draw (1,-3) arc (180:0:1/2);
\draw (1,-3) to [out=-90, in=-90] (.5,-3)--(.5,0) to [out=90, in=90] (1,0);
\draw (2,-3) to [out=-90, in=-90] (0,-3)--(0,0) to [out=90, in=90] (2,0);

\draw (2.5,-3) arc (180:0:1/2);
\draw (3.5,-3) to [out=-90, in=-90] (4,-3)--(4,0) to [out=90, in=90] (3.5,0);
\draw (2.5,-3) to [out=-90, in=-90] (4.5,-3)--(4.5,0) to [out=90, in=90] (2.5,0);
\node (second) [right=of first, xshift=2.25cm] {};
\begin{scope}[shift={(second)}]
\draw (-1.25,-1.5) node {$=\frac{1}{2}\cdot  A$};

\draw (1,-1.5) arc (-180:0:1/2);
\draw (1,0)--(1,-1.5);
\draw (2,0)--(2,-1.5);

\draw (2.5,-1.5) arc (-180:0:1/2);
\draw (2.5,0)--(2.5,-1.5);
\draw (3.5,0)--(3.5,-1.5);

\draw (1,-3) arc (180:0:1/2);
\draw (1,-3) to [out=-90, in=-90] (.5,-3)--(.5,0) to [out=90, in=90] (1,0);
\draw (2,-3) to [out=-90, in=-90] (0,-3)--(0,0) to [out=90, in=90] (2,0);

\draw (2.5,-3) arc (180:0:1/2);
\draw (3.5,-3) to [out=-90, in=-90] (4,-3)--(4,0) to [out=90, in=90] (3.5,0);
\draw (2.5,-3) to [out=-90, in=-90] (4.5,-3)--(4.5,0) to [out=90, in=90] (2.5,0);
\end{scope}
\node (third) [right=of second, xshift=2.5cm] {};
\begin{scope}[shift={(third)}]
\draw (-1.5,-1.5) node {$+ \frac{1}{2}\cdot  A^{-1}$};

\draw (1,-1.5) arc (-180:0:1/2);
\draw (1,-1.5) arc (180:0:1/2);

\draw (2.5,-1.5) arc (-180:0:1/2);
\draw (1,0) arc (-180:0:1/2);

\draw (2.5,0)--(2.5,-1.5);
\draw (3.5,0)--(3.5,-1.5);

\draw (1,-3) arc (180:0:1/2);
\draw (1,-3) to [out=-90, in=-90] (.5,-3)--(.5,0) to [out=90, in=90] (1,0);
\draw (2,-3) to [out=-90, in=-90] (0,-3)--(0,0) to [out=90, in=90] (2,0);

\draw (2.5,-3) arc (180:0:1/2);
\draw (3.5,-3) to [out=-90, in=-90] (4,-3)--(4,0) to [out=90, in=90] (3.5,0);
\draw (2.5,-3) to [out=-90, in=-90] (4.5,-3)--(4.5,0) to [out=90, in=90] (2.5,0);
\draw (9,-1.5) node {$= \frac{1}{2}(A \cdot d^2 + A^{-1} d^3) = -A^{-3}$};
\end{scope}
\end{tikzpicture},$$
where the crossings were resolved using the Kauffman bracket. Similar calculations for the remaining matrix entries show that $$\rho(\sigma_1) = \begin{pmatrix} -A^{-3} & 0 \\ 0 & A \end{pmatrix}.$$ By repeating the same method to find the remaining generators $\rho(\sigma_2)$ and $\rho(\sigma_3)$, one can calculate the image $\rho_{4,2,0}(b)$ for any $b\in \B_4$.

This outlines an elementary way to find the Jones representation. While it has the benefit that it uses only knowledge of the Kauffman bracket and arithmetic, as $n$ gets larger it becomes inefficient to do by hand. Using additional structure in $TLJ(A)$ as a UMC, the \emph{$\theta$-symbols}, \emph{$R$-symbols}, and \emph{$F$-symbols} provide more tools to find $\rho_{n,k,i}$ using the graphical calculus.

\subsection{$\theta$-symbols, $R$-symbols, and $F$-symbols}
In this section we will use graphs with open circles at the vertices to indicate that it has been normalized.

Take any admissibly-labeled trivalent vertex $e^{ab}_c$, i.e. an element of $\text{Hom}(c,a\otimes b)$. Then the $\theta$-symbol $\theta(a,b,c)$ is defined to be the the inner product $\langle e^{ab}_c,e^{ab}_c\rangle$ of the trivalent tree basis $e^{ab}_c$.

In the next two sections we introduce the $R$ and $F$-symbols. Once these symbols are determined for an anyon model they can be used to calculate the desired braid group representation.

\subsubsection{R-matrices}
Braiding gives a linear map
$$\begin{tikzpicture}[scale=.5]
\braid[number of strands={2}] a_1^-1;
\draw (1,-1.5) to [out=-90, in=-90] (2,-1.5);
\draw (1.5,-1.8)--(1.5,-2.5);
\draw (7,-1) node {$ :  \text{Hom}(c, a \otimes b) \to \text{Hom}(c, b \otimes a)$};
 \end{tikzpicture}.$$ It is in fact an isomorphism, with inverse given by the opposite crossing. Since $\text{Hom}(c, b \otimes a)$ is one-dimensional in TLJ theory, and we already have a preferred basis for it, namely the trivalent vertex labeled by $b$, $a$, and $c$, the equation

$$\begin{tikzpicture}[scale=.5]
\node (one) {};
\braid[number of strands={2}] a_1^-1;
\draw (1,-1.5) to [out=-90, in=-90] (2,-1.5);
\draw (1,.5) node {$a$};
\draw (2,.5) node {$b$};
\draw (1.5,-1.8)--(1.5,-2.5);
\draw (1.5,-1.8) circle (1/4);
\draw (1.5, -3) node {$c$};
\draw (3,-1) node {$=R_c^{ab}$};
\node[shape=coordinate, right=of one, xshift=1.75cm, yshift=-.5cm] (trivalent) {};
\begin{scope}[shift={(trivalent)}]
\draw (trivalent) -- node[auto] {$c$} (-90:2cm);
\draw (trivalent) -- node[auto,swap] {$b$} (45:2cm);
\draw (trivalent) -- node[auto] {$a$} (135:2cm);
\draw (0,0) circle (1/4);
\end{scope}
\end{tikzpicture}$$

holds, where $R_c^{ba}$ is a scalar, which we call the \emph{braiding eigenvalue} or $R$-symbol.
There is a general formula that gives the $R$-symbols, for any Kaufmann variable $A$ and any Temperley-Lieb category, given by

$$R_c^{ab} = (-1)^{\frac{a+b-c}{2}}A^{\frac{-[a(a+2) + b(b+2)-c(c+2)]}{2}}.$$

One can also calculate the $R$-symbols from their defining relation by taking the inner product of both sides of the equation with the trivalent vertex labeled by $a$, $b$, and $c$, which we demonstrate in the following example with $R^{22}_2$, where the labels ``2" on the right hand side of the equation indicate that every edge in the diagram is labeled by a 2.

$$\begin{tikzpicture}[scale=.5]
\node (top) {};
\draw (3,-1) node {2};
\draw (3,-5) node {2};
\braid[number of strands={2}] a_1^-1;
\draw (1,-1.5) arc (-180:0:1/2);
\draw (1.5,.5) -- (1.5,1);
\draw (1,0) -- (1.5, .5);
\draw (2,0) -- (1.5, .5);
\draw (1.5,-2)--(1.5,-2.5);
\draw (1.5,1) to [out=90, in =90] (2.5,1) -- (2.5,-2.5) to [out=-90, in=-90] (1.5,-2.5);
\draw (0,-3)--(4,-3);
\draw (-1.5,-3) node {$R^{22}_2=$};
\node (bottom) [below=of top, yshift=-1cm] {};
\begin{scope}[shift={(bottom)}]
\draw (1.5,0) -- (1.5,1);
\draw (1,-.75) -- (1.5, 0);
\draw (2,-.75) -- (1.5, 0);
\draw (1.5,-1.5)--(1.5,-2.5);
\draw (1.5,-1.5)--(2,-.75);
\draw (1.5,-1.5)--(1, -.75);
\draw (1.5,1) to [out=90, in =90] (2.5,1) -- (2.5,-2.5) to [out=-90, in=-90] (1.5,-2.5);
\end{scope}
\end{tikzpicture}$$

By resolving the top diagram via the Kauffman bracket and inserting Jones-Wenzl projectors $p_2$ at the appropriate vertices, $R^{22}_2$ can be calculated by hand. However, doing so efficiently requires familiarity with the graphical calculus and its shortcuts.

\begin{ex} Calculate $R^{22}_2$ using the graphical calculus and compare with the formula.
\end{ex}

\subsubsection{F-matrices}
The other hero in the theory is the $F$-matrices, which tells us how to make changes of basis. We consider the space $\text{Hom}(d, (a \otimes b) \otimes c)$, with orthonormal basis given by $\{e^{(ab)c}_{d,m}\}$, where

$$\begin{tikzpicture}
\node (first) {};
\draw (-2.5,1.5) node {$e^{(ab)c}_{d,m} =$};
\draw (-3/2, 2)--(-3/2,13/8);
\draw (-1,2)--(-1, 10/8);
\draw (-1,10/8) circle (1/8);
\draw (-2,2.25) node {$a$};
\draw (-3/2,13/8) circle (1/8);
\draw (-3/2,2.25) node {$b$};
\draw (-1,2.25) node {$c$};
\draw (-2,2)--(-1/2,7/8);
\draw (-11/8, 1) node {$m$};
\draw (-1/2, 5/8) node {$d$};
\end{tikzpicture}.$$

Similarly, we can consider the orthonormal basis $\{e^{a(bc)}_{d,n}\}$ of $\text{Hom}(d, a\otimes (b \otimes c))$, where

$$\begin{tikzpicture}
\draw (-1,1.5) node {$e^{a(bc)}_{d,n} =$};
\raisebox{\depth}{\scalebox{-1}[1]{\draw (-2/2, 2)--(-2/2,13/8);
\draw (-1/2,2)--(-1/2, 10/8);
\draw (-1/2,10/8) circle (1/8);
\draw (-2/2,13/8) circle (1/8);
\draw (-3/2,2)--(0,7/8);
\draw (-11/8, 1) node {};}}
\draw (1/2,2.25) node {$a$};
\draw (2/2,2.25) node {$b$};
\draw (3/2,2.25) node {$c$};
\draw (0, 5/8) node {$d$};
\draw (7/8, 1) node {$n$};
\end{tikzpicture}.$$

Then $F: \{e^{(ab)c}_{d,m}\} \to \{e^{a(bc)}_{d,m}\}$ is the change of basis matrix, satisfying

$$e^{(ab)c}_{d,m}= \sum F^{abc}_{d,nm}e^{a(bc)}_{d,n}.$$

The $F$-matrices are notoriously hard to find, although there is a general formula. Similar to the method described in the previous section for computing the $R$-symbols by tracing out their defining relation, there is a way to calculate them using the graphical calculus.

For example, the two fusion trees \begin{tikzpicture}[baseline=0] \node[shape=coordinate, xshift=-2.5mm] (trivalent) {};
\begin{scope}[shift={(trivalent)}]
\draw (1/4,1/8) node[draw=none] {\tiny $0/2$};
\draw (trivalent) -- node[draw=none] {} (0:.5cm);
\draw (trivalent) -- node[draw=none] {} (120:.5cm);
\draw (trivalent) -- node[draw=none] {} (240:.5cm);
\end{scope}
\node[shape=coordinate, xshift=2.5mm] (trivalent2) {};
\begin{scope}[shift={(trivalent2)}]

\draw (trivalent2) -- node[draw=none] {} (60:.5cm);
\draw (trivalent2) -- node[draw=none] {} (-60:.5cm);
\end{scope}\end{tikzpicture}
and
\begin{tikzpicture}[baseline=0]
\node[shape=coordinate, yshift=2.5mm] (trivalent) {};
\begin{scope}[shift={(trivalent)}]
\draw (1/4,-1/4) node[draw=none] {\tiny $0/2$};
\draw (trivalent) -- node[draw=none] {} (30:.5cm);
\draw (trivalent) -- node[draw=none] {} (150:.5cm);
\draw (trivalent) -- node[draw=none] {} (-90:.5cm);
\end{scope}
\node[shape=coordinate, xshift=0mm,yshift=-.25cm] (trivalent2) {};
\begin{scope}[shift={(trivalent2)}]
\draw (trivalent2) -- node[draw=none] {} (-30:.5cm);
\draw (trivalent2) -- node[draw=none] {} (-150:.5cm);
\end{scope}
\end{tikzpicture} each give a basis for Hom$(2 \otimes 2, 2 \otimes 2)$, where each edge that is not explicitly labeled is understood to be labeled with a 2. Therefore they satisfy the following equations.

$$\begin{tikzpicture}
\matrix[nodes={draw},
        row sep=0.3cm,column sep=0.5cm] {

\node[shape=coordinate, xshift=-2.5mm] (trivalent) {};
\begin{scope}[shift={(trivalent)}]
\draw (0,0) circle (1/10);
\draw[dashed] (trivalent) -- node[draw=none] {} (0:.5cm);
\draw (trivalent) -- node[draw=none] {} (120:.5cm);
\draw (trivalent) -- node[draw=none] {} (240:.5cm);
\end{scope}
\node[shape=coordinate, xshift=2.5mm] (trivalent2) {};
\begin{scope}[shift={(trivalent2)}]
\draw (0,0) circle (1/10);
\draw (trivalent2) -- node[draw=none] {} (60:.5cm);
\draw (trivalent2) -- node[draw=none] {} (-60:.5cm);
\end{scope}
&
\draw (0,0) node[draw=none] {$=$};
&
\draw (-1,0) node[draw=none] {$\alpha$};
\node[shape=coordinate, yshift=2.5mm] (trivalent) {};
\begin{scope}[shift={(trivalent)}]
\draw (0,0) circle (1/10);
\draw (trivalent) -- node[draw=none] {} (30:.5cm);
\draw (trivalent) -- node[draw=none] {} (150:.5cm);
\draw[dashed] (trivalent) -- node[draw=none] {} (-90:.5cm);
\end{scope}
\node[shape=coordinate, yshift=-2.5mm] (trivalent2) {};
\begin{scope}[shift={(trivalent2)}]
\draw (0,0) circle (1/10);
\draw (trivalent2) -- node[draw=none] {} (-30:.5cm);
\draw (trivalent2) -- node[draw=none] {} (-150:.5cm);
\end{scope}
&
\draw (0,0) node[draw=none] {$+$};
&
\draw (-1,0) node[draw=none] {$\beta$};
&
\node[shape=coordinate, yshift=2.5mm] (trivalent) {};
\begin{scope}[shift={(trivalent)}]
\draw (0,0) circle (1/10);
\draw (trivalent) -- node[draw=none] {} (30:.5cm);
\draw (trivalent) -- node[draw=none] {} (150:.5cm);
\draw (trivalent) -- node[draw=none] {} (-90:.5cm);
\end{scope}
\begin{scope}[shift={(trivalent2)}]
\draw (0,0) circle (1/10);
\draw (trivalent2) -- node[draw=none] {} (-30:.5cm);
\draw (trivalent2) -- node[draw=none] {} (-150:.5cm);
\end{scope}
\\};
\end{tikzpicture}$$

$$\begin{tikzpicture}
\matrix[nodes={draw},
        row sep=0.3cm,column sep=0.5cm] {

\node[shape=coordinate, xshift=-2.5mm] (trivalent) {};
\begin{scope}[shift={(trivalent)}]
\draw (0,0) circle (1/10);
\draw (trivalent) -- node[draw=none] {} (0:.5cm);
\draw (trivalent) -- node[draw=none] {} (120:.5cm);
\draw (trivalent) -- node[draw=none] {} (240:.5cm);
\end{scope}
\node[shape=coordinate, xshift=2.5mm] (trivalent2) {};
\begin{scope}[shift={(trivalent2)}]
\draw (0,0) circle (1/10);
\draw (trivalent2) -- node[draw=none] {} (60:.5cm);
\draw (trivalent2) -- node[draw=none] {} (-60:.5cm);
\end{scope}
&
\draw (0,0) node[draw=none] {$=$};
&
\draw (-1,0) node[draw=none] {$\gamma$};
\node[shape=coordinate, yshift=2.5mm] (trivalent) {};
\begin{scope}[shift={(trivalent)}]
\draw (0,0) circle (1/10);
\draw (trivalent) -- node[draw=none] {} (30:.5cm);
\draw (trivalent) -- node[draw=none] {} (150:.5cm);
\draw[dashed] (trivalent) -- node[draw=none] {} (-90:.5cm);
\end{scope}
\node[shape=coordinate, yshift=-2.5mm] (trivalent2) {};
\begin{scope}[shift={(trivalent2)}]
\draw (0,0) circle (1/10);
\draw (trivalent2) -- node[draw=none] {} (-30:.5cm);
\draw (trivalent2) -- node[draw=none] {} (-150:.5cm);
\end{scope}
&
\draw (0,0) node[draw=none] {$+$};
&
\draw (-1,0) node[draw=none] {$\delta$};
&
\node[shape=coordinate, yshift=2.5mm] (trivalent) {};
\begin{scope}[shift={(trivalent)}]
\draw (0,0) circle (1/10);
\draw (trivalent) -- node[draw=none] {} (30:.5cm);
\draw (trivalent) -- node[draw=none] {} (150:.5cm);
\draw (trivalent) -- node[draw=none] {} (-90:.5cm);
\end{scope}
\begin{scope}[shift={(trivalent2)}]
\draw (0,0) circle (1/10);
\draw (trivalent2) -- node[draw=none] {} (-30:.5cm);
\draw (trivalent2) -- node[draw=none] {} (-150:.5cm);
\end{scope}
\\};
\end{tikzpicture}$$

The constants can be determined by taking the trace of each equation in two different ways: once vertically and once horizontally. That is, once connecting top and bottom edges and once connecting left and right edges. The calculations can be simplified using the $\theta$-symbols. The $F$-matrix is then given by $\begin{pmatrix} \alpha & \beta \\ \gamma & \delta \end{pmatrix}$.

\begin{ex} Use the graphical calculus to show that the $F$-matrix is given by $\begin{pmatrix} \phi^{-1} & \phi^{-1/2} \\ \phi^{-1/2} & -\phi^{-1} \end{pmatrix}$.
\end{ex}

\subsection{Calculating the representation $\rho_{3,\tau^{\otimes 3}, \tau}$ for the Fibonacci theory}

\subsubsection{The Fibonacci $R$ and $F$-matrix}
When $A=\pm ie^{2\pi i/10}$ using the formula $R_c^{ab} = (-1)^{\frac{a+b-c}{2}}A^{\frac{-[a(a+2) + b(b+2)-c(c+2)]}{2}},$ we find $$R^{22}_0=A^{-8}=e^{-4\pi i/5}, R^{22}_2=-A^{-4}=-e^{-2\pi i/5}.$$

The Fibonacci $F$-matrix is given by

$$F= \begin{pmatrix} \phi^{-1} & \phi^{-1/2} \\ \phi^{-1/2} & -\phi^{-1} \end{pmatrix}.$$
Either of these quantities can be found using the method outlined in the previous section.

\subsubsection{The Jones representation $\rho_{3,\tau^{\otimes 3}, \tau}$}

Given the $R$ and $F$-symbols, $\rho(\sigma_1)$ and $\rho(\sigma_2)$ take the form

$$\rho(\sigma_1) = \begin{pmatrix} R_1^{\tau\tau} & 0 \\ 0 & R_{\tau}^{\tau\tau} \end{pmatrix}, \hspace{5pt} \rho(\sigma_2)=F\rho(\sigma_1)F^{-1}$$

Explicitly, the generators $\sigma_1$ and $\sigma_2$ have representations
$$\rho(\sigma_1)= \begin{pmatrix} \xi^{-2} & 0 \\ 0 & -\xi^{-1}\end{pmatrix}, \rho(\sigma_2) =  \begin{pmatrix} \phi^{-1}\xi^{2} & -\phi^{-1/2}\xi \\ -\phi^{-1/2}\xi &  - \phi^{-1}\end{pmatrix},$$
where $\xi=e^{2\pi i/5}$ and $\phi$ is the golden ratio.

Having introduced the main tools for calculating braid group representations, we turn to studying their images.


\subsection{The image of the braid group representation}
The basic questions we must have the answers to in order to apply the representations to quantum computation are the following:
\begin{ques}
\begin{enumerate}
\item Is the image $\rho_{k,n,i}(\B_n)$ in $U(V_{k,n,i})$ finite or infinite?
\item If it is infinite, what is the compact Lie group $\overline{\rho_{k,n,i}(\B_n)} \subset U(V_{k,n,i})$?
\end{enumerate}
\end{ques}

The first question was answered by Jones in his seminal 1984 paper \cite{JJ}, and the second by Freedman, Larsen, and Wang in 2002 \cite{FLWu}.
\begin{thm}[Jones] For $k \in \{1,2,4\}$, $\rho_{k,n,i}(\B_n)$ is a finite group. For other values of $k$ and $n\ge 3$, $\rho_{n,k,i}(\B_n)$ is infinite, except for when $k=8$ and $n=4$.
\end{thm}
The following theorem characterizes the closed image of the Jones representation.

\begin{thm}[Freedman, Larsen, W.]
When $\overline{\rho_{k,n,i}(\B_n)}$ is infinite, $SU(V_{k,n,i}) \subset \overline{\rho_{k,n,i}(\B_n)}$.
\end{thm}

Special unitary matrices are the stuff of which quantum gates are made, so this result has an important application to quantum computation.

While we have stated very general results about the images of braid group representations whose proofs are beyond the scope of these notes, there are more elementary ways that we can reproduce these results in the $k=2$ and $k=3$ case, to address whether the Ising and Fibonacci models can be useful for quantum computation.

\subsubsection{The image of $\rho_{2,4,0}$}

\begin{thm}\label{k=2JonesThm}
\begin{enumerate}
\item For $k=2$, the  Temperley-Lieb-Jones algebra $TLJ_n(A)$ is isomorphic to a Clifford algebra.
\item The short exact sequence
$$1 \to \mathds{Z}_2^n \to \rho_{2,n,i}(\B_n) \to S_n \to 1$$ is exact projectively.
\end{enumerate}
\end{thm}
To prove the first part of the theorem, we'll need to know that the Majorana version of a Clifford algebra is the $\mathds{C}$-span of vectors $\{e_1, \ldots, e_{n-1}\}$, subject to the relation $e_ie_j+e_je_i=2\delta_{ij}$.
The $k=2$ Temperley-Lieb-Jones algebra in terms of generators and relations is the $\mathds{C}$-span of the diagrams $\{u_1, \ldots, u_{n-1}\}$, modulo the relation $p_3=0$: $$TLJ_n(A) = \{u_1, \ldots, u_{n-1}\} / p_3=0.$$

In order to show that these two algebras are isomorphic, we need a conversion between the generators $e_i$ of the Clifford algebra and the $u_i$ of the Temperley-Lieb-Jones algebra. Recall $\sigma_i = A +A^{-1} u_i$, and define $g_i=-A^{-1}\sigma_i=-1-A^{-2}u_i$. Since $A$ is an eighth root of unity when $k=2$, $g_i^2=1-du_i$. On the other hand, the $e_i$ can be written as $e_i=(\sqrt{-1})^{i-1}g_i^2g_{i-1}^2\cdots g_1^2$, so that $g_i=\sqrt{-1}e_ie_{i+1}$. Then one can check that their mutual definitions with respect to the $g_i$ agree.

The following proposition contains the relations needed in order to prove the second part of the theorem.
\begin{prop}\label{k=2Identities}
\begin{enumerate}
\item $g_i^2g_{i+1}^2+g_{i+1}^2g_i^2=0$
\item $g_ig_{i\pm1}^2g_i^{-1}=ig_i^2g_{i\pm 1}^2$
\end{enumerate}
\end{prop}
The proof is an exercise for the reader in the graphical calculus. For the second part of the proof, we reduce the problem to that of showing that the image of the pure braid group is finite projectively (up to a scalar which is a root of unity).

Recall the pure braid group $P\B_n$, which is defined implicitly through the short exact sequence
$$1 \to P\B_n \to \B_n \to S_n \to 1,$$
and can be interpreted as a group of braid diagrams which start and end in the same position. For example, the following braid on four strands is a pure braid: $$
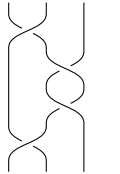
\begin{tikzpicture}[scale=.5]
 \braid[number of strands={4}] (braid)  a_1^{-1}a_2 a_2 a_1^{-1} ;
 \end{tikzpicture}$$  An important class of pure braids are those of the form $\sigma_i^2$ for any generator $\sigma_i$ of $\B_n$.

In order to show $\rho(\B_n)$ is finite, it is enough to show that $\rho(P\B_n)$ is finite. Certain conjugates of elements of the form $\sigma_i^2$ form a generating set $A_{ij}$ of $P\B_n$, where $$A_{ij} = (\sigma_j\sigma_{j-1}\cdots \sigma_{i+1})\sigma_i^2(\sigma_j\sigma_{j-1}\cdots \sigma_{i+1})^{-1},\hspace{10pt} i < j.$$

It is also enough to show that $\rho(P\B_n)$ is finite projectively. Now the first part of Proposition \ref{k=2Identities} tells us that the $g_i^2$'s commute up to an overall minus sign, and furthermore we can deduce that $g_i^{16}=1$. Then it follows from the second part of the proposition that the image of the pure braid group is generated by $g_i^2$.

Thus when $k=2$ the Jones representation of the braid group has finite image.

The physical consequence of \ref{k=2JonesThm} is that
\begin{cor} The Jones representation at level 2 cannot be used for universal quantum computation by braiding alone.
\end{cor}
This means that the set of quantum gates that come from the matrix representations $\rho_{2,n,i}$ of the braiding of the anyons $\{1, \sigma, \psi\}$ of the Ising model is not powerful enough to build a universal quantum computer. In order to prove this corollary, one needs to know about the mathematical formalism of quantum computation, which will be discussed shortly.

One the other hand, we claim that all of $SU(2)$ is contained in $\overline{\rho_{3,\tau^{\otimes 3}, i}(\B_3)}$.

\subsubsection{The image of $\rho_{3,\tau^{\otimes 3}, \tau}$ }\label{imagefib}The following theorem characterizes the closed images of $\rho_{3,\tau^{\otimes 3}, i}$ in $U(V_{3,\tau^{\otimes 3}, i})$.

\begin{thm}\label{1qubit} 
$\overline{\rho_{3,\tau^{\otimes 3}, i}(\B_n)} \supset SU(V_{3,\tau^{\otimes 3}, i})$.
\end{thm}
Recall that the generators $\sigma_1$ and $\sigma_2$ have representations
$$\rho(\sigma_1)= \begin{pmatrix} \xi^{-2} & 0 \\ 0 & -\xi^{-1}\end{pmatrix}, \rho(\sigma_2) =  \begin{pmatrix} \phi^{-1}\xi^{2} & -\phi^{-1/2}\xi \\ -\phi^{-1/2}\xi &  - \phi^{-1}\end{pmatrix}.$$

\subsubsection{Proof of the theorem} We begin by showing that the image is infinite. Of course, it suffices to demonstrate the existence of an element in the image which is of infinite order. Since the matrix representation of $\sigma_1$ is diagonal, it is easy to see that its order is just the least common multiple of the orders of the roots of unity appearing on the diagonal, and hence $\sigma_1$ has order 10. Since all elements of the braid group are in the same conjugacy class, it follows that $\sigma_2$ has order 10 as well.
Let $\sigma_1=\begin{pmatrix} \xi^{-2} & 0 \\ 0 & -\xi^{-1} \end{pmatrix}$ and $\sigma_2=F\sigma_1F$, where $\xi = e^{2\pi i/5}$.
Consider $\sigma_1^m\sigma_2$, where $\sigma_1^{10}=1$ and $m=1,2,\ldots,9$. We claim that when $m=4$ and $m=9$, both elements $\sigma_1^m\sigma_2$ are of infinite order, and moreover, they don't commute. To see that they do not commute, suppose
 $$(\sigma_1^4\sigma_2)(\sigma_1^9\sigma_2)=(\sigma_1^9\sigma_2)(\sigma_1^4\sigma_2)$$
 then using that $\sigma_1$ has order ten, it follows that
 $$\sigma_1^4\sigma_2\sigma_1^{-1}\sigma_2=\sigma_1^{-1}\sigma_2\sigma_1^4\sigma_2$$
and hence $$\sigma_1^5\sigma_2=\sigma_2\sigma_1^5.$$
But one can check using the definitions of $\sigma_1$ and $\sigma_2$ that this is impossible.

To prove that $\sigma_1^4\sigma_2$ and $\sigma_1^9\sigma_2$ are of infinite order, we use the result in \cite{Conway} that all $0$ sum of a few roots of unity are known.  Suppose $\sigma_1^4\sigma_2$ is of finite order, then there will be two roots of unity $\lambda_i, i=1, 2$ such that
$\text{Tr}(\sigma_1^4\sigma_2)=\lambda_1 + \lambda_2$ and $\lambda_1\lambda_2=\det(\sigma_1^4\sigma_2)=\xi$. The trace $\text{Tr}(\sigma_1^4\sigma_2)$ can be written as an integral identity of $\xi$ and $\lambda_1$.  But this cannot happen by inspecting such sums in \cite{Conway}.  Similarly, we can show $\sigma_1^9\sigma_2$ is of infinite order.
Putting $g_1=\sigma_1^4\sigma_2$ and $g_2 =\sigma_1^9\sigma_2$,
this shows that $\overline{\{g_1^m\}}=SO(2)$ and $\overline{\{g_2^m\}}=SO(2)$ both inject into $SU(2)$, showing that this representation is one-qubit universal.  This theorem tells us that $\{\rho(\sigma_1),\rho(\sigma_2)\}$ is a universal gate set, and so in a sense is ``as large as it can get''.

\begin{cor}
Braiding Fibonacci anyons is enough to get any single qubit quantum gate.
\end{cor}
But what about $n$-qubit gates? In general, we need $n$-qubit space $(\mathds{C}^2)^{\otimes n}$ to be contained in $\text{Hom}(3, \tau^{\otimes n}, 1)$. It will be enough so show we can get all two-qubit gates.

In the next section we provide the background necessary to assess the power of the images of the Jones representations as quantum gates. Once the algebra of universality has been established the physical theorems stated for the Ising and Fibonacci theories follow.


\subsection{Quantum gates and universal quantum computation}
Classically, a decision problem is the following: given a sequence of functions $\{f_n\}:\mathds{Z}_2^n$ to $\mathds{Z}_2^n$ on $n$-bit strings, compute $f_n(x)$ for all $x \in \mathds{Z}_2^n$. ``Quantizing" this set up, we have a quantum decision problem - given a sequence of $\{f_n\}$ on $\mathds{C}[\mathds{Z}_2^n] \cong (\mathds{C}^2)^{\otimes n}$, the space of $n$-qubits, find a unitary matrix $U$ such that $U|x\rangle = | f_n(x) \rangle$.
Such matrices are written with respect to the \emph{computational basis} $\mathds{Z}_2^n$ of $(\mathds{C}^2)^{\otimes n}$.

We are always concerned with \emph{efficient approximation} when performing computation. The correct notion of efficiency is that $U$ should be a composition of gates, of length polynomial in $n$.


The building blocks of which such a $U$ is composed are elements of a small \emph{gate set}, say, $S=\{g_1, \ldots, g_m\}$, where each $g_i$ is a $2 \times 2$ or $4 \times 4$ unitary matrix, i.e. each acts on a one qubit ($\mathds{C}^{2}$ ) or two-qubit $(\mathds{C}^{2} \otimes \mathds{C}^{2})$ subspace of $(\mathds{C}^2)^{\otimes n}$. These gate sets, while acting on a few qubits at a time, are extended trivially on the remaining qubits by tensoring with the identity. A gate set is said to be \emph{universal} if we can build any unitary matrix to arbitrary accuracy with a finite number of elements of our gate set. More precisely, if we consider the set of all quantum circuits on $(\mathds{C}^2)^{\otimes n}$ that can be built from our gate set, then it is universal if it is dense in $SU(2^n)$.

\subsubsection{Single-qubit gates}
Some of the most foundational results in quantum computation are the following theorems concerning the gates $$H=\frac{1}{\sqrt{2}} \begin{pmatrix} 1 & 1 \\ 1 & -1\end{pmatrix}, T=\begin{pmatrix}1 & 0 \\ 0 & e^{i \pi/4} \end{pmatrix}, \text{ and } CNOT=\begin{pmatrix} 1 & 0 & 0 & 0 \\ 0 & 1 & 0 & 0 \\ 0 & 0 & 0 & 1 \\ 0 & 0 & 1 & 0 \end{pmatrix}.$$
\begin{thm} The gate set $\{H,T, CNOT\}$ consisting of the Hadamard gate, phase-shift, and controlled-not gates is universal for quantum computation.
\end{thm}

\begin{thm} Let $G$ be the set of all compositions of $H$ and $T$. Then the closure $\overline{G} \supset SU(2)$.
\end{thm}
Therefore, if the Hadamard and $\frac{\pi}{8}$ matrices can either be realized exactly or efficiently approximated by matrices coming from the images of Jones representations, then the corresponding anyon model is sufficient to perform any single-qubit computation. If in addition an {\it entangling} gate like $CNOT$ can be realized, then the anyon model can be used to build a universal quantum computer.

\begin{ques}
Which matrices in $U(2)$ can be realized exactly by a braid up to an overall phase?
\end{ques}

The following theorem is an answer to this question for $\B_3$ in the Fibonacci theory \cite{KBS}.
\begin{thm}
Let $\omega = e^{2\pi i/10}$ and let $u,v \in \mathds{Z}[\omega]$ satisfying $|u|^2 + \frac{|v|^2}{\phi}=1$. Then any matrix of the form

$$M= \begin{pmatrix} u & \bar{v}\phi^{1/2} \\ v\phi^{-1/2} & - \bar{u} \end{pmatrix} \begin{pmatrix} 1 & 0 \\ 0 & \omega^k \end{pmatrix}$$
can be realized exactly by a braid in $\B_3$ in the Fibonacci theory.
\end{thm}


%

\subsubsection{Two-qubit gates and entanglement}
The notion of entanglement is key to understanding universality.
\begin{defn} A gate $g$ in $U(4)$ is not entangling if either $g=A \otimes B$ or SWAP $g = A \otimes B$, where $A,B \in U(2)$ and SWAP $=\begin{pmatrix} 1 & 0 & 0 & 0 \\ 0 & 0 & 1 & 0 \\ 0 & 1 & 0 & 0 \\ 0 & 0 & 0 & 1 \end{pmatrix}$. If a gate is not of this form, then it is called entangling.
\end{defn}

The simplest example of an entangling gate is the CNOT gate. To show that CNOT is entangling, we must prove that neither CNOT nor SWAP CNOT can be written as a tensor product $A \otimes B$, with $A,B \in U(2)$. Recall that the CNOT has the following matrix with respect to the computational basis

$$CNOT = \begin{pmatrix} 1 & 0 & 0 & 0 \\ 0 & 1 & 0 & 0 \\ 0 & 0 & 0 & 1 \\ 0 & 0 & 1 & 0 \end{pmatrix}$$

If a matrix $M$ can be written $M=M_1 \otimes M_2$, where $M_1$ has eigenvalues $\lambda_i$ and $M_2$ has eigenvalues $\mu_j$, then the eigenvalues of $M$ are of the form $\{\lambda_i\mu_j\}$. CNOT has eigenvalues $1,1,1,-1$. So if there were matrices $A$ and $B$ with eigenvalues $\lambda_1, \lambda_2$ and $\mu_1, \mu_2$ respectively, then they would have to satisfy the system of equations

$$\begin{cases} \lambda_1\mu_1=1 \\ \lambda_1\mu_2=1\\ \lambda_2\mu_1=1\\ \lambda_2\mu_2=-1 \end{cases}.$$ But $\det \lambda_{ij} =0$, so this cannot happen. The same argument applies to show that SWAP CNOT cannot be written as a tensor product.

Any four-by-four unitary matrix, that is, any two-qubit quantum gate, can be written as a tensor product of single-qubit gates and an entangling gate, as the following theorem states.

\begin{thm} Given any entangling gate $E$, any matrix in $U(4)$ can be written as a finite product of some number of $E$'s and matrices in $U(2)$ up to an overall phase.
\end{thm}

Thus far we understand how to produce two-qubit gates. To be able to construct $n$-qubit gates, we need the following definition.

\begin{defn} A matrix is $2$-level if it is not the identity only on a $2$-dimensional subspace.
\end{defn}

An important example of a $2$-level matrix for our purposes is of the kind $\begin{pmatrix} \begin{matrix} 1 & 0 \\ 0 & 1 \end{matrix} & 0 \\ 0 & U \\ \end{pmatrix}$, where $U \in U(2)$.

The following lemma collects the facts that allow local unitary computation.
\begin{lem}Every unitary matrix $M$ is a product of $2$-level unitary matrices. Every $2$-level matrix can be realized as a product of $1$-qubit gates and CNOTs.
\end{lem}

The results we have collected thus far imply the following theorem.
\begin{thm} The CNOT gate, together with SU(2) forms a universal gate set. That is, if $U\in SU(2^n)$, then $U$ can be written as a tensor product of CNOT gates and $2 \times 2$ unitary matrices.
\end{thm}
All of the linear algebra is now in place to understand what is needed for a small gate set to be universal.

\subsection{Ising and Fibonacci quantum computers}

Since the image of $\rho_{2,4,0}$ is finite in $U(2)$, we can't get a universal gate set from the Ising theory. On the other hand, we have shown that the closure of the image of $\rho_{3,\tau^{\otimes 3}, \tau}$ contains $SU(2)$, and hence can produce any single-qubit gate. Moreover, the image of $\rho_{3,\tau^{\otimes 6}, 1}:\B_6 \to U(5)$ can be used to approximate an entangling gate, as implied by the theorem below.

In the dense encoding, we choose the two-qubit computation subspace in $V_{3,\tau^{\otimes 6},1}$ as follows by the following fusion tree:  

$$
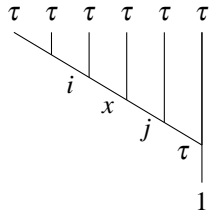
\begin{tikzpicture}
\draw (0,2)--(0,0); 
\draw (-2.5,2)--(0,.5);
\draw (-2, 2)--(-2, 1.7);
\draw (-1.5,2)--(-1.5,1.4);
\draw (-1,2)--(-1,1.1);
\draw (-.5,2)--(-.5,.8);
\draw (0,2)--(0,.5);

\draw (0,-1/4) node {$1$};
\draw (-2.5,2.25) node {$\tau$};
\draw (-2,2.25) node {$\tau$};
\draw (-1.5,2.25) node {$\tau$};
\draw (-1/2,2.25) node {$\tau$};
\draw (-1,2.25) node {$\tau$};
\draw (0,2.25) node {$\tau$};

\draw(-1.75, 1.3) node {$i$};
\draw(-1.25, 1) node {$x$};
\draw(-.75, .7) node {$j$};
\draw(-.25, .4) node {$\tau$};
\end{tikzpicture}$$

$V_{3,\tau^{\otimes 6},1}$ is $5$-dimensional.  When $x=1$, there is only one admissible labeling, so it spans a one-dimensional subspace of $V_{3,\tau^{\otimes 6},1}$.  We will not use this subspace for computation, so it will be called a non-computational subspace.  When $x=\tau$, then all choices of $i,j\in \{1,\tau\}$ are admissible, so we obtain a natural two-qubit subspace of $V_{3,\tau^{\otimes 6},1}$.  We will denote the $4$ basis elements as $\{e_{ij}\}$.  Ideally, we would like to find an entangling braid $b\in \B_6$ on $V_{3,\tau^{\otimes 6},1}$ so that the resulting matrix $\rho(b)$ is the identity on the non-computation subspace, but an entangling gate on the two-qubit subspace.  But we do not know the existence of such braiding gates.  It would be extremely interesting to know if such no leakage entangling braiding gates exist or not.  But in practice, we will use the following density theorem to find entangling braiding gates with arbitrary small leakage to the non-computational subspace.

\begin{thm}\label{2qubits}
\begin{enumerate}
\item $SU(5) \subset \overline{\rho_{3,\tau^{\otimes 6}, 1}(\B_6)}$.
\item Any matrix in $SU(4)$ can be approximated to any precision by the gate set $\{\rho(\sigma_i),i=1,...,5\}$ on the computational subspace $\mathds{C}[\{e_{ij}\}], i,j \in \{1,\tau\}$.
\end{enumerate}
\end{thm}

The proof of this theorem is not elementary so we omit the details.  Presumably we can use an inductive argument using irreducibility and density of one-qubit gates Thm. \ref{1qubit}.  
Therefore a universal gate set can be built from braiding Fibonacci anyons, proving Theorem 1.1 for $r=5$.

\begin{prob} Is there a two-qubit entangling gate that can be realized by braiding exactly in the Fibonacci theory?
\end{prob}

\subsection{General TLJ for quantum computing}
In earlier sections, we explain Ising and Fibonacci theories.  In general any TLJ theory can be used for anyonic quantum computation.

The Jones-Kauffman theory  
at level $k=r-2 \ge 1$ is the TQFT associated to the TLJ theory with the choice of
$$A=\begin{cases} i e^{-2\pi i/4r} & k\equiv 0 \mod 2\\
ie^{2\pi i/4} & k \equiv 1 \mod 2 \text{ and } k \equiv 1 \mod 4 \\
-ie^{-2\pi i /4}  & k \equiv 1 \mod 2 \text{ and } k \equiv -1 \mod 4.
\end{cases}$$

When $k$ is even, the TLJ category is a unitary modular category modeling anyons, while when $k$ is odd, it is a unitary pre-modular category modeling fermions. This generalizes the Ising and Fibonacci theories.

\section{Approximation of The Jones polynomial}\label{jonespolyalg}
In this section we discuss the approximation of the Jones evaluations $J(L,e^{\pm 2\pi i/r})$, i.e. Jones polynomial at roots of unity, by a quantum computer.  This approximation algorithm is a consequence of the efficient simulation of TQFTs with several clarifications \cite{FKW}.  Our approximation goes through the Jones representations of the braid group, therefore we need to choose a closure from braids to links.  For our algorithm, we use the plat-closure of braids with an even number of strands.  It was observed that if instead the braid closure is used, the approximation is potentially easier \cite{SJ}.  Approximations have variations \cite{Ku}, and our approximation is an additive one.  Strictly speaking, we approximate the normalized Jones evaluations: $J(L,e^{\pm 2\pi i/r})$ divided by $d^n$.  For the plat-closures of braids $b \in \B_{2n}$, the unlink of $n$-components has the largest absolute value $d^n$.   It is known that the distributions of Jones evaluations $J(L,e^{\pm 2\pi i/r})$ for $r\ne 1,
2,3,4,6$ are limiting to a Gaussian as $n\rightarrow \infty$ \cite{FLWt}.  Hence, a Jones evaluation is typically small.  Our BQP-complete theorem for an additive approximation with an error scaling as inverse $\textrm{poly}(n,m)$, where $m$ is the length of braids, implies that the normalized Jones evaluation cannot be always exponentially small because otherwise, we could just set the approximation to be $0$.

Recall that the Jones evaluations $J(L,e^{\pm 2\pi i/r})$ is a map from the set of oriented links to  $\mathds{Z}[q^{\pm1/2}]$ for $q=e^{\pm 2\pi i/r}$. In order to turn the evaluation at roots of unity into a computation problem, we must encode a link $L$ and the Jones evaluations $J(L,e^{\pm 2\pi i/r})$ as bit strings, whereupon it becomes Boolean maps $\mathds{Z}_2^n \to \mathds{Z}_2^{m(n)}$, that send the bit strings encoding of the input $L$ to the bit strings encoding of the output $J(L;e^{\pm 2\pi i/r})$.

How does one turn a link into a bit string? First we present $L$ as the plat-closure of some braid, say $L= \widehat{\sigma_{i_k}^{s_k}\cdots \sigma_{i_1}^{s_1}}$. Then the integers $\{i_j\}$ and $\{s_j\}$ can be written in terms of their binary expansions and finally converted into bit strings.

As for encoding $J(L,e^{\pm 2\pi i/r})$ as a bit string, we can use their binary expansions of the real and imaginary parts.  In general these binary expansions will be infinitely long.  In our algorithm we are going to approximate the evaluations.  Therefore, once we are given the error $\epsilon$, we can decide where to truncate the infinite bit strings.

The classical complexity of computing the Jones polynomial exactly at roots of unity is summarized in the following theorem \cite{Ver, FLWt}.

\begin{uthm} For $r\ne 1,2,3,4,6$, computing the Jones evaluations $J(L,e^{2\pi i/r})$ exactly is \#P hard.  Moreover, the Jones evaluations ${\{J(L,e^{2\pi i/r})\}}$ for all link $L$ at $r\ne 1,2,3,4,6$ is dense in $\mathds{C}$.
\end{uthm}

The following table organizes the known complexity results concerning the computation and approximation of the Jones polynomial at roots of unity of order $r\neq 1,2,3,4,6$ \cite{Ver,Ku,FKW,FLWu}.

\begin{figure}[h!]
\centering
\begin{tabular}{l||l|l}
& Exactly & Approximately \\ \hline \hline
Classically & \#P & No FPRAS \\ \hline
Quantum mechanically & ? & BQP-complete \\
\end{tabular}
\end{figure}

\subsection{Approximating Jones evaluations by a quantum computer}
How is the Jones polynomial of a link $L$ evaluated by an anyonic quantum computer at a root of unity? Suppose $b \in \B_{2n}$ is a braid whose closure gives the link $L$. Physically, a ``cup" state is prepared by creating $n$ pairs of anyons from the vacuum. Then the anyons are braided by $b$. Then measurement is performed by projecting onto a ``cap state". This computes $|\langle \text{ cap } | \rho(b) | \text{ cup }\rangle|^2$, which recovers normalized the Jones evaluations. The figure below illustrates the process, which corresponds to the mathematical operation of taking the \emph{plat closure} of $b$.

$$\begin{tikzpicture}
\draw (-1.5,1/2) node {$\langle \text{ cap } | \rho(b) | \text{ cup }\rangle =$};
\draw (0,0) rectangle (3,1);
\draw (3/2,1/2) node {$b$};
\draw (1/4,1) arc (180:0:1/4);
\draw (1,1) arc (180:0:1/4);
\draw (15/8, 9/8) node {$\cdots$};
\draw (9/4,1) arc (180:0:1/4);
\draw (1/4,0) arc (-180:0:1/4);
\draw (1,0) arc (-180:0:1/4);
\draw (15/8, -1/8) node {$\cdots$};
\draw (9/4,0) arc (-180:0:1/4);
\end{tikzpicture}$$

Classically, this computation is hard as the exponential size of $\rho(b)$ hints.  However, there exists an efficient quantum algorithm to approximate the evaluations of the Jones polynomial \cite{FKW}, cf. Theorem \ref{FKWthm}.  Precisely, we have:

\begin{thm} Let $q=e^{\pm 2\pi i /r}$, $d=2\cos \pi/r$,$\hat{\sigma}^P$ the plat closure of $\sigma \in \B_{2n}$, and $J(\hat{\sigma}^P,q): \bigsqcup_{n=1}^{\infty} \B_{2n} \to \mathds{Z}[q^{\pm 1/2}]$ the Jones evaluation. Given $m=|\sigma|$, $n$, there exists a quantum circuit of size polynomial in $n$, $m$ and $1/\epsilon=\textrm{poly}(n,m)$ such that $U_L$ outputs a random variable $Z(\sigma)$, where $0 \le Z(\sigma) \le 1$, and
$$\Biggl | \frac{|J(\hat{\sigma}^P, q)|^2}{d^{n}} - Z(\sigma)\Biggr |  < \epsilon.$$
\end{thm}
Such an approximation is called an \lq\lq additive'' scheme.  The full details of the proof are not provided here and can be found in \cite{cbms,AJZ}.  We will simply illustrate the main steps and ideas.

\subsection{Turn basis vectors into bit strings}
For an illustration of this step, that of converting basis vectors into bit strings, we consider the case $k=3$ and the six-strand braid group $\B_6$. Given the fusion tree in $\text{Hom}(0,1^{\otimes 6})$, we attach a qubit at each vertex with basis $|i_1 i_2 i_3 i_4 \rangle$,

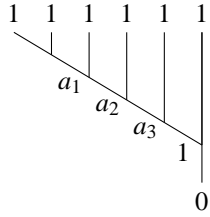
$$\begin{tikzpicture}
\draw (0,2)--(0,0); 
\draw (-2.5,2)--(0,.5);
\draw (-2, 2)--(-2, 1.7);
\draw (-1.5,2)--(-1.5,1.4);
\draw (-1,2)--(-1,1.1);
\draw (-.5,2)--(-.5,.8);
\draw (0,2)--(0,.5);

\draw (0,-1/4) node {$0$};
\draw (-2.5,2.25) node {$1$};
\draw (-2,2.25) node {$1$};
\draw (-1.5,2.25) node {$1$};
\draw (-1/2,2.25) node {$1$};
\draw (-1,2.25) node {$1$};
\draw (0,2.25) node {$1$};

\draw(-1.75, 1.3) node {$a_1$};
\draw(-1.25, 1) node {$a_2$};
\draw(-.75, .7) node {$a_3$};
\draw(-.25, .4) node {$1$};
\end{tikzpicture}$$

where $i_j \in \{0,1\}$ and $a_j \in \{0,2\}$.
If $a_j=0$, we take $i_j=0$, and if $a_j=2$, then we take $i_j=1$. Then we can define the map $a_1a_2a_3 \mapsto |i_1 i_2 i_3 i_4 \rangle$.
This gives an efficient embedding of a basis of $V_{3, 6, 0}$ into bit strings.

\subsection{Simulating the Jones representation}

When we use anyons for quantum computation, we choose a computational subspace $(\mathds{C}^2)^{\otimes l}$ in $V_{k,1^{\otimes m},0}$.  Now for the simulation of the Jones representation $V_{k,m,0}$, we seek a quantum circuit $U$ which makes the following diagram commute.  By turning basis into bit strings, $V_{k,m,0}$ is embedded as a subspace in $(\mathds{C}^2)^{\otimes (m-2)}$.  The Jones representation of $\B_{m}$ is extended to $(\mathds{C}^2)^{\otimes (m-2)}$ by the identity on the orthonormal complement of embedded $V_{k,m,0}$.

$$\xymatrix{(\mathds{C}^2)^{\otimes l}  \ar@{^{(}->}[r]&V_{k,1^{\otimes m},0}\ar@{^{(}->}[r]\ar[d]^{\rho(\sigma)}& (\mathds{C}^2)^{\otimes (m-2)}\ar[d]^{U} \\ &V_{k,1^{\otimes m},0}\ar@{^{(}->}[r]&(\mathds{C}^2)^{\otimes (m-2)}}$$

We must compute the braid group action on the basis of $V_{k,m,0}$.  By thinking about how the braid group generators act on $V_{k,m,0}$, we are led to pieces of fusion trees like the ones below, those diagrams now being drawn horizontally.

$$
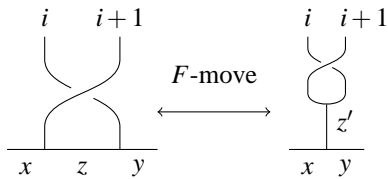
\begin{tikzpicture}
\node (first) {};
\draw (.5,-1.5)--(2.5,-1.5) {};
\braid[number of strands={2}] a_1^-1;
\draw (.75,-1.75) node {$x$};
\draw (2.25,-1.75) node {$y$};
\draw (1.5, -1.75) node {$z$};
\draw (1, .25) node {$i$};
\draw (2, .25) node {$i+1$};
\draw (3.25,-.5) node {$F$-move};
\path[<->] (2.5,-1) edge (4,-1);
\node (second) [right=of first] {};
\begin{scope}[shift={(second)},scale=.5, xshift=5.5cm]
\draw (.5,-3)--(2.5,-3) {};
\braid[number of strands={2}] a_1^-1;
\draw (1,-1.5) to [out=-90, in=-90] (2,-1.5);
\draw (1.5,-1.8)--(1.5,-3);
\draw (1,-3.5) node {$x$};
\draw (1, .5) node {$i$};
\draw (2.5, .5) node {$i+1$};
\draw (2,-3.5) node {$y$};
\draw (2, -2.25) node {$z'$};
\end{scope}
\end{tikzpicture}$$

The diagram on the left can be changed to the form on the right via an $F$-move, and then the braiding can be removed using an $R$-symbol. In this manner, by stacking the braid group generators on a basis element of $V_{k,1^{\otimes m},0}$ and then using the graphical calculus to resolve it into a linear combination of our computational basis elements, the Jones representation can be calculated for any $n$.  The important observation is that the extended Jones representation is now local: the $2$-level matrix in the definition of the Jones representation acts now on only two qubits.  In physical terms, since everything is localized, we only need to concern ourselves with two qubits at a time.

Finally, due to the fusion rules for basis elements of $V_{k,1^{\otimes m},0}$, the Jones representation is a composition of a sequence of multi-qubit controlled $2$-qubit gates on $(\mathds{C}^2)^{\otimes (m-2)}$.  Such controlled gates can be implemented efficiently on a quantum computer using a few ancillary qubits.  Therefore, we can approximate the Jones evaluations efficiently by a quantum computer.

\section{Localization of braid group representations}\label{localization}
The quantum circuit model is explicitly local: the $n$-qubit spaces are tensor powers $(\C^{2})^{\ot n}$ and the $n$-qubit circuits are composed of gates (e.g. promotions of SWAP, CNOT) that act non-trivially on just a few adjacent qubit spaces.

In contrast, the topological model relies upon gates that are not explicitly local, coming from representations of the braid group.
So far we have met several families of $\B_n$ representations: the local representations $\rho^R$ associated with an $R$-matrix, the Burau representations $\rho$ and their reduced versions $\tilde{\rho}$ and the Jones representations (generic and specialized).  In subsection \ref{imagefib} we gave a detailed version of Theorem \ref{FLWthm}, which is the main result of \cite{FLWu}.  Consequently, the quantum circuit model (hence a quantum Turing machine) can be efficiently simulated on a topological quantum computer via certain level $k$ Jones representations of the braid group.

The main theorem (i.e. Theorem \ref{FKWthm}) of \cite{FKW} is a partial converse: the specialized Jones polynomial can be efficiently approximated on the quantum circuit model (cf. Section \ref{jonespolyalg}).  This is achieved by exploiting a \emph{hidden locality} in topological quantum field theory, which we now outline.  Recall from Subsection \ref{TLJmatrix} that the labels for the Temperley-Lieb-Jones category at level $k$ are $\mathcal{L}:=\{0,\ldots,k\}$.  The $\B_n$ representation obtained from the standard faithful $TL_n(A)$-module is $\cH:=\bigoplus_j \Hom(j,1^{\ot n})$.  Here each direct summand $\cH_j=\Hom(j,1^{\ot n})$ is an irreducible $\B_n$ representation associated with a disk with $n$ interior points marked with the anyon type $1$ and the boundary labelled $j$. Now we decompose our $n$-punctured disks into $n-1$ pairs of pants by making $n-2$ concentric circular cuts for each boundary label $j$.  The gluing and disjoint union axioms then show that
$$\cH=\bigoplus_{(i_1,\ldots,i_{n-1})\in\mathcal{L}^{n-1}}\Hom(1^{\ot 2},{i_1})\ot\Hom(1\ot {i_1},{i_2})\ot\cdots\ot\Hom(1\ot {i_{n-2}},{i_{n-1}}).$$
Here the boundary label is $j=i_{n-1}$.  Now we set $V=\bigoplus_{(a,b,c)}\Hom(a\ot b,c)$ and distribute $\ot$ over $\oplus$ to realize $\cH$ inside $V^{\ot (n-1)}$. The complement $\cH^{\perp}$ of $\cH$ inside $V^{\ot (n-1)}$ is does not typically admit a braid group action. Alternatively we can be slightly more efficient and take $U=\bigoplus_{(b,c)}\Hom(1\ot b,c)$.
The upshot is that the specialized Jones representations of $\B_n$ can be realized inside a vector space of the form $V^{\ot f(n)}$, but with $\B_n$ only acting on a certain \emph{hidden} subspace. One may employ the same technique for the Fibonacci representations.
\begin{ex}
 Set $k=2$ and show that $\dim(V)=10$, while $\dim(U)=4$.  Observe that for the $\B_3$ representation we have $\dim(V_{3,1})=2$ which is embedded into either $V^{\ot 2}$ or $U^{\ot 2}$ which has a very large complement.
\end{ex}

This gross inefficiency motivates the following:
\begin{ques}
 When can a family of $\B_n$ representations be realized locally (uniformly for all $n$, and ``on the nose'')?
\end{ques}
Eventually we will restrict to unitary representations, but
first we must make sense of what sort of families we are interested in.
\subsection{Sequences of $\B_n$ Representations}
Notice that we have natural injective group homomorphisms $\iota:\B_{n}\rightarrow \B_{n+1}$ given by $\iota(\sigma_i)=\sigma_i$, for $1\leq i\leq n-1$ allowing us to identify $\B_n$ as a subgroup of $\B_{n+1}$\footnote{Of course there are many less natural injective homomorphisms, for example $\sigma_i\mapsto (\sigma_{n-i+1})^{-1}$ can be verified as an injective homomorphism.}.  Which families of representations respect these identifications?  For precision's sake we phrase the following in terms of group algebras \cite{rowellwang2010}:
\begin{definition}\label{seq}
 An indexed family of complex $\B_n$-representations $(\rho_n,V_n)$ is a \emph{sequence of braid representations} if there exist injective algebra homomorphisms $\varphi_n:\C\rho_n(\B_n)\rightarrow \C\rho_{n+1}(\B_{n+1})$ such that the following diagram commutes:

$$\xymatrix{ \C\B_n\ar[r]^{\rho_n}\ar@{^{(}->}[d]^\iota &
\C\rho_n(\B_n)\ar@{^{(}->}[d]^{\varphi_n} \\ \C\B_{n+1}\ar[r]^{\rho_{n+1}} &
\C\rho_{n+1}(\B_{n+1})}$$
\end{definition}

\begin{example}
 If $R$ is a solution to the Yang-Baxter equation on a vector space $V$, then it is easy to see that $\rho_n:\C\B_n\rightarrow \End(V^{\ot n})$ given by $\rho_n(\sigma_i)=I_V^{\ot i-1}\ot R\ot I_V^{\ot n-i-1}$ is a sequence of braid representations: take $\varphi_n:\End(V^{\ot n})\rightarrow\End(V^{\ot n+1})$ to be $\varphi_n(f)=f\ot I_V$.
 \end{example}

The Jones representations (specialized or not) and the (related) Fibonacci representations $\rho_{n,\tau}$ are sequences in this sense.  For example, for the generic Jones representation we have $\C\rho_n(\B_n)=TL_n(A)$ and $\C\rho_{n+1}(\B_{n+1})=TL_{n+1}(A)$.  Letting $\varphi_n(u_i)=u_i$ we see that the appropriate diagram commutes.

\begin{example}
 The Burau (reduced or unreduced) representations do not form a sequence of $\B_n$ representations in our sense.  Indeed in the case where $\tilde{\rho}$ is irreducible, we have $\C\tilde{\rho}(\B_n)\cong M_{n-1}(\C)$ (for all $n$).  Since there are no injective homomorphisms from $M_{n-1}(\C)$ to $M_n(\C)$, the required map $\varphi_n$ does not exist.
\end{example}

\begin{ex}
 Show that the standard permutation representations of $S_n$, lifted to $\B_n$ in the obvious way via $(i\;i+1)\mapsto \sigma_i$ is not a sequence in our sense.
\end{ex}

Now we can describe what we mean by a \emph{localization} of a sequence of braid group representations.

\begin{definition}\label{localdef}
Suppose $(\rho_n,V_n)$ is a sequence of braid representations.  A
\emph{localization} of $(\rho_n,V_n)$ is a braided vector space $(W,R)$ with
$R\in\U(W^{\ot 2})$ such that for all $n\geq 2$ there exist injective algebra
homomorphisms $\psi_n:\C\rho(\B_n)\rightarrow \End(W^{\ot n})$ satisfying
$\psi_n\circ\rho(b)=\rho_R(b)$ for $b\in\B_n$.

\end{definition}

This definition may seem a bit complicated at first, but encapsulates the notion of ``on the nose'' local realizations of a sequence of $\B_n$ representations.  From the point of view of quantum computation, we are trying to discover when the singleton gate set $\{R\}$ can simulate all braiding gates.  In spite of the slightly mystifying definition, the idea is quite simple: we want to find a single solution to the Yang-Baxter equation $R$ on a vector space $W$ so that

\begin{enumerate}
\item For each $n$, $(\rho_n,V_n)$ is a sub-representation of $(\rho_R,W^{\ot n})$.  Notice that we distinguish between equivalent irreducible sub-representations of $V_n$: if $\ell$ isomorphic copies of some fixed irreducible $U$ appears in $V_n$ then $W^{\ot n}$ must contain at least $\ell$ copies of $U$.  This is a distinction at the level of algebras: $\C^2$ is not a faithful representation of $M_2(\C)\oplus M_2(\C)$, but $\C^2\oplus \C^2$ is.

\item There are no irreducible $\B_n$-subrepresentations of $W^{\ot n}$ that do not appear in $V_n$.  Whereas the hidden locality of \cite{FKW} has a large non-computational space upon which the braid group does not act, we are asking that there is no such axillary space.
\end{enumerate}

Comparing with the property $\mathbf{F}$ conjecture, we obtain:
\begin{conj}\label{uybconj}

Suppose $(V,R)$ is a unitary solution to the YBE such that $R$ has finite (projective) order, with corresponding $\B_n$-representations $(\rho_R,V^{\ot n})$.  Then $\rho_R(\B_n)$ is a finite group (projectively).
\end{conj}
 If the words \emph{unitary} or \emph{finite order} are omitted Conjecture \ref{uybconj}(a) is false, see \cite{rowellwang2010} for examples.

\subsection{Non-localizable representations}

In the braided fusion category setting Conjecture \ref{uybconj} is closely related to another fairly recent conjecture (see \cite[Conjecture 6.6]{RSW}).  Braided fusion categories are naturally divided into two classes according to the algebraic complexity of their fusion rules.  In detail, one defines the \emph{Frobenius-Perron dimension} $\FPdim(X)$ of an object $X$ in a fusion category $\CC$ to be the largest eigenvalue of the fusion matrix $N_X$ corresponding to tensoring with $X$ on the left.    If $\FPdim(X_i)\in\N$ for all simple $X_i$ then one says $\CC$ is \emph{integral} while if $\FPdim(X_i)^2\in\N$ for all simple $X_i$ then $\CC$ is said to be \emph{weakly integral}.  An object $X$ in a braided fusion category $\CC$ is said to have \emph{property \textbf{F}} if the $\B_n$-representations on $\End(X^{\ot n})$ have finite image for all $n$.  Then a version of Conjecture 6.6 of \cite{RSW} states: \emph{an object $X$ has property $\textbf{F}$ if, and only if, $\FPdim(X)^2\in\N$}.  Some recent
progress towards this conjecture can be found in \cite{propF,RUMA} and further evidence can be found in \cite{LRW,LR1}.  Combining with Conjecture \ref{uybconj} we make the following:

\begin{conj}\label{localintegral}
Let $X$ be a simple object in a braided fusion category $\CC$. The representations $(\rho_X,\End(X^{\ot n}))$ are localizable if, and only if, $\FPdim(X)^2\in\N$.
\end{conj}

The main result of this section is to indicate that the specialized Jones representations are not localizable unless $r=1,2,3,4$ or $6$.  That is, the sequence of representations coming from the Temperley-Lieb algebras at all other roots of unity are not localizable.  In the next section we will show the converse for $r=4$ and $r=6$, with the other cases being trivial since the corresponding categories are pointed (and hence the representations are $1$-dimensional).

To any sequence of multi-matrix algebras $\mathcal{S}:=\C=A_1\subset\cdots\subset A_j\subset A_{j+1}\subset\cdots$ with the same identity one associates the \emph{Bratteli diagram} which encodes the combinatorial structure of the inclusions.  The Bratteli diagram for a pair $M\subset N$ of multi-matrix algebras is a bipartite digraph $\Gamma$ encoding the decomposition of the simple $N$-modules into simple $M$-modules, and the inclusion matrix $G$ is the adjacency matrix of $\Gamma$.  More precisely, if $N\cong\bigoplus_{j=1}^t \End(V_j)$ and $M\cong\bigoplus_{i=1}^s \End(W_i)$ the inclusion matrix $G$ is an $s\times t$ integer matrix with entries:
$$G_{i,j}=\dim\Hom_M(\Res_M^NV_j,W_i)$$ i.e. the multiplicity of $W_i$ in the restriction of $V_j$ to $M$.
For an example, denote by $M_n(\C)$ the $n\times n$ matrices over $\C$ and let $N=M_4(\C)\oplus M_2(\C)$ and $M\cong\C\oplus\C\oplus M_2(\C)$ embedded in $N$ as matrices of the form:
$$\begin{pmatrix} a & 0 & 0\\ 0 & a&0\\ 0 & 0 & A\end{pmatrix}\oplus \begin{pmatrix} a & 0\\ 0 & b\end{pmatrix}$$
where $a,b\in\C$ and $A\in M_2(\C)$.  Let $V_1$ and $V_2$ be the simple $4$- and $2$-dimensional $N$-modules respectively, and $W_1$, $W_2$ and $W_3$ be the simple $M$-modules of dimension $1$, $1$ and $2$.  Then the Bratteli diagram and corresponding inclusion matrix for $M\subset N$ are:

$$
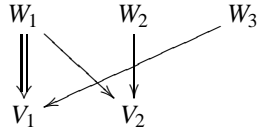
\xymatrix{ W_1\ar@{=>}[d]\ar[dr] & W_2\ar[d] & W_3\ar[dll]\\ V_1 & V_2}$$
and
$$\begin{pmatrix} 2& 1\\ 0 & 1\\ 1 &0\end{pmatrix}.$$

The Bratteli diagram for the sequence $\mathcal{S}$ is the concatenation of the Bratteli diagrams for each pair $(A_k,A_{k+1})$, with corresponding inclusion matrix $G_k$.  We organize this graph into levels (or stories) corresponding to each algebra $A_k$ so that the Bratteli diagram $(A_{k-1},A_k)$ is placed above the vertices labelled by simple $A_k$-modules, and that of $(A_k,A_{k+1})$ is placed below.   Having fixed an order on the simple $A_k$-modules we record the corresponding dimensions in a vector $\mathbf{d}_k$.  Observe that $\mathbf{d}_{k+1}=G_k^T\mathbf{d}_k$.

Let us illustrate this for the Fibonacci theory corresponding to the colored TLJ-theory at $A=ie^{2\pi i/20}$.  Here we have two labels $\one$ and $\tau$ as in subsection \ref{fibthy}.  Decomposing the simple $TLJ_n(A)$ modules for this theory as $TLJ_{n-1}(A)$ modules for $n=1,2,\ldots$ we have:

$$
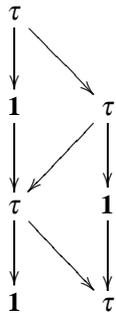
\xymatrix{\tau\ar[d]\ar[dr]\\ \one\ar[d] & \tau\ar[dl]\ar[d]\\ \tau\ar[d]\ar[dr]&\one\ar[d]\\ \one&\tau}$$
For each $n>1$, $TLJ_n(A)$ has two simple modules: $V_{\one,n}:=\Hom(\one,\tau^{\ot n})$ and $V_{\tau,n}:=\Hom(\tau,\tau^{\ot n})$. Moreover, $\C\rho_{n}(\B_n)=TLJ_n(A)$ so these are irreducible $\B_n$-representations.
Let us compute the inclusion matrices as above, ordering the modules $[V_{\one,n},V_{\tau,n}]$ in spite of the alternating arrangement of the Bratteli diagram. Since $V_{\one,n}|_{\B_{n-1}}\cong V_{\tau,n-1}$ and $V_{\tau,n}|_{\B_{n-1}}\cong V_{\one,n-1}\oplus V_{\tau,n-1}$ we have: are $G=G_n=\begin{pmatrix}0 & 1\\ 1& 1\end{pmatrix}$ for all $n>1$.
By computing powers of $G$ one can see that $\dim(V_{\one,n})=f_{n-1}$ and $\dim(V_{\tau,n})=f_n$ where $f_0=1,f_1=1,f_2=1,f_3=2,\ldots$ is the Fibonacci sequence.

Now suppose that we could find a Yang-Baxter matrix $R$ on a space $W$ of dimension $d$ localizing the sequence of $\B_n$-representations $(\rho_n,V_{\one,n}\oplus V_{\tau,n})$.  Using the algebra injections $\psi_n$, the space $W^{\ot n}$ becomes a $TLJ_n(A)$-module and hence $W^{\ot n}\cong a_nV_{\one,n}\oplus b_nV_{\tau,n}$ as $TLJ_n(A)$ (or $\B_n$) modules with $a_n\geq 1$, $b_n\geq 1$ multiplicities.  Notice that $d^n=a_nf_{n-1}+b_nf_n$ for all $n>1$.  We can use $G$ to inductively express the multiplicities $(a_n,b_n)$. Indeed, since restricting $a_nV_{\one,n}\oplus b_nV_{\tau,n}$ to $\B_{n-1}$ we get $b_nV_{\one,n-1}\oplus (a_n+b_n)V_{\tau,n-1}$, we have $G(a_n,b_n)^T=(a_{n-1},b_{n-1})$.  Notice also that $G(f_{n-2},f_{n-1})^T=(f_{n-1},f_n)$.  Thus the formula $d^n=a_nf_{n-1}+b_nf_n$ valid for all $n>1$ gives us the two equations: $\langle(a_n,b_n),G(f_{n-2},f_{n-1})\rangle=d^n$ and $\langle G(a_n,b_n),(f_{n-1},f_n)\rangle=d^{n-1}$.  But since $G^T=G$ we have $$d^n=\langle(a_n,b_n),G(f_{n-2},f_{n-1})\
rangle=\langle G(a_n,b_n),(f_{n-1},f_n)\rangle=d^{n-1},$$ a contradiction.  So the Fibonacci theory cannot be localized.

This may seem a bit \emph{ad hoc}, but in fact this can be generalized whenever the Bratteli diagram for $\C\rho_n(\B_n)\subset\C\rho_{n+1}(\B_{n+1})$ is periodic of some period $k$.  In this case there is a strictly positive integer-valued square matrix $G$ that describes the inclusion of $\C\rho_n(\B_n)\subset\C\rho_{n+k}(\B_{n+k})$.  One then applies the the Perron-Frobenius theorem to see that some vector of multiplicities $\mathbf{b}_n$ is an eigenvector of $G$ corresponding to the largest eigenvalue $\lambda$ of $G$.  This implies that $\lambda\in\Z$, since $G$ and $\mathbf{b}_n$ are integral, which often leads to a contradiction (see \cite{rowellwang2010} for details).

\subsection{Jones representation at levels $2$ and $4$}
In this section we give explicit localizations for the Jones representations at levels $2$ and $4$.

For the Ising theory (level $2$) an explicit localization appears in \cite{FRW}.    The objects are $\one,\sigma$ and $\psi$ where $\FPdim(\sigma)=\sqrt{2}$ and $\FPdim(\psi)=1$.    The Bratteli diagram is:
$$
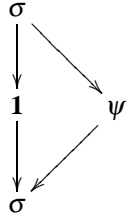
\xymatrix{\sigma\ar[d]\ar[dr]\\ \one\ar[d] & \psi\ar[dl]\\ \sigma}$$ and the matrix $$\frac{-e^{-\pi i/4}}{\sqrt{2}}\begin{pmatrix}1 & 0 & 0 & 1\\0 & 1& -1 & 0\\0 & 1&1& 0\\-1 & 0 & 0
&1\end{pmatrix}$$ gives an explicit localization (see \cite[Section 5]{FRW}).

At level $4$ ($A=ie^{-2\pi i/24}$), the categorical model is a rank $5$ category with simple objects $\one,Z$ of dimension $1$, $Y$ of dimension $2$ and $X,X^\prime$ of dimension $\sqrt{3}$.  The fusion rules for this category are determined by:
\begin{eqnarray}
&& X\ot X\cong \one\oplus Y,\quad X\ot X^\prime\cong Z\oplus Y\\
&& X\ot Y\cong X\oplus X^\prime, \quad Z\ot X\cong X^\prime.
\end{eqnarray}
The Bratteli diagram (starting at level 1) is shown in below.$$
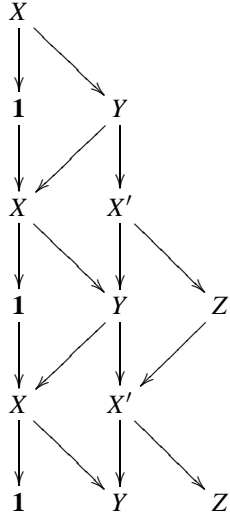
\xymatrix{ X\ar[d]\ar[dr] & \\ \one\ar[d] & Y\ar[dl]\ar[d]\\X\ar[d]\ar[dr] & X^\prime\ar[d]\ar[dr]\\ \one\ar[d] & Y\ar[dl]\ar[d] & Z\ar[dl]\\X\ar[d]\ar[dr] & X^\prime\ar[d]\ar[dr]\\\one & Y & Z }$$
We have $TLJ_n(A)\cong\End(X^{\ot n})$ for each $n$, where the isomorphism is induced by
$$g_i\leftrightarrow Id_X^{\ot i-1}\ot c_{X,X}\ot Id_X^{\ot n-i-1}\in\End(X^{\ot n}).$$
Here $c_{X,X}$ is the (categorical) braiding on the object $X$.
The irreducible sectors of ${TLJ}_n(A)$ under this isomorphism are the $\End(X^{\ot n})$-modules $H_{n,W}:=\Hom(W,X^{\ot n})$ where $W$ is one of the $5$ simple objects in $\CC$.  Observe that for $n$ even $W$ must be one of $\one,Y$ or $Z$ while for $n$ odd $W$ is either $X$ or $X^\prime$.  We have the following formulae for the dimensions of these irreducible representations (for $n$ odd):

$$
\dim\Hom(X,X^{\ot n})=\frac{3^{\frac{n-1}{2}}+1}{2},
\quad \dim\Hom(X^\prime,X^{\ot n})=\frac{3^{\frac{n-1}{2}}-1}{2},$$
$$
\dim\Hom(\one,X^{\ot n+1})=\frac{3^{\frac{n-1}{2}}+1}{2}, \quad \dim\Hom(Y,X^{\ot n+1})=3^{\frac{n-1}{2}},
$$$$ \dim\Hom(Z,X^{\ot n+1})=\frac{3^{\frac{n-1}{2}}-1}{2}
$$

We present the explicit localization, referring the reader to \cite{rowellwang2010} for a complete proof.  Here $\omega=e^{2\pi i/3}$ is a 3rd root of unity.

$R=\frac{i}{\sqrt{3}}\begin{pmatrix}\om&0&0&0&1&0&0&0&\om\\ 0&\om&0&0
&0&\om&1&0&0\\ 0&0&\om&
\om^2&0&0&0&\om^2&0
\\ 0&0&\om^2&\om&0&0&0&\om^2&0
\\ \om&0&0&0&
\om&0&0&0&1\\ 0&1&0
&0&0&\om&\om&0&0\\ 0&
\om&0&0&0&1&\om&0&0
\\ 0&0&\om^2&\om^2&0&0&0&\om&0\\ 1&0
&0&0&\om&0&0&0&\om\end {pmatrix}$

 In fact we have a complete characterization of localizable (uncolored) Jones representations, verifying Conjecture \ref{localintegral} in these cases:
\begin{thm}\label{cor:nolocal}
 The Jones representation at level $k$ can be localized if and only if $k\in\{1,2,4\}$.
\end{thm}

This should be compared with Theorem \ref{FLWthm} from which it follows that the Jones representations are universal for quantum computation precisely when the representations are not localizable.


\begin{thebibliography}{9}

\bibitem{AJZ}D.\ Aharonov, V.\ Jones, and Z.\ Landau. \emph{A polynomial quantum algorithm for approximating the Jones polynomial.} Algorithmica 55.3 (2009): 395-421.

\bibitem{bigelow}
S. Bigelow,
The Burau representation is not faithful for $n = 5$,
\url{http://arxiv.org/pdf/math/9904100v2.pdf}
1999.

\bibitem{Conway}J.\ Conway, and A.\ Jones. \emph{Trigonometric Diophantine equations (On vanishing sums of roots of unity).} Acta Arithmetica 30.3 (1976): 229-240.

\bibitem[E]{Epple} M.\ Epple,  \emph{Orbits of asteroids, a braid, and the first link invariant.} Math. Intelligencer \textbf{20} (1998), no. 1, 45–-52.
\bibitem{FRW}  J.\ Franko; E.\ C.\ Rowell; Z.\ Wang, Extraspecial 2-groups
and images of braid group
representations, \emph{J. Knot Theory Ramifications} \textbf{15} (2006) no. 4,
1--15.

\bibitem{FKW}M.\ H.\ Freedman; A.\ Kitaev; Z.\ Wang,
{Simulation of topological field theories by quantum
computers.} \emph{Comm. Math. Phys.} \textbf{227} (2002) no. 3, 587--603.

\bibitem{FLWu} M.\ H.\ Freedman; M.\ J.\ Larsen; Z.\ Wang,
    {A modular functor which is universal for quantum computation.} \emph{Comm.
Math. Phys.} \textbf{227} (2002) no. 3, 605--622.

\bibitem{FLWt} M.\ H.\ Freedman; M.\ J.\ Larsen; Z.\ Wang,\emph{The Two-Eigenvalue Problem and Density of Jones Representation of Braid Groups.}
 Communications in mathematical physics 228.1 (2002): 177-199.

\bibitem{fultonharris}
W. Fulton and J. Harris,
Representation Theory, a First Course. Springer, New York,
1991.

\bibitem{JJ}V.F.R.\ Jones, \emph{Braid groups, Hecke algebras and type $II_1$ factors.} Geometric methods in operator algebras (Kyoto, 1983) 123 (1983): 242-273.

\bibitem{jones}
V.F.R. Jones,
Hecke-algebra representations of braid groups and link polynomials.
\emph{Ann. Math.} 126 (1987) 335-288

\bibitem{KBS}V.\ Kliuchnikov, A.\ Bocharov, K.\ M.\ Svore. \emph{Asymptotically optimal topological quantum compiling.} Physical review letters 112.14 (2014): 140504.

\bibitem{Ku}G.\ Kuperberg. \emph{How hard is it to approximate the Jones polynomial?.} arXiv preprint arXiv:0908.0512 (2009).


\bibitem{LR1}  M.\ J.\ Larsen, E.\ C.\ Rowell, {An algebra-level
version of a link-polynomial identity of Lickorish}. \emph{Math.
Proc. Cambridge Philos. Soc.} \textbf{144} no. 3 (2008), 623--638.

\bibitem{LRW}  M.\ J.\ Larsen, E.\ C.\ Rowell, and Z.\ Wang:
    {The $N$-eigenvalue problem and two applications.}
    \textit{Int.\ Math.\ Res.\ Not.} \textbf{2005} (2005) no.\ 64, 3987--4018.

\bibitem{morrison}
S. Morrison,
A formula for the Jones-Wenzl projectors
\url{http://arxiv.org/pdf/1503.00384v1.pdf}
\bibitem{propF} D.\ Naidu; E.\ C.\ Rowell, \emph{A finiteness property
for braided fusion categories}, \emph{Algebr. Represent. Theory}. \textbf{15} (2011) no. 5, 837--855.
    \bibitem{RUMA} E.\ C.\ Rowell, {Braid representations from quantum
groups of exceptional Lie type},
\emph{Rev. Un. Mat. Argentina}
\textbf{51} (2010) no. 1, 165--175.
\bibitem{rowellwang2010}
  E. C. Rowell, Z. Wang ,
  Localization of unitary braid representations,
\emph{Comm. Math. Phys.} \textbf{311} (2012) no. 3, 595--615.
  \url{http://arxiv.org/pdf/1009.0241v2.pdf},
  2010.
\bibitem{RSW} E.\ Rowell; R.\ Stong; Z.\ Wang, {On classification of
modular tensor categories},  \emph{Comm. Math. Phys.}
\textbf{292}  (2009),  no. 2,
343--389.

\bibitem{SJ}P.\ W.\ Shor; S.\ P.\ Jordan, \emph{Estimating Jones polynomials is a complete problem for one clean qubit.}  Quantum Information and Computation 8.8 (2008): 681-714.

\bibitem{Ver}D.\ L.\ Vertigan, On the computational complexity of Tutte, Jones, Homfly and Kauffman invariants. Diss. University of Oxford, 1991.

\bibitem{cbms}
Z. Wang,
Topological Quantum Computation, Amer. Math. Soc., Providence.
\url{http://www.math.ucsb.edu/~zhenghwa/data/course/cbms.pdf}
2010.


\end{thebibliography}
\end{document}